\documentclass[12pt]{article}

\usepackage{amsmath}

\usepackage{color}
\usepackage{epsfig}
\usepackage{graphicx}
\usepackage{amsbsy}
\usepackage{latexsym}
\usepackage[mathscr]{eucal}
\usepackage{amsfonts,amsmath}
\usepackage{amssymb}
\usepackage{stmaryrd}
\usepackage{hyperref}
\usepackage{algorithmicx}
\usepackage{algorithm}
\usepackage{algpseudocode}

\usepackage{hyperref}
\usepackage[latin1]{inputenc}
\usepackage{fullpage}

\begin{document}
\title{Compositional Sparsity, Approximation Classes, and Parametric Transport Equations\thanks{This work has been supported in
part by National Science Foundation grant  NSF-DMS-2012469, and the  SFB 1481,
funded by the German Research Foundation
}\\
{\small \it{Dedicated to Ronald DeVore on the occasion of his 80th birthday}}}
%

\author{Wolfgang Dahmen\thanks{Mathematics Department, University of South Carolina, Columbia, SC, USA}
} %

\textheight=213 mm
\textwidth=220 mm
\renewcommand{\arraystretch}{1.5}
\hbadness=10000
\vbadness=10000
\newtheorem{lemma}{Lemma}[section]
\newtheorem{prop}{Proposition}[section]
\newtheorem{cor}[prop]{Corollary}
\newtheorem{theorem}[prop]{Theorem}
\newtheorem{rem}[prop]{Remark}
\newtheorem{example}[prop]{Example}
\newtheorem{definition}[prop]{Definition}
\def\ds{\displaystyle}
\def\ev#1{\vec{#1}}     
\newcommand{\lt}{\ell_{2}(\cJ)}
\def\Supp#1{{\rm supp\,}{#1}}
\def\e{\varepsilon}
\def\nl{\newline}
\def\T{{\relax\ifmmode I\!\!\hspace{-1pt}T\else$I\!\!\hspace{-1pt}T$\fi}}
\def\N{\mathbb{N}}
\def\R{\mathbb{R}}
\def\Z{{\relax\ifmmode Z\!\!\!Z\else$Z\!\!\!Z$\fi}}
\def\Zd{\Z^d}
\def\Q{{\relax\ifmmode I\!\!\!\!Q\else$I\!\!\!\!Q$\fi}}
\def\C{{\relax\ifmmode I\!\!\!\!C\else$I\!\!\!\!C$\fi}}
\def\Rd{\R^d}
\def\gsim{\mathrel{\raisebox{-4pt}{$\stackrel{\textstyle>}{\sim}$}}}
\def\sime{\raisebox{0ex}{$~\stackrel{\textstyle\sim}{=}~$}}
\def\lsim{\raisebox{-1ex}{$~\stackrel{\textstyle<}{\sim}~$}}
\def\div{\mbox{div}}
\def\gr{\nabla}
\def\curl{\mbox{ curl }}
\def\M{M}  \def\NN{N}                  
\def\L{{\ell}}               
\def\Le{{\ell_1}}            
\def\Lz{{\ell_2}}
\def\Let{{\tilde\ell_1}}     
\def\Lzt{{\tilde\ell_2}}
\def\Ltw{\ell_\tau^w(\cJ)}
\def\tLtw{{}_t\Asl}
\def\t#1{\tilde{#1}}
\def\la{\lambda}
\def\La{\Lambda}
\def\ga{\gamma}
\def\Ga{\eta}
\def\al{\alpha}
\def\cO{{\cal O}}
\def\cA{{\cal A}}
\def\cC{{\cal C}}
\def\cF{{\cal F}}
\def\bu{{\bf u}}
\def\bbf{{\bf f}}
\def\bz{{\bf z}}
\def\bZ{{\bf Z}}
\def\bI{{\bf I}}
\def\bN{{\bf N}}
\def\cE{{\cal E}}
\def\cD{{\cal D}}
\def\cI{{\cal I}}
\def\cJ{{\cal J}}
\def\cM{{\cal M}}
\def\cN{{\cal N}}
\def\cT{{\cal T}}
\def\cU{{\cal U}}
\def\cV{{\cal V}}
\def\cL{{\cal L}}
\def\cB{{\cal B}}
\def\cG{{\cal G}}
\def\cK{{\cal K}}
\def\cS{{\cal S}}
\def\cP{{\cal P}}
\def\cQ{{\cal Q}}
\def\cR{{\cal R}}
\def\cU{{\cal U}}
\def\bL{{\bf L}}
\def\bK{{\bf K}}
\def\bC{{\bf C}}
\def\X{X\in\{L,R\}}
\def\ph{{\varphi}}
\def\D{{\Delta}}
\def\H{{\cal H}}
\def\bM{{\bf M}}
\def\bQ{{\bf Q}}
\def\bG{{\bf G}}
\def\bP{{\bf P}}
\def\bW{{\bf W}}
\def\bT{{\bf T}}
\def\bV{{\bf V}}
\def\bv{{\bf v}}
\def\bz{{\bf z}}
\def\bw{{\bf w}}
\def\bQ{{\bf Q}}
\def\t{\tilde}
\def\lll{\langle}
\def\rr{\rangle}
\def\dJ{\nabla}
\newcommand{\ba}{{\bf a}}
\newcommand{\bb}{{\bf b}}
\newcommand{\bc}{{\bf c}}
\newcommand{\bd}{{\bf d}}
\newcommand{\bs}{{\bf s}}
\newcommand{\bff}{{\bf f}}
\newcommand{\bp}{{\bf p}}
\newcommand{\bg}{{\bf g}}
\newcommand{\bq}{{\bf q}}
\newcommand{\br}{{\bf r}}
\newcommand{\beze}{{\bf 0}}
\newcommand{\PP}{\mathbb{P}}
\newcommand{\gkl}{\Gamma_{kl}}
\newcommand{\Hoo}{H^{1/2}_{00}}
\newcommand{\Hoogkl}{\Hoo(\gkl)}
\newcommand{\Asl}{{\cA^s}}
\newcommand{\AsH}{{\cA^s(\cH)}}
\newcommand{\Aslt}{{\cA^s_{\rm tree}}}
\newcommand{\AsHt}{{\cA^s_{\rm tree}(\cH)}}
\newcommand{\beqn}{\begin{equation}}
\newcommand{\eeqn}{\end{equation}}
\newcommand{\be}{\begin{equation}}
\newcommand{\ee}{\end{equation}}
\def\beginproof{\noindent{\bf Proof:}~ }
\def\endproof{\hfill\rule{1.5mm}{1.5mm}\\[2mm]}

\newenvironment{Proof}{\noindent{\bf Proof:}\quad}{\endproof}

\renewcommand{\theequation}{\thesection.\arabic{equation}}
\renewcommand{\thefigure}{\thesection.\arabic{figure}}

\makeatletter
\@addtoreset{equation}{section}
\makeatother

\newcommand\abs[1]{\left|#1\right|}
\newcommand\clos{\mathop{\rm clos}\nolimits}
\newcommand\trunc{\mathop{\rm trunc}\nolimits}
\renewcommand\d{d}
\newcommand\dd{d}
\newcommand\dist{\mathop{\rm dist}}
\newcommand\diam{\mathop{\rm diam}}
\newcommand\cond{\mathop{\rm cond}\nolimits}
\newcommand\eref[1]{(\ref{#1})}
\newcommand\Hnorm[1]{\norm{#1}_{H^s([0,1])}}
\def\int{\intop\limits}
\renewcommand\labelenumi{(\roman{enumi})}
\newcommand\lnorm[1]{\norm{#1}_{\ell_2(\Z)}}
\newcommand\Lnorm[1]{\norm{#1}_{L_2([0,1])}}
\newcommand\LR{{L_2(\R)}}
\newcommand\LRnorm[1]{\norm{#1}_\LR}
\newcommand\Matrix[2]{\hphantom{#1}_#2#1}
\newcommand\norm[1]{\left\|#1\right\|}
\newcommand\ogauss[1]{\left\lceil#1\right\rceil}
\newcommand{\QED}{\hfill \raisebox{-2pt}{\rule{5.6pt}{8pt}\rule{4pt}{0pt}}%
  \smallskip\par}
\newcommand\Rscalar[1]{\scalar{#1}_\R}
\newcommand\scalar[1]{\left(#1\right)}
\newcommand\Scalar[1]{\scalar{#1}_{[0,1]}}
\newcommand\Span{\mathop{\rm span}}
\newcommand\ugauss[1]{\left\lfloor#1\right\rfloor}
\newcommand\with{\, : \,}

\newcommand\bnull{{\bf 0}}
\newcommand\bA{{\bf A}}
\newcommand\bB{{\bf B}}
\newcommand\bR{{\bf R}}
\newcommand\bD{{\bf D}}
\newcommand\bE{{\bf E}}
\newcommand\bF{{\bf F}}
\newcommand\bH{{\bf H}}
\newcommand\bU{{\bf U}}
\newcommand\cH{{\cal H}}
\newcommand\sk{{sk}}
\newcommand{\XX}{\mathbb{X}}
\newcommand{\id}{{\rm id}}
\newcommand{\Hom}{\Omega}

\newif\ifNZB
\newcommand\NZB[1]{\ifNZB \marginpar{\raggedright \scriptsize NZB:\\#1}
 \else \fi}
\newcommand{\FText}[1]{\mbox{#1}}
\makeatletter
\newcommand{\tr}{{\mathop{\operator@font T}\nolimits}}
\makeatother
\newcommand{\UArrow}[4]{
 \begin{array}{ll}
  #1&\stackrel{#2}\longrightarrow\\
  #3&\;\raisebox{1ex}{$\nearrow$}\mkern-14mu_{#4}
 \end{array}}
\newcommand{\DArrow}[4]{
 \begin{array}{ll}
  \stackrel{#1}\longrightarrow&#2\\
  \mkern-10mu_{#3}\mkern-22mu\raisebox{1ex}{$\searrow$}&#4
 \end{array}}
\newcommand{\fig}[3]{\par\begin{figure}[ht]
  \centerline{\epsfbox{#1.eps}}\caption{#3}\label{fig#2}\end{figure}}

\newcommand{\ubk}{u^{\rm bk}}
\newcommand{\utr}{u^{\rm true}}
\newcommand{\Atr}{A^{\rm true}}
\newcommand{\dI}{\Delta}
\newcommand{\cY}{\mathcal{Y}}
\newcommand{\U}{\mathbb{U}}
\newcommand{\V}{\mathbb{V}}
\newcommand{\W}{\mathbb{W}}
\newcommand{\pP}{\mathbb{P}}
\newcommand{\Um}{\mathbb{U}_{\rm mixed}}
\newcommand{\Vm}{\mathbb{V}_{\rm mixed}}
\newcommand{\oVm}{\tilde{\mathbb{V}}_{\rm mixed}}
\newcommand{\Y}{\mathbb{Y}}

\newcommand{\E}{\mathbb{E}}
\newcommand{\Lip}{{\rm Lip}}
\newcommand{\cX}{\mathcal{X}}

\def\argmin{\mathop{\rm argmin}}

\def\argmax{\mathop{\rm argmax}}

 \newcommand{\tripnorm}[1]{|\!|\!| #1|\!|\!|}
 \newcommand{\w}{\mathsf{w}}
  \newcommand{\di}{\mathsf{d}}
   \newcommand{\Di}{\mathsf{D}}
 
 \newcommand{\dN}{\mathcal{N}\!\mathcal{N}}
\newcommand{\cZ}{\mathcal{Z}}
\renewcommand{\tr}{{\rm train}}
\newcommand{\hx}{{\hat{x}}}
\renewcommand{\H}{\mathbb{H}}
\newcommand{\s}{\mathsf{s}}
\newcommand{\utau}{{\underline{\tau}}}
\newcommand{\otau}{{\overline{\tau}}}

\newcommand{\F}{\mathbb{F}}

\newcommand{\ep}{\epsilon}

\renewcommand{\ss}{\textsc{s}}

\newcommand{\cAg}{{\cA^{\gamma,\ss}}}
\newcommand{\cAgL}{{\cA^{\gamma,\ss}}}
\newcommand{\vt}{\vartheta}
\newcommand{\Agt}{{\cA^{\gamma,\vartheta}}}

\renewcommand{\X}{\XX}
\newcommand{\fN}{\mathfrak{N}}
\newcommand{\QAL}{Q_{\mbox{\tiny A,L,T}}}

\newcommand{\hcN}{\widehat{\cN}}
\newcommand{\mfr}{\mathfrak{M}}
\newcommand{\eps}{\epsilon}
 \newcommand{\bxi}{\bar\xi}
\newcommand{\uxi}{\underline{\xi}}
\newcommand{\uomega}{{\underline{\omega}}}

\newcommand{\In}{{\widehat{I}}}
\newcommand{\wt}{\widetilde}
\newcommand{\Ca}{C_{a}}
\newcommand{\Ce}{C_e}
\newcommand{\Tin}{{\widehat{T}}}

\maketitle

\begin{abstract}
Approximating functions of a large number of variables poses particular challenges
often subsumed under the term ``Curse of Dimensionality'' (CoD). Unless the approximated 
function exhibits a very high level of smoothness the CoD can be avoided only
by exploiting some typically hidden {\em structural sparsity}. In this paper we propose a general framework for new model classes of functions in high dimensions.
They are based on suitable notions of {\em compositional dimension-sparsity} quantifying, on a continuous level, approximability by compositions with
certain structural properties. In particular, this describes scenarios  where  deep  
neural networks can avoid the CoD.
The relevance of these concepts is demonstrated for {\em solution manifolds} 
of parametric transport equations. 
For such PDEs parameter-to-solution maps do not   enjoy the type of high order regularity  that helps to avoid the CoD by more conventional methods
 in other model scenarios. Compositional sparsity is shown to serve as the key
 mechanism for proving that sparsity of problem data is inherited in a quantifiable
 way by the solution manifold. In particular, one obtains  convergence rates for deep neural network realizations   showing that
 the CoD is indeed avoided.

\end{abstract}

\vspace*{-4mm}
\section{What this  is about}\label{sec:1}

{\bf Background:} Like hardly any other topic have Deep Neural Networks (DNNs) been recently
influencing   numerous vibrant research activities 
constantly broadening the scope of applications of machine learning.
 In a nutshell,
DNNs are compositions of simple finitely parametrized mappings, each one being
in turn a composition of an affine map and a componentwise acting nonlinear 
``activation'' function. Thus, DNNs are highly nonlinear mappings differing in many respects
from  classical approximation systems, see \cite{DDFHP,DHP} for DNNs in the general context of nonlinear approximation.

The empirically observed stunning success of deep networks in unsupervised and supervised
learning,   in particular, in  ``Big Data'' application scenarios, has raised  high expectations  
of a similar impact in other, more science related areas which are, however, typically less error-tolerant. 

A  challenge common to either regime is the need
to recover or approximate functions of a large number of variables. 
In the latter science oriented sector examples are partial differential equations (PDEs) in high-dimensional phase space like Schr\"odinger equations or Fokker-Planck equations,   describing the evolution of 
probability distributions. Similarly,  the formulation of ``background models''
in terms of parameter dependent families of PDEs is  ubiquitous in Uncertainty Quantification, especially, in connection with   inverse tasks like state- or parameter-estimation. Here, the states of interest are functions of time, space, and 
parametric variables. They form what is often called the ``solution manifold'' 
comprised of all parameter dependent solutions obtained when traversing the 
parameter domain. We will interpret this notion in a broader sense including the
dependence of solutions on   ``problem data'' such as
 coefficients, initial or boundary conditions, and right hand sides, ranging 
 over suitable compact sets.
Thus, such data  either have already a parametric form or can be approximated within any given tolerance by expressions depending on a finite number of parameters. In either case, the states in solution manifolds can be viewed as functions of spatio-temporal as well as (possibly infitely) many parametric variables. Being able to efficiently explore such solution manifolds, i.e., to
approximate the {\em parameter-to-solution map} is at the heart of numerous applications.

Generally, the so-called {\em Curse of Dimensionality} (CoD)
describes an  obstruction, encountered in such high-dimensional approximation tasks
(meaning approximation of functions of many variables).  It roughly expresses 
an exponential dependence of computational complexity on the number of variables.
A startling early result by Novak and Wojniakowsky \cite{NW} states, for instance, that 
the ``information complexity'' suffers from the CoD even when  the model class
consists of all infinitely differentiable functions on $[0,1]^d$ with
all derivatives bounded by one and  approximation cost is measured by the
number of functional evaluations used to construct a reconstruction algorithm.  
This means the classical paradigm of characterizing approximability of functions of
a few variables by smoothness is no longer relevant in high-dimensional approximation.
It also hints at the fact that avoiding the CoD does not depend on a particular method alone but refers to the interplay or combination of a particular {\em model class} 
of approximands and a particular {\em approximation/information system} applied
to this class.

A general perception is that an increasing level of nonlinerarity of an approximation
method increases the chance to avoid or mitigate the CoD hopefully for a larger scope of model classes. A by now famous result concerns the greedy construction of shallow neural networks with 
a single hidden layer \cite{Barron}. It is known to realize dimension independent approximation rates
when the model class is a {\em Barron class}, see e.g. \cite{BCDD,SX}. Note though that this model class becomes ``smaller'' with increasing dimension as smoothness needs to be proportional to   the spatial dimension.   In view of their celebrated universality and even
higher level of nonlinearity, deep neural networks (DNNs) are therefore often perceived as providing. the ``silver bullet'' for high-dimensional
approximation. This is one motivation behind  intense research over the past years on the {\em expressive power} of DNNs which can roughly be grouped  as follows: (i) ```emulation'' strategies, (ii) super-convergence 
results, (iii) approximating {\em parameter/data-to-solution maps}. 

(i) is based on proving that classical approximation tools have efficient representations 
in terms of DNNs. As a consequence the approximation power of nearly all known ``classical'' methods are (perhaps up to some log-factors) matched by deep neural networks (see e.g. \cite{DDFHP,KPRS,Yar,Schmid-Hieber,DHP,PV,GPEB,GKNV,GKP,OpSchwZe}   and the general exposition \cite{DHP} in the context of nonlinear approximation). 
This applies to classical generic model classes, defined through smoothness properties
as well as to known ways of capturing solution manifolds. Regarding the CoD these
results show success only when the emulated schemes do too.

What goes beyond? It is remarkable that DNNs can well approximate very rough (fractal-like)  functions
as well as very smooth ones like holomorphic mappings, \cite{Yar,OpSchwZe}.
Even more stunning are the results under (ii)
establishing strictly better DNN approximation rates for certain smoothness classes than those obtained by classical methods, see e.g. \cite{Yar,Shen}. But even these results come with a grain of salt. The finding in \cite{CDPW} on
so-called  {\em stable widths}  show
that algebraic convergence rates of any (nonlinear) method  that is stable in a certain sense, applied to a given compact class,  is lower bounded by the rates of {\em entropy numbers} for that class. The entropy numbers of classical smoothness classes, in turn, are known to suffer from the CoD. This questions whether
the super-convergence rates (which wouldn't break the CoD for smoothness
classes either) can be stably realized. 

Nevertheless,
the prospect of having a single approximation system that can address essentially any 
approximation problem with near-optimal expressive power, otherwise realized only
by specialized systems, may strike one as a significant advantage.
However, 
  classical methods come typically with well-understood {\em stable}
algorithmic realizations. By contrast,   training DNNs
is haunted by a significant remaining uncertainty of optimization success let alone 
a satisfactory stability assessment.
In fact, it is shown in \cite{GV} that (for ReLU networks) theoretically possible high approximation orders cannot be realized based only on point samples as training input. 
So, one may wonder: \\[1.2mm]  
{\em for which model classes of interest can DNNs 
avoid the Curse of Dimensionality while state of the art methods don't?}\\[1.2mm]
in which case one may see good reasons beyond ``black box comfort'' to put  up with the disadvantages of DNN approximation. This turns the spotlight to less generic model classes like those under (iii). In fact,
 interesting rigorous results have  been  obtained 
for DNN approximations of solutions of certain {\em spatially high-dimensional} 
partial differential equations (PDEs) like Kolmogorov's of Hamilton-Jacobi-Bellman equations
where the use of DNNs is shown to break the Curse of Dimensionality,  
see e.g. in \cite{BBGJJ,GHJW}. Especially    
 \cite{GH} quantifies DNN approximability of solutions when the problem data
 are given as   DNNs which is close to the point of view taken here as well. 
 
 Comparable results about DNN approximation of solution manifolds
 in the context of parameter dependent PDEs seem to be so confined to emulation,
 see \cite{KPRS}. That is, 
  there exist alternate methods performing at least as well.
 These include   low-rank approximation, reduced basis methods,  or sparse polynomial
approximations, \cite{BD,CDacta,CDDN}. These methods exploit   the fact that the Kolmogorov $n$-widths of the solution manifolds decay robustly with respect to the parametric dimension. This, in turn, hinges on the holomorphy of the underlying parameter-to-solution maps warranting an efficient approximability from {\em linear spaces}. This is a typical asset of {\em elliptic} models and their close
relatives, see \cite{CDacta} and references therein.

Examples of PDEs,  so to speak, at the extreme other end of the spectrum are {\em transport equations} where no dissipative effects come to aid. Unless imposing
rather restrictive assumptions on domain smoothness or initial conditions subject to compatibility constraints, for affine parameter representations (see \cite{HS}), 
parameter-to-solution maps 
are ``generically'' no longer holomorphic and have low regularity.
In fact, the approximation of solutions to parameter dependent
transport equations by   (DNNs) has been analyzed already in the interesting work  \cite{PP}. 
It is important to note though that the 
convergence rates established  there still  reflect the full Curse of Dimensionality.
Roughly, achieving accuracy $\e$ requires a network complexity of the
order $\e^{-\frac{m+1+d_y}{\alpha}}$ where $m+1, d_y$ denote the number of
space-time and parametric variables and $\alpha$ represents the 
smoothness of the solutions resulting from suitable smoothness
assumptions on the problem data (initial conditions, right hand sides, convection
coefficients). In the light of the preceding discussion, this is, Under such generic assumptions,  best possible.
\\[1.4mm]
{\bf Objectives:}
The angle taken in this paper, regarding solution manifolds of parametric transport equations, is therefore quite different. It is more in the spirit of classical regularity theory for PDEs where smoothness of solutions is deduced from
smoothness of problem data. Only, as indicated earlier, in high-dimensions, smoothness needs to be replaced by a new type of ``regularity'' that manifests itself by some suitable notion of ``structural sparsity''. The objective of this paper is therefore twofold: (1)
establish such ``heredity results'' - sparsity of problem data implies sparsity
of solutions - for parametric transport equations forming an extreme adversary to
existing model reduction concepts; (2) develop corresponding structural sparsity 
notions, whose relevance is argued in the end to go beyond the particular case study
of transport equations.\\[1.4mm]
\noindent
{\bf Layout:}  The   relevance of {\em compositions} for neural approximation is by no means new and rather apparent, see
e.g.  \cite{E,Shen,MasPogg,Schmid-Hieber}. Although somewhat inspired by the work in \cite{Schmid-Hieber} in the context
of nonlinear regression and statistical estimation, the present approach is quite
different.  It is more in line with \cite{DDGS} where it is shown that a certain notion
of tensor-sparsity of data is inherited the solutions to high-dimensional diffusion equations. It is important to note that this sparsity notion does {\em not} require the data to have a fixed finite rank or can be described by a finite number of parameters. This is also the case here where 
the central objects   are {\em compositional
approximation classes}. These are based on what we call {\em tamed dimension-sparse compositional representations}, see \S~\ref{sec:compappr}, especially \S~\ref{ssec:tamed}, \S~\ref{ssec:appr}.
The notion of {compositional approximability} is then formulated first on a {\em continuous level} in terms of instrinsic structural properties of functions.
The principal mechanism of how to obtain from this in the end {\em finitely parametrized} DNN approximations is shown
by Theorem \ref{thm:implant} (primarily for the purpose of conceptual orientation).

With these   preparatory sections in place we formulate in \S~\ref{sec:main}
our results  concerning the DNN approximability of
solution manifolds for transport equations. In particular, we show that 
dimension-sparse compositional approximability of the data (convection field, initial conditions, right hand side) is inherited by the parameter dependent solutions.
Under such circumstances the
parameter-to-solution map admits DNN approximations that do not suffer from the
Curse of Dimensionality.  Since solutions result from compositions with
initial data and right hand side a key is to establish first efficient 
approximability of characteristics. For instance,
when the parameter dependence of the convection field is affine, the characterisfics, as functions of $m$ spatial and $d_y$ parametric variables, can
be approximated by DNNs as follows\vspace*{-1.5mm}
\be
\label{Neconv}
\|z - \cN_\e\|_{L_\infty([0,T]\times D \times \cY;\R^m)} \le \e,\quad \#\cN_\e
\lsim  d_y   \Big(\frac{e^{LT}}{\e}\Big)^{m+1}  \Big|\log_2\frac{e^{LT}}{\e}\Big|^{2}.\vspace*{-1.5mm}
\ee
Our results yield, in particular, upper bounds for {\em nonlinear manifold widths}
of corresponding solution manifolds.
All proofs  are given in \S \ref{sec:proofs}. They are at times technical and tedious which is the price for keeping track of how problem parameters effect complexity.  

In \S~\ref{sec:outlook} we close with indicating several directions of future research suggested by the present findings and their bearing on a wider problem scope.
 \vspace*{-3mm} 
\paragraph{Notational Conventions:}  
In what follows we often write $a \lsim b$ to indicate that $a$ is bounded by a constant
multiple of $b$ where the constant is independent of any parameters $a$ and $b$ may depend on, unless specified otherwise. 
Accordingly $a \eqsim b$ means $a\lsim b$ and $b\lsim a$.

For notational brevity and convenience we will use, for any pair of finite dimensional metric spaces
$X,Y$ and any
continuous function $g: X\subset \R^{d_0} \to Y\subset \R^{d_1}$,  the shorthand notation 
$$
\|g\|_\X = \|g\|_{L_\infty(X;Y)} =\sup_{x\in X}\sup_{i=1,\ldots, d_1}|g_i(x)|
=\sup_{x\in X}|g(x)|_\infty,
$$
when the 
  particular domains and ranges don't matter.
Likewise we use the domain- and dimension independent notation
$\|g\|_{\Lip_1}, |\cdot|_{\Lip_1}$ to denote the full Lipschitz norm, respectively 
semi-norm\vspace*{-0.5mm}
\be
\label{Lipnorm}
|g|_{\Lip_1} := \sup_{x,z\in X}\frac{|g(x)-g(z)|}{|x-z|_\infty}, \quad
\|g\|_{\Lip_1}:= \max\{\|g\|_\infty,|g|_\Lip\}.
\ee
Our default meaning of $|\cdot|$ for vector-valued arguments is the max-norm.
\vspace*{-3mm} 
\section{Compositions}\label{sec:compappr}
\newcommand{\Co}{\mathfrak{C}}
\newcommand{\Rs}{\cR^{res}}
\newcommand{\eqrep}{\stackrel{{\rm rep}}{=}}
\newcommand{\bi}{{\bf i}}
\vspace*{-1mm}
The role of  compositions for DNN approximation has been addressed from different perspectives in a number of earlier studies,
 see e.g., \cite{E,MasPogg,Shen,Schmid-Hieber}. 
From a rather different angle the aim here  is to put forward a general framework
that allows one to quantify {\em compositional approximability}, first on a continuous infinite-dimensional functional level. Here approximants are
not yet determined by a finite number of degrees of freedom but are subjected
to certain structural conditions revolving on the notion of {\em dimension-sparsity}, see in \S~\ref{ssec:basic}. 
In \S~\ref{ssec:tamed}
we introduce metric properties of (finite) compositions, that allow us 
to introduce in \S~\ref{ssec:appr} {\em compositional approximation classes}.
Later those are to serve  as model classes for characterizing approximability by DNNs.\vspace*{-2mm} 
%
\subsection{Basic and Notation and Structural Properties}\label{ssec:basic}
\subsubsection{Compositional Representations}
We will interpret ``compositions'' - denoted in what follows (at times in slight
abuse of notation) by the symbol
 ``$\circ$'' in $(g\circ h)(z):= g(h(z))$ in a broad sense, covering iterated applications
of global operators. A prominent albeit not exclusive role is, however, played by compositions in a 
{\em pointwise sense}, defined for continuous functions.
Specifically, we  consider compositions of   mappings \vspace*{-1.5mm}
$$
g^1: D_0\subset \R^{d_0} \to \R^{d_1},\quad g^\ell: \R^{d_{\ell-1}}\to   \R^{d_\ell}, \quad \ell =2,\ldots, n, \vspace*{-1.5mm}
$$
where we always require that the last factor is linear, i.e.,
for some $\alpha_i\in \R$, $i=1,\ldots, d_{n-1}$\vspace*{-2mm}
$$
g^n=\sum_{i=1}^{d_{n-1}}\alpha_ig^{n-1}_i.\vspace*{-2mm}
$$
When  writing $h\circ g$  it will always be implicitly assumed that the two factors $g, h$ are ``dimensionally compatible'' in the above sense.
It will be convenient to abbreviate the ordered array of mappings by $\bg= (g^j)_{j=1}^n$ to provide a particular {\em realization} 
\vspace*{-1.5mm} 
 \be
 \label{G}
G (z) = (g^n\circ \cdots\circ g^1)(z)  =: G_\bg(z),\vspace*{-2mm}
\ee
of a mapping from $D_0\subset \R^{d_0}\to \R^{d_n}$.  
We sometimes simply refer to $\bg$ instead of $G_\bg$. 
 A 
prominent example of compositional representations are DNNs whose formal definition can be found in numerous 
texts, see e.g. \cite{DHP,GKP,GH,{OpSchwZe},Yar}. Here we are content with mentioning that DNN realizations (denoted by $\cN$) are (in their simplest feed-forward version) of the form \eref{G} with factors
$g^j(\cdot)= \sigma(A^j\cdot + b^j)$ where $A^j\in \R^{d_j\times d_{j-1}}$,
$b^j\in \R^{d_j}$, and the {\em activation function} or {\em rectifier} $\sigma$ acts componentwise.
An important examples is the ReLU rectifier $\sigma(t)= \max\{0,t\}=:t_+$.  The entries in $A^j, b^j$ or in the linear output layer are called {\em weights} and
their number $\#\cN$ is the {\em size} or {\em complexity} of the DNN $\cN$.
There are numerous important ``architectural variants''like ResNet structures. We will address those as well as slight generalizations
  later below when
the need arises.
\vspace*{-2.5mm}
\subsubsection{Lipschitz-Stability}
When dealing with pointwise compositions, the first important constraint on \eref{G} in what follows is that we require
all factors $g^j$ (and hence $G$ itself) in $\bg$ to be Lipschitz continuous
with Lipshitz constants $L_j$. One reason lies 
in the following simple folklore perturbation bound.
Consider dimensionally compatible Lipschitz functions $g,h$ with constants $L_g,L_h$ and $\e_g, \e_h$ accurate approximations $\t g,\t h$.  Then\vspace*{-1.5mm}
\begin{align}
\label{2fac}
\|g\circ h - \t g\circ \t h\|_\XX& \le \|g\circ h - g\circ \t h\|_\XX + \|g\circ \t h- \t g\circ \t h\|_\XX\le L_g\e_h + \e_g.\vspace*{-1.5mm}
\end{align}
Given a compositional representation $G_\bg$ of the
form \eref{G}, we denote for $k\le n$ by $L_{[n,k]}=L_{[n,k]}(\bg)$ the Lipschitz constant of 
partial compositions $g^n\circ g^{n-1}\circ \cdots \circ g^k$. It will be convenient
to set $L_{[n,n+1]}=1$.  Then, using the above argument inductively 
yields the following familiar facts.  see e.g. \cite{Schmid-Hieber}. 
\vspace*{-1.5mm}
\begin{rem}
\label{rem:error}
Assume that 
we have mappings $\t g^j: \R^{d_{j-1}}\to \R^{d_j}$, $j=1,\ldots,n,$ such that
$$
\|g^j - \t g^j\|_\XX \le \e_j,\quad j=1,\ldots, n.
$$
Then 
$$
\|g^n\circ \cdots \circ g^1 - \t g^n\circ \cdots \circ \t g^1\|_\XX \le \e_n +  \sum_{j=1}^{n-1} \e_j L_{[n,j+1]} 
=\sum_{j=1}^{n} \e_j L_{[n,j+1]}.
$$
By a symmetric argument we can replace $L_{[n,j+1]}(\bg)$ by $L_{[n,j+1]}(\t \bg)$.
 
In terms of the individual Lipschitz constants, one has, of course, 
 $L_{[n,k]} \le \prod_{j=k}^n L_{j}.
$ 
Finally, the estimates remain valid for more general iterated applications 
of operators as long as an estimate like \eref{2fac} holds.
\end{rem}

An important context where this will be used later is the following result.

\begin{prop}
\label{prop:Lip}
For any $\delta>0$, and any $g\in \Lip_1((0,1)^\ss)$, there exsist a ReLU network
$\cN_\delta$ such that 
\be
\label{Lipbound}
\|g- \cN_\delta\|_{\infty}\le \delta ,\quad  \|\cN_\delta\|_{\Lip_1}\le c_3
(1+\|g\|_\infty)
\|g\|_{\Lip_1},
\ee
and 
\be
\label{Ndcompl}
\#\cN_\delta \le c_1 \|g\|_{\Lip_1}^\ss \delta^{-\ss}\log_2\frac{1}{\delta},
\quad \mbox{depth of } \, \cN_\delta \le c_2 \log_2\frac{1}{\delta},
\ee
 where the constants $c_1,c_2.c_3$ depend only on $\ss$. Using if necessary Lipschitz-stable continuation from bounded domains to hypercubes, analogous 
 results hold for more general domains under mild geometric constraints
 (see e.g. \cite{GKP}).
\end{prop}
%
  \cite{GKP} establishes the existence of ReLU networks
 that approximate functions of higher Sobolev regularity in
weaker Sobolev norms without compromizing the standard complexity bounds.
If we imposed more regularity  than just
  $g\in \Lip_1$, 
These results would imply \eref{Lipbound} with a constant tending to one. In the present context
we prefer to avoid assuming such ``excess regularity'' and sketch corresponding arguments in Appendix A 
for completeness,  building on
some of the concepts   in \cite{GKP}.\vspace*{-3mm}
\subsubsection{Compositional Complexity}
Since 
\be
\label{insert}
f\circ g = (f \circ h)\circ (h^{-1}\circ g) =: \tilde f \circ \tilde g,
\ee
 a mapping $G$ may have infintely many {compositional 
representations}. In slight abuse of terminology 
we sometimes write $\bg_G$ to express that $\bg$ is a representation of the 
mapping $G$. Hence, one can always ``reshape'' compositional factors by writing
where the factors in the new represenation could have unfavorable regularity or stability properties, or in fact vice versa. 
In particular, we are interested in {\em large} input-dimension $d_0$ and ask
whether specific compositional structures $\bg_G$ may help avoiding the 
Curse of Dimensionality. In this regard, the role of sparsely connected neural
networks has been already observed  in \cite{Schmid-Hieber} in the context of
nonlinear regression and has been a source of motivation for what follows.

To that end, for a scalar valued function $v:\R^d\to \R$, let $\ss (v)\le d$ denote
the actual number of variables, $v$ depends on explicitly. Then,  for $g^j:\R^{d_{j-1}}\to \R^{d_j}$ we define
\be 
\label{ssj}
\begin{array}{c}
\ss_\infty(g^j)= 0, \quad\mbox{if $g^j$ is the identity}, \quad 
\ss_\infty(g^j)= 1,\quad\mbox{if $g^j$ is multi-linear},\\
 \ss_\infty(g^j)= \max_{i\le d_j}\,\ss (g^j_i), \quad\mbox{else}.
 \end{array}
%
 \ee
Low weighting of linear, bilinear factors or the identity ``rewards'' the fact that
such factors are already finitely parametrized and further approximations are
not needed.
 
\begin{definition}
\label{def:dimsparse}
For a representation $\bg$ of a mapping $G$ let 
$$
\ss_\infty(\bg):= \max_{j=1,\ldots,n(\bg)-1}\ss_\infty(g^j).
$$
A representation $\bg_G$ of the form \eref{G} is called $\ss${\em-dimension sparse} if $\ss_\infty(\bg)\le \ss$.  
\end{definition}

This gives rise to the following   measure for the {\em compositional complexity} of
a representation $\bg$\vspace*{-2mm}
$$
\fN(\bg )  :=  \sum_{j=1}^n \sum_{i=1}^{d_j}\ss(g_i^j)   \le 
  \sum_{j=1}^{n} d_{j-1}d_j.\vspace*{-2mm}
$$
\begin{rem}
\label{rem:sharp}
Thus, when $\bg$ represents a DNN $\cN$ dimension-sparsity corresponds to sparse
connectivity and    $\fN(\bg)\eqsim \#\cN$
which will be frequently used in what follows.
\end{rem}

 We say that two representations $\bg, \bg'$ are dimensionally compatible (in this order) if $G_\bg, G_{\bg'}$ are, i.e.,
if the output dimension $d_{n(\bg)}$ agrees with the input-dimension $d_0(\bg')$.

\begin{rem}
\label{rem:comp}
For any two dimensionally compatible representations $\bg,\bg'$,
the composition $\hat G:=G_{\bg'}\circ G_\bg$ has a representation 
$\hat\bg = (\bg'|\bg)$  satisfying
\be
\label{Nplus}
\fN(\hat\bg )= \fN(\bg)+\fN(\bg').
\ee
Likewise, when $\bg, \bg'$ have equal in- and output dimension, i.e., $d_0(\bg)=d_0(\bg')$ and $d_{n(\bg)} = d_{n(\bg')}$, the sum $\hat G:=G_{\bg}+G_{\bg'}$ has a representation $\hat\bg= {\bg\choose\bg'}$ satisfying 
\be
\label{Gsum}
\fN(\hat\bg)=  \fN(\bg)+\fN(\bg').  
\ee
This   follows 
  easily by (what in the DNN context is called) parallelization (as hinted at by the notation ${\bg\choose\bg'}$), upon possibly inserting 
identity factors in the representation of smaller depth.   
\end{rem}

Given any integer $\ss$,  
$$
 \Co_{ N,\ss}(D,d') := \Big\{G=G_\bg: D\to \R^{d'}:    \fN(\bg)\le N,\,\ss_\infty(\bg)\le \ss,\, 
   g^j\in \Lip_1,\,\, 1\le j\le n(\bg)
    \Big\} 
$$
 denotes then the collection of mappings with $\ss$-dimension-sparse representations
 of complexity at most $N$.
 When $D,d'$ is  clear from the context we  write $\Co_{ N,\ss}$.

 $\Co_{ N,\ss}$ is of course not a linear set but, by Remark \ref{rem:comp},
 one has (for compatible in- and output dimensions)
\be
\label{Coplus}
\Co_{ N,\ss}+\Co_{ N',\ss}\subset   \Co_{N+ N',\ss}.
\ee
Likewise for mappings $G\in \Co_{N,\ss}(D,d')$, $\t G\in \Co_{N',\ss}(\R^{d'},d'')$
we infer from \eref{Nplus} that
\be
\label{compCo}
\t G\circ G \in \Co_{N+N',\ss}(D, d''),
\ee


The following
 remarks    motivate the discussion in the 
next section. 
\begin{rem}
\label{rem:meaning}
When $\ss\ge d_0$, $D\subset \R^{d_0}$, the constraint of $\ss$-dimension sparsity is, of course void.
In this case one simply has that $\Co_{N,\ss}(D,d')= \Lip_1(D;\R^{d'})$.
In fact, any $G\in \Lip_1(D,d')$ has a trivial representation\vspace*{-1.5mm}
\be
\label{triv}
G(x)= (g^2\circ g^1)(x), \quad g^1_i= G_i,\,\,i=1,\ldots,d'\quad g^2={\rm id}_{d'}.
\vspace*{-1.5mm}
\ee 
 So, $\ss$-dimension-sparsity with $\ss<d=d_0$ is essential for such a framework to offer interesting information. 
 \end{rem}
 \vspace*{-1.5mm}

On the other hand, dimension-sparsity per se is also not yet a sufficient
instrument which can be seen from the famous Arnold-Kolmogorov Superposition Theorem. In a variant
established by G.G. Lorentz (\cite{Lorentz}), it states that every continuous function
$G$ on $D_0= [0,1]^d$ has a representation\vspace*{-1.5mm}
\be
\label{KolLo}
G(x_1,\ldots,x_d) = \sum_{q=0}^{2d}\Phi\Big(\sum_{p=1}^d\phi_{q,p}(x_p)\Big),
\vspace*{-1.5mm}
\ee
where $\Phi, \phi_{q,p}$ are continuous functions. Thus,
$\,
G= g^4\circ g^3\circ g^2\circ g^1,\,
$
 where $g^1(x)= (\phi_{q,p}(x_p))_{q,p=0,1}^{2d,d}$, i.e.,
$d_1 = (2d+1)d$, $g^2(z)= \Big(\sum_{p=1}^d z_{q,p} \Big)_{q=0}^{2d}$, i.e., $d_2=
2d+1$, $g^3(z)= (\Phi(z_q))_{q=0}^{2d}$,  $\,d_3= 2d+1$, $g^4(z)= \sum_{q=0}^{2d}
z_q$. Hence, every $G\in \Lip_1([0,1]^d)$ belongs to $\Co_{N_d,1}([0,1]^d,1)$
where $N_d= d(2d+1)+ (2d+1) + (2d+1)+1$. In the above terminology the representation is $1$-dimension sparse.  Unfortunately, $\Phi$ could have arbitrarily low regularity beyond continuity (even for smooth $G$),
an approximation of these functions by finitely parametrized objects could still
be arbitrarily expensive. This limits the direct use of the Superposition Theorem
for practical purposes, see   
\cite{GP,Kur} for  somewhat
controversial views in this regard. In summary,
 $\ss$-dimension sparsity just by itself
is still not sufficient either. In addition,
 one needs
to {\em ``tame''} $\ss$-dimension-sparse representations to arrive at a meaningful
sparsity concept.\vspace*{-3mm}
%
\subsection{Tamed Compositions and   Compactness}\label{ssec:tamed}
 In view of the preceding comments, we   need to {\em regularize} 
 representations of elements in $\Co_{N,\ss}$. 
Assume that $\cR: \bg \mapsto \cR(\bg)\in \R_+$ complies with addition and composition in the sense that \vspace*{-1.5mm}
\be
\label{Rproperties}
\cR{\bg\choose \bg'}\le \max\{\cR(\bg),\cR(\bg')\},\quad \cR(G_\bg\circ G_{\bg'})
\le \max\big\{\cR(\bg),\cR(\bg'),\cR(\bg)\cdot\cR(\bg')\big\},\vspace*{-1.5mm}
\ee 
where we assume dimensional compatibility in the second relation. Then
 for any $G\in \Co_{N,\ss}$ let
\be
\label{tripnorm}
\tripnorm{G}_{N,\ss,\cR}=\tripnorm{G}_{N,\ss}:= \inf 
\Big\{\cR(\bg):
 G_\bg =G,\, \ss(\bg)\le\ss,\, \fN(\bg)\le N \Big\}.
\ee
We supress   reference to $\cR$ when this is clear from the context.
We refer to $\tripnorm{\cdot}_{N,\ss,\cR}$ as  {\em ``compositional norm''}
although it is not a norm 
but close to one. In fact, the following relations follow from
Remark \ref{rem:comp}, \eref{Coplus} and \eref{compCo}, combined with \eref{Rproperties}.

\begin{rem}
\label{rem:propR}
For any $G\in \Co_{N,\ss}$ and $\t G\in \Co_{\t N,\ss}$ with the same in- and output dimensions one has
\be
\label{maxplus}
\tripnorm{G +\t G }_{N+\t N,\ss}   
\le \max\{\tripnorm{G}_{N,\ss},\tripnorm{\t G}_{\t N,\ss}\}.
\ee
  Similarly, for dimensionally compatible mappings $G_i\in \Co_{N_i,\ss} $, $i=1,2$,   
 one has 
\be
  \label{G1G2Lip}
   \tripnorm{G_2\circ G_1}_{N_1+N_2,\ss  }
  \le \max\,\Big\{\tripnorm{G_1}_{N_1,\ss }, \tripnorm{G_2}_{N_2,\ss },
 \tripnorm{G_1}_{N_1,\ss }\cdot\tripnorm{G_2}_{N_2,\ss } \Big\}.
  \ee
\end{rem}

Although other variants of $\cR$ are conceivable we focus here on\vspace*{-1.5mm}
\be
\label{RLip}
\cR(\bg):= \max \big\{\| g^\ell\|_{\Lip_1}, L_{[n(\bg),\ell+1]}(\bg) :\ell=1,\ldots,n(\bg) \big\}, \vspace*{-1.5mm}
\ee
 where $L_{[n(\bg),\ell+1]}(\bg)$ denotes again the Lipschitz constant of
the {\em partial compositions} $g^{n(\bg)}\circ\cdots \circ g^{\ell+1}$. 
One easily verifies that $\cR$ satisfies \eref{Rproperties}.
Controlling $\tripnorm{\cdot}_{N,\ss }$,   obviously constrains representations of elements further. In particular,   for $\cR$, given by \eref{RLip},
\be
\label{Lipimply}
\|G\|_{\Lip_1}\le \tripnorm{G}_{N,\ss },\quad G\in \Co_{N,\ss}.
\ee

\begin{rem}
\label{rem:R0}
It will be at times useful to consider alternate weaker regularizers. A natural alternative would be 
\be
\label{R0}
\cR^\circ(\bg):= \max \big\{\| g^\ell\|_{\Lip_1}:\ell=1,\ldots,n(\bg) \big\},
\ee
so that trivially 
\be
\label{triprel}
\tripnorm{G}_{N,\ss,\cR^\circ}\le \tripnorm{G}_{N,\ss,\cR}\le \max\{1, \tripnorm{G}_{N,\ss,\cR^\circ}^N\}.
\ee

Another way of weakening
$\tripnorm{\cdot}_{N,\ss }$ is to replace the compositional
Lipschitz constants $L_{[n(\bg),\ell+1]}(\bg)$ in \eref{RLip}
by $\zeta(N) L_{[n(\bg),\ell+1]}(\bg)$, where $\zeta(N) \to 0$. This  
would permit some growth of the $L_{[n(\bg),\ell+1]}(\bg)$
and hence no longer confines compositions to stay in $\Lip_1$. In view of the 
applications to come we focus in what follows
 on the stronger version, which implies   \eref{Lipimply}, see also Remark \ref{rem:general}.   
 \end{rem}

  %
 \begin{prop}
 \label{prop:comp}
 For any fixed constant $B<\infty$ and any $\ss\in\N$ the set
\be
\label{BY}
\Co_{N,\ss}( B ) := \{ G\in \Co_{N,\ss} : \tripnorm{G}_{N,\ss }\le B\}
\ee
is compact in $\XX_0= L_\infty(D)$ for either regularizor $\cR$ or $\cR^\circ$.
Hence, a minimizing representation in \eref{tripnorm} exist. 
 \end{prop}
\vspace*{-1.5mm}
The proof of Proposition \ref{prop:comp} is given in Appendix A.

Again, in view of Remark \ref{rem:meaning} and \eref{Lipimply}, compactness of $\Co_{N,\ss}(B)$
 is trivial when $\ss\ge d_0$ or $\cR$ is given by \eref{RLip}.
 
 We are now prepared to formulate {\em compositional approximability}.

\subsection{Approximation Classes}\label{ssec:appr}
%
Unless stated otherwise \eref{RLip} is used in the definition of \eref{tripnorm}
in what follows.
In the spirit of \cite{DDGS},  consider the ``K-functional''  
\be
\label{K}
K_\ss(v,N,\delta) 
 := \inf_{G\in \Co_{N,\ss} } \| v - G\|_\X + \delta \tripnorm{G}_{N,\ss }.
 \vspace*{-1.5mm}
\ee
%
%
Obviously, $K_{\ss }(v,N,\delta)\le C\delta$ means that for some $\bar G\in \Co_{N,\ss } $   one has $\|v-\bar G\|_\X\le C\delta$
and $\tripnorm{\bar G}_{N,\ss  } \le C$,
 i.e., accuracy $\delta$ is achieved with a controlled ``composition-norm''.  

Interrelating $N$ and $\delta$ is then 
a way to define collections of functions with a certain quantifiable 
{\em compositional approximability}. A (smooth) strictly increasing function $\gamma :\R_+ \to \R_+$ with $\lim_{s\to \infty}\gamma(s) =\infty$, is called 
a {\em growth function}. Its inverse $\gamma^{-1}$ exists and is also a growth function. This is to be distinguished from $\gamma(\cdot)^{-1}= 1/\gamma(\cdot)$. Given such a growth function $\gamma$, 
we consider  the class
\be
\label{apprclass}
\cA^{\gamma,\ss}  
:= \{ v\in \mathbb{X} : \|v\|_{\cA^{\gamma,\ss} } :=\|v\|_\infty +|v|_{\cA^{\gamma,\ss} }<\infty\},
\ee
where
$\, |v|_{\cA^{\gamma,\ss} }:=   \sup_{N\in\N} \gamma(N) K_{\ss }(v,N,\gamma(N)^{-1}).
$
 \begin{rem}
 \label{rem:Adiff}
 The definition of $\cA^{\gamma,\ss}$ makes sense for any regularizer satisfying
 \eref{Rproperties}. If we want to specify any other regularizer than the default one
 \eref{RLip} we indicate this by a corresponding subscript such as e.g. $\cA^{\gamma,\ss}_{\cR^\circ}$.
 \end{rem}
 
 Since trivially
 $$
 \|v- G_N\|_\X\le \gamma(N)^{-1}\gamma(N)\big\{\|v- G_N\|_\X + \gamma(N)^{-1}
 \tripnorm{G_N}_{N,\ss  }\big\}\le \gamma(N)^{-1}|v|_{\cA^{\gamma,\ss} },
 $$
we will use the information $v\in \cA^{\gamma,\ss} $ often in the form
 that for each $N\in\N$ there exists $G_N\in \Co_{N,\ss} $ such that 
 \be
 \label{estA}
 \|v- G_N\|_\X\le \gamma(N)^{-1}|v|_{\cA^{\gamma,\ss} },\quad \tripnorm{G_N}_{N,\ss }\le \|v\|_{\cA^{\gamma,\ss} },\quad N\in \N.
 \ee
 To put it differently,  realizing target accuracy $\e$ by
a composition $G_{\bg_\e}$ is achievable within a complexity $N=N_\e$ of the order
\be
\label{Ne}
N_\e = \ \left\lceil \gamma^{-1}\Big(\frac{ |v|_{\cA^{\gamma,\ss} }}{\e}\Big) \right\rceil.
\ee

A natural candidate for a {\em model class} can now be defined as follows.
Let for any fixed $B<\infty$  
\be
\label{K}
\cK_{\gamma,\ss }(B)  := \{v\in \cA^{\gamma,\ss}  : \|v\|_{{\cA^{\gamma,\ss} }}\le B \}.
\ee
 
Again, on account of Remark \ref{rem:meaning} the classes $\cA^{\gamma,\ss}_\cR,
\cA^{\gamma,\ss}_{\cR^\circ}$ agree with $\Lip_1$ when $\ss\ge d_0$ and hence do not
provide any useful information. In particular, by \eref{Lipimply} and \eref{estA},
balls in these classes are compact when $\ss\ge d_0$. Since $\|v\|_{\Lip_1}\le \|v\|_{\cA^{\gamma,\ss}_\cR}$ for all $\ss\in \N$, precompactness of $\cK_{\gamma,\ss }(B)$ is also immediate.  Since in general $\cA^{\gamma,\ss}_{\cR^\circ}$ need not be contained in $\Lip_1$ for $\ss<d$ the following claim requires an argument provided 
in Appendix A.
\begin{rem}
\label{rem:Acomp}
The sets   $\cK_{\gamma,\ss,\cR^\circ }(B)$
(see 
\eref{K}) are for $B<\infty$ and $\ss\in \N$ compact.
\end{rem}
\begin{rem}
\label{rem:general}
Finer scales $\cA^{\gamma,\beta,\ss}$, $\beta\in (0,1]$, of approximation classes can be obtained by defining $|v|_{\cA^{\gamma,\beta,\ss}}:= \sup_{n\in\N}
\gamma(N)^\beta K_\ss(v,N,\gamma(N)^{-1})$. This implies the existence of $G_N\in \Co_{N,\ss}$ such that $\|v-G_N\|_\X\le \gamma(N)^{-\beta}|v|_{\cA^{\gamma,\beta,\ss}}$, while $\gamma(N)^{\beta-1}\tripnorm{G_N}_{N,\ss }\le \|v\|_{\cA^{\gamma,\ss} }$. So, for $\beta <1$ the composition norms are allowed to grow like at most
$\gamma(N)^{1-\beta}$, $N\in \N$. Since in subsequent applications a uniform control
on composition norms matter, we confine the discussion henceforth to the special case $\beta=1$
\end{rem}
%
We record a few elementary properties of compositional approximation classes.
\begin{rem}
\label{rem:linear}
(a)  Obviously one has
 \be
 \label{contained}
 \gamma (\cdot)\lsim \t\gamma(\cdot)\quad \Rightarrow \quad \cA^{\t\gamma,\ss} \subseteq \cA^{\gamma,\ss} .
 \ee
(b)   Whenever $\gamma$ satisfies $\gamma(N)\ge c_\gamma  \gamma(2N)$ for   $N\in \N$, one has $\cA^{\gamma,\ss} + \cA^{\gamma,\ss}\subset \cA^{\gamma,\ss}$.
\end{rem}

\begin{rem}
\label{rem:closedness}
(a) For each $N\in \N$ (fixed) $\Co_{N,\ss}\subset \cA^{\gamma,\ss} $ for
{\em every} growth function $\gamma$.\\[2mm]
(b) Assume that $\gamma \lsim \t\gamma$ and   $\gamma(N)\ge c_\gamma  \gamma(2N)$
for a fixed $c_\gamma>0$. For respective equal input- and output-dimension and $\cR$ according to \eref{RLip}, $\cA^{\gamma,\ss}$ is closed under composition with elements from $\cA^{\t\gamma,\ss}$. An analogous statement holds when $v\in \cA^{\gamma,\ss}(D,d')$, $w\in \cA^{\t\gamma,\ss}_{\cR^\circ}(D',d'')$, ${\rm dim}\,D'=d'$.
The growth range  covers  any polynomial growth.
\end{rem}
{\bf Proof:} (a) is obvious. 

Regarding (b), 
 for  $v \in \cA^{\gamma,\ss} $ and $w\in \cA^{\t\gamma,\ss} $, $N\in \N$,
there exist $G_w, G_v\in \Co_{  N,\ss}$ (with respective in- and output-dimensions), 
such that $\|v- G_v\|_\XX 
\le \|v\|_{\cA^{\gamma,\ss} } \gamma(N)^{-1}$, and 
$\|w- G_w\|_\XX 
\le \|w\|_{\cA^{\t\gamma,\ss} } \t\gamma(  N)^{-1}$. Thus, since by \eref{compCo},
$G_w\circ G_v\in \Co_{2N,\ss}$, we use \eref{Lipimply} to conclude  
\begin{align}
\label{compclosed}
\| w\circ v - G_w\circ G_v\|_\XX &\le \|(w- G_w)\circ v\|_\XX + \|G_w(v)-G_w(G_v)\|_\XX
\nonumber\\
& \le   \|w\|_{\cA^{\t\gamma,\ss}}\t\gamma(N)^{-1}   + |G_w|_{\Lip_1}\|v\|_\cAg  \gamma(N)^{-1}
\nonumber\\
&\le c_\gamma^{-1} \|w\|_{\cA^{\t\gamma,\ss}} \Big\{1+ \|v\|_{\cA^{\gamma,\ss}}\Big\}\gamma(2N)^{-1}.
\end{align}
Since by \eref{G1G2Lip} and \eref{estA}, 
$\tripnorm{G_w\circ G_v}_{2N,\ss}\le \max\big\{\|w\|_{\cA^{\t\gamma,\ss}}, 
\|v\|_{\cA^{\gamma,\ss}}, \|w\|_{\cA^{\t\gamma,\ss}}\cdot 
\|v\|_{\cA^{\gamma,\ss}}\big\},
$
 (b) follows.\hfill $\Box$ 
\vspace*{-2mm}
%
 \vspace*{-2mm} 
\subsubsection{Implanting finitely-parametrized components:}
The approximation classes introduced above characterize approximability by ``tamed'' or
``regularized'' compositions which themselves are not yet
described by finitely many parameters. 
However, the compositional structure leads, in a second step, to a {\em finitely parametrized approximation}, see also \cite{Schmid-Hieber} for an approach in the same spirit.
We refer to this as ``implanting Lipschitz-stable neural networks'' which means 
that every component $g^j_i$ in a factor $g^j$ with $\ss(g^j)>1$ is replaced by a 
neural network (with in- and output-dimension one), so as to produce an expanded
composition whose factors are either at most bilinear or neural networks. So in total
they form a finitely parametrized function which we still refer to as a neural network.
Note that this does preserve $\ss$-dimension-sparsity.

The next result, although not directly applied in this form later in applications,
reflects the underlying guiding principle.
\begin{theorem}
\label{thm:implant}
Assume that $v\in \cA^{\gamma,\ss}$, then for every $\e>0$ there exists a DNN $\cN_\e$, such that
\be
\label{proto}
\|v-\cN_\e\|_\X\le \e,  
\ee
where
\be
\label{proto2}
\#\cN_\e \lsim \Big(\frac{\|a\|_{\cA^{\gamma,\ss}}}{\e}\Big)^\ss
\Big(\gamma^{-1}\Big(\frac{2|v|_{\cA^{\gamma,\ss}}}{\e}\Big)\Big)^{\ss+1}\Big|\log_2\gamma^{-1}\Big(\frac{2|v|_{\cA^{\gamma,\ss}}}{\e}\Big)\Big|.
\ee
Moreover, one has 
\be
\label{proto3}
\tripnorm{\cN_\e}_{ N_\e,\ss}\lsim \|v\|_{\cA^{\gamma,\ss}}^{N_\e},\quad 
N_\e \eqsim \gamma^{-1}\big( 2|v|_{\cA^{\gamma,\ss}}
/\e\big),
\ee
while for the weaker regularization $\cR^\circ$ one has $\tripnorm{\cN_\e}_{ N_\e,\ss,\cR^\circ}\lsim \|v\|_{\cA^{\gamma,\ss}_\cR}$.
\end{theorem}
Thus, unless $\|v\|_{\cA^{\gamma,\ss}}$ hides an exponential dependence on the
large input dimension $d_0$, unit balls $U\cA^{\gamma,\ss}$ are model classes
for which DNN approximation avoids the Curse of Dimensionality when $\ss\ll d_0$.

It is instructive to specialize these estimates for two types of growth functions
\be
\label{assa}
\mbox{(alg):}\quad \gamma(r)\eqsim \Ca r^\alpha,\quad \mbox{or}\quad
\quad \mbox{(exp):}\quad \gamma(r)\eqsim \Ce e^{\alpha r},
 \ee
for some $\alpha >0$. (With a bit more technical effort the arguments extend to
more refined scales like $\gamma(r)\sim e^{\alpha r^\beta}$, for some $0<\beta \le 1$.)
For convenience we record for frequent future use\vspace*{-1.5mm}   
\be
\label{ginverse}
\gamma^{-1}(s)\eqsim 
\left\{
\begin{array}{lll}
\Ca^{-1/\alpha}s^{1/\alpha},& \gamma \sim \mbox{(alg)},& \alpha\ge \alpha_0>0;\\
\frac{1}{\alpha} \ln \frac{s}{\Ce},& \gamma \sim \mbox{(exp)},& \alpha >0.
\end{array}\right.
\ee
Hence \eref{proto2} takes the form
\be
\label{protospec}
\cN_\e \lsim \left\{
\begin{array}{ll}
\|v\|^\ss_{\cA^{\gamma,\ss}} \Big(\frac{|v|^\ss_{\cA^{\gamma,\ss}}}{\e}\Big)^{\frac{\ss(\alpha +1)+1}{\alpha}}\Big|\log_2 \frac{|v|_{\cA^{\gamma,\ss}}}{\e}\Big|,& \gamma\sim \mbox{(alg)},\\
 \Big(\frac{\|v\|_{\cA^{\gamma,\ss}}}{\e}\Big)^\ss
\Big|\ln \frac{|v|_{\cA^{\gamma,\ss}}}{\e}\Big|^{1+\ss} \Big|\ln\ln \frac{|v|_{\cA^{\gamma,\ss}}}{\e}\Big|,& \gamma\sim \mbox{(exp)}.
\end{array}\right.
\ee
When $\gamma \sim $ (exp)  strong stability in \eref{proto3}   deteriorates only slowly 
according to $\|v\|_{\cA^{\gamma,\ss}}^{|\ln \e|}$.

In general, the stronger the algebraic growth order  the closer the dominating complexity factor comes to the rate $\e^{-\ss}$ which is attained for
exponential growth (up to logarithmic factors). This rate   is what one can expect 
for Lipschitz functions of $\ss$ variables.

The assertion of Theorem \ref{thm:implant}  hinges on the following Lemma which we
state here for later reference in several applications of similar type.
\begin{lemma}
\label{lem:import}
Assume that for some $\ss\le d_0$, the mapping $G$ belongs to $\Co_{N,\ss}$, i.e.,
is $\ss$-dimension-sparse (see Definition \ref{def:dimsparse}). Let $G=G_\bg\in \Co_{N,\ss}$. Let the DNN $\cN $ be obtained by 
 replacing each   component $g^j_i$ of each factor $g^j$ with $\ss(g^j)>1$,  by
a $\delta_j$ accurate Lipschitz stable network $\cN^j_i$, i.e., 
$$
\|g^j_i-\cN^j_i\|_{\X_i}\le
\delta_j,\quad  \|\cN_j^i\|_{\Lip_1}\le c_3
(1+\|g^j_i\|_\infty)
\|g^j_i\|_{\Lip_1},
\quad i=1,\ldots, d_j,\,j=1,\ldots, n(\bg)-1.
$$
Then one has
\be
\label{NapprLip}
\|G - \cN\|_{L_\infty} \le  \delta_{n}
+ \sum_{j=1}^{n(\Di)-1}\delta_j L_{[n(\Di),j+1]}, 
\quad \#\cN \le   \sum_{j=1}^{n(\Di)-1} \|g^j\|^\ss_{\Lip}  
 {d_j}\delta_j^{-\ss}|\log_2\delta_j| .
\ee
and 
\be
\label{Nsharp}
\#\cN \le N\tripnorm{G}_{N,\ss}^\ss 
   \max_{j=1,\ldots,n(\Di)}\delta_j^{-\ss}|\log_2\delta_j| .
\ee
 \end{lemma}
{\bf Proof:} The first relation in \eref{NapprLip} follows from  Remark \ref{rem:error} while the second one is a consequence of Proposition \ref{prop:Lip}
together with Remark \ref{rem:sharp}. 
Since  $\sum_{j=1}^{n(\bg)}d_j\le \fN(\bg)$,    \eref{Nsharp} follows from the definition
of the compositional norms.\hfill $\Box$\\

We return to the proof of Theorem \ref{thm:implant}. By \eref{estA} we can find for each $N\in\N$ a $G_N\in \Co_{N,\ss}$ satisfying \eref{estA}.
Given $\e>0$, \eref{Ne} says that $N_\e= \gamma^{-1}
\big(2|v|_{\cA^{\gamma,\ss}}/\e\big)$ (we ignore the ceil-opertor) suffices to ensure
that\vspace*{-1.5mm}
$$
\|v- G_{N_\e}\|_\X \le \frac{\e}2.\vspace*{-1.5mm}
$$
Let $\cN_{N_\e,\delta}$ denote the DNN, obtained by Lemma \ref{lem:import},
with implant-tolerances $\delta=\delta_j$ all equal. to conclude that for any
minimizing representation $\bg_{N_\e}$ of $G_{N,\e}$\vspace*{-1.5mm}
\begin{align*}
\|v-\cN_{N_\e,\delta}\|_\X &\le \frac{\e}2 + \delta n(\Di(\bg_{N_\e}))\tripnorm{
G_{N_\e}}_{N_\e,\ss} 
\le \frac{\e}2 + N_{\e}\delta\|v\|_{\cA^{\gamma,\ss}}.\vspace*{-1.5mm}
\end{align*}
 Choosing 
$
\delta = \delta(\e):=  {\e}({2\|a\|_{\cA^{\gamma,\ss}}\gamma^{-1}(2|v|_{\cA^{\gamma,\ss}}/\e)})^{-1} 
$,
produces a network $\cN_\e:= \cN_{N_\e,\delta(\e)}$ satisfying \eref{proto}, on account of \eref{NapprLip}.
Regarding \eref{proto2}, we infer now from $\tripnorm{G_{N_\e}}_{N_\e,\ss}
\le \|v\|_{\cA^{\gamma,\ss}}$ and \eref{Nsharp} that
\begin{align*}
\#\cN_\e &\le N_\e\|v\|^\ss_{\cA^{\gamma,\ss}} \delta(\e)^{-\ss}|\log_2\delta(\e)|
\\
&\lsim \gamma^{-1}\Big(\frac{2|v|_{\cA^{\gamma,\ss}}}{\e}\Big)\Big(\frac{\|a\|_{\cA^{\gamma,\ss}}
\gamma^{-1}(2|v|_{\cA^{\gamma,\ss}}/\e)}{\e}\Big)^{\ss}
\Big|\log_2
\gamma^{-1}(2|v|_{\cA^{\gamma,\ss}}/\e)\Big|,
\end{align*}
which is \eref{proto2}.

Finally, regarding the stability of the networks $\cN_\e$, 
we employ the (possibly over-pessimistic) estimate \eref{triprel}
to obtain \eref{proto3}.\hfill $\Box$

\begin{rem}
\label{rem:scope}
 Time-stepping in discretized dynamical systems is not the only context where one can expect to encounter
compositional sparsity. More generally, solutions to operator equations 
can often be approximated by {\em iterative processes} such as fixed-point iterations
that may help to assert membership to a com positional approximation class.
This is exemplified next   for a specific scenario 
where standard reduced modeling concepts suffer from    slowly decaying Kolmogorov widths.\vspace*{-2mm}
\end{rem}
 %
\section{Linear Parametric Transport Equations}\label{sec:transport}
\vspace*{-1.5mm}
\subsection{A Model Problem }\label{ssec:model}
\vspace*{-1.5mm}
We consider  the Cauchy problem for a linear  (scalar) transport equation in $m$ spatial dimensions ($m\in \{1,2,3\}$, say)
with parameter dependent data  \vspace*{-1.5mm}
\be
\label{3.1}
\begin{array}{c}
\partial_t u(t,x)+ a(t,x,y)\cdot \nabla_x u (t,x) -f(t,x,y) =0,\quad x\in \R^m,\, t\in  [0,T_\infty),\,\, y\in \cY,\\
u(0,x,y)= u_0(x,y),\quad x\in \R^m,\,\,y\in \cY,
\end{array}\vspace*{-1.5mm}
\ee
which is a standard format for models with uncertain data.
We assume for convenience that ${\rm supp}\,u_0 = \overline{D}\times \cY$ where $D$ is a bounded domain. Hence, for a fixed time horizon $T_\infty$ the solution,
as a function of $t,x $ can take values different from zero only in a bounded subset   of $[0,T_\infty)\times \R^m$.
In what follows $T_\infty$ should be viewed as a fixed finite but possibly large 
time horizon whose order of magnitude is expected to affect the complexity of the envisaged
parameter-to-solution maps $y\mapsto u(y)$. 

We shall sometimes view $u(y)$ for each $y\in\cY$  as a function of $(t,x)\in [0,\Tin)\times \R^m$, i.e.,
 as a ``point'' in $L_\infty([0,\Tin)\times \R^m)$.
 To generate in the end efficient surrogates for 
the parameter-to-solution map, it will nevertheless be useful to view
$u$ as a function of all variables $(t,x,y)$ 
\be
\label{Hom}
u:\Hom:= [0,\Tin)\times \R^m\times \cY \to \R,
\ee
so that
 for {\em high parameter-dimensionality} $\cY
\subset \R^{d_y}$,   $d_y\gg 1$,
one faces approximation problems
in high dimensions. Since we are interested in conditions other than smoothness
that may help  avoiding
the Curse of Dimensionality  we impose only low or moderate smoothness conditions on the problem data.
Specifically,  we assume throughout:\vspace*{-1.5mm}
%
\be
\label{Lip}
\begin{array}{c}
a\in C(0,\Tin;\Lip_1(\R^m\times \cY)), \quad \|a\|_{L_\infty(\Hom ;\R^m)}\le A,\\
|a(t, z;y)- a(t,z',y')| \le L\max\{ |z-z'| , |y-y'| \},\quad (t, z,y), (t, z',y')\in \Hom. 
\end{array}
\ee
In addition  we require at times in addition Lipschitz-in-time continuity
of $a$
\be
\label{Lipt}
a \in \Lip_1(0,\Tin;C(\R^m\times\cY)), \quad |a(\cdot,w)|_{\Lip_1([0,\Tin])}\le L_t.
\ee
We separate the Lipschitz-conditions \eref{Lip} and \eref{Lipt}  because    \eref{Lipt}
 is, under certain circumstances not necessary, see Remark \ref{rem:odd} in the next section.
%

To see what one can expect regarding sparsity of solutions, the very special case, where $a$ is   independent of $t, x$, is  instructive. 
E.g. when $f=0$ 
the solution $
u(t, x, y) = u_0(x- t a(y),y)$  
 is a simple {\em composition}
of $u_0$  with a linear function in $(t,x)$ involving, however,
a $y$-dependent coefficient. Even when $u_0$ does not depend on $y$ but $a(y)$
can be any element in $\Lip_1(\cY)$ the solution can be in essence  an arbitrary
Lipschitz function and stable approximations will suffer from the CoD.
The same holds, if $a(y)=a$ is constant but $u_0(\cdot,y)$ is an arbitrary element 
in a $\Lip_1(\cY)$-ball. Analogous considerations apply to the right hand side
$f$ when only smoothness conditions are imposed.  This is in agreement with the findings in \cite{PP}
where the only conditions
on the convection field are   given in terms of classical smoothness properties.

In conclusion, more specific structural constraints on the data  are needed to ensure that   $u$ can be approximated without suffering from  
 the Curse. In brief, all one can expect is a ``heredity'' effect where some structural 
sparsity of the data leads to a structural solution sparsity that allows one to
avoid the CoD. 
 \vspace*{-2mm}
\subsection{Dimension-Sparse Compositional Convection Fields}\label{ssec:a-dimsparse}
In the light of the preceding comments we consider convection fields $a$
that belong to the Bochner-type space  of functions that are continuous in time with uniformly controlled values in 
$\cA^{\gamma,m}= \cA^{\gamma,\ss}_\cR$,   $\cR$ given by \eref{RLip} 
\be
\label{a-sparse}
a\in L_\infty([0,T_\infty);\cAgL(\R^m\times\cY;\R^m)).
\ee
The regularization \eref{RLip}  is used in the definition of $\cA^{\gamma,\ss}$ because, under the above assumptions, solutions and characteristics 
belong to $\Lip_1(\Omega)$.

Here and below $m\le \ss \le m+1 + d_y$ marks 
some dimension-sparsity with respect to the total number of variables. 
The fact that we do not assume $a\in \cAgL(\Hom;\R^m)$
indicates that the time variable receives a special treatment.\vspace*{-1.5mm}
\begin{rem}
\label{rem:odd}
Time-Lipschitz continuity \eref{Lipt} is not always necessary. 
For our purposes it would suffice to know that compositional approximability 
is inherited by time-averages $ a_I(\bar z;y):= |I|^{-1}\int_Ia(s,\bar z;y)ds\in \cA^{\gamma,\ss}$, i.e., \vspace*{-2mm}
 \be
\label{av}
\|a_I\|_{\cA^{\gamma,\ss} }\lsim \|a\|_{L_\infty(0,\Tin;\cA^{\gamma,\ss} )},
 \quad \forall\,\, I\subset [0,\Tin]. \vspace*{-1mm} 
\ee
 This condition holds e.g. 
in the case of affine parameter dependence, introduced next.
\end{rem}

A particular case of interest concerns {\em affine parametric} expansions for the convection field \vspace*{-1.5mm}
\be
\label{parconv}
a(t,x;y) =  \sum_{j=1}^{d_y} y_j a_j(t,x), \quad \cY = [-1,1]^{d_y},\vspace*{-1.5mm}
\ee
 i.e., $a : \R^{m+1+d_y} \to \R^m$. Such representations arise, for instance,
 from Karhunen-Lo\'{e}ve expansions of random convection fields in which case the 
 $a_j$ have some decay properties.
 Notice   that   
 the second relation in \eref{Lip} now reads
   \vspace*{-1.5mm}
 \be
 \label{affinebound}
 \|a\|_{L_\infty( \Hom)} = \sup_{(t,x)\in \Hom}\sum_{j=1}^{d_y}|a_j(t,x)|\le A.
 \vspace*{-1.5mm}
 \ee
More specifically, we   choose the following representation format
that allows us later to explore several possible regimes\vspace*{-1.5mm}
\be
\label{aspecial}
a_j(t,w) = \omega_j a_j^\circ(t,w) ,\quad \|a_j^\circ\|_{L_\infty([0,\Tin]\times \R^m)}\le A^\circ,\quad \Lambda := \max_{\substack{ j=1,\ldots,d_y\\ 
t\in[0,\Tin]} }|a^\circ_j(t,\cdot ) |_{\Lip_1(\R^m )},\vspace*{-1.5mm}
\ee
where\vspace*{-1.5mm}
$$
\uomega =(\omega_1,\ldots,\omega_{d_y})\in \R_+^{d_y},\quad |\uomega|_1:= \sum_{j=1}^{d_y}\omega_j,\quad\mbox{so that }\,\, A= |\uomega|_1 A^\circ.\vspace*{-1.5mm}
$$
Then, one has for all $(x,y), (x',y')\in \R^m\times \cY$
\begin{align*}
|a(t,x;y)-a(t,x';y')|&\le \sum_{j=1}^{d_y}|y_j-y_j'||a_j(t,x)|+ |y_i'||a_j(t,x)-a_j(t,x')|\\
&\le A|y-y'|  + \Lambda |\uomega|_1|x-x'|,
\end{align*}
and therefore
\be
\label{aL}
 \sup_{t\in[0,\Tin]}|a(t,\cdot;\cdot)|_{\Lip_1(\R^m\times \cY)} \le A+\Lambda |\uomega|_1 =: L
\ee
is a valid Lipschitz constant permitted in \eref{Lip}.
Note that we {\em do not} require here the validity of \eref{Lipt}.

Note also that \eref{affinebound} is possible even when $\omega_j=1$, i.e., $|\uomega|_1
=d_y$ in which case $A$ and $L$ are   proportional to $d_y$.
The case $|\uomega|_1=1$ ensures a dimension-independent boundedness and regularity of the 
convection field. 

However, in either of the two ``extreme''
regimes 
{\rm (R1)}: $\,\, |\uomega|_1 = 1$ , $\,{\rm (R2)}:\,\, |\uomega|_1=d_y$,
the convection field
 $a$ is $m$-dimension sparse. More specifically, one has:
%
%
\begin{rem}
\label{rem:a}
Assume that the convection field $a$ is of the form \eref{parconv}, satisfying
\eref{affinebound} and \eref{aspecial}. 
Then $a$ has an $m$-dimension sparse compositional
representation of depth two
\be
\label{Coa}
a(t,\cdot;\cdot)\in \Co_{N_a,m },\quad \rm{depth}(a)=2,\quad N_a:=\fN(a)=d_y(1+m^2)+1,
\ee
and
\be
\label{tripa}
\tripnorm{a(t,\cdot,\cdot)}_{N,m }\le  A+\Lambda |\uomega|_1,\quad N\ge N_a.
\ee
Hence $a$  belongs to $L_\infty(0,\Tin;\cA^{\gamma,m} )$ for {\em every} growth function
$\gamma$.
\end{rem}

To see this, 
note that
\be
\label{fac}
a(t,x;y)= (g^2\circ g^1)(t,x;y),
\ee
where, in view of \eref{ssj}, for $S_A:= \{(r^1,\ldots,r^{d_y})\in \R^{md_y}:
\sum_{j=1}^{d_y}|r^j|\le A\}$,
\be
\label{gj}
\left.\begin{array}{c}
g^1:(t,x,y)\mapsto (y,a_1(t,x),\ldots,a_{d_y}(t;x))\in \R^{d_y(1+ m)},\\
 g^2:(y,r^1,\ldots,r^{d_y})\in \cY\times S_A\mapsto \sum_{j=1}^{d_y}y_j r^j,
 \end{array}\right\}\quad \fN(a)= d_y(1+ m^2)+1.
\ee
This shows \eref{Coa}. Moreover, we infer from
\eref{aL} and \eref{affinebound} that
\be
\label{Lipg}
|g^1|_{\Lip_1}\le \max\{1,\Lambda\} ,\quad, |g^2|_{\Lip_1(\cY\times S_A)}\le \big(A+\Lambda  |\uomega|_1\big)
\max\big\{|r-r'| ,|y-y'| \big\}.
\ee
Since by assumption $|g^2\circ g^1|_{\Lip_1}= |a|_{\Lip_1}\le \big(A+\Lambda  |\uomega|_1\big)=L$ we see that
uniformly in $t$, as a function of $x,y$, one has
$a\in \Co_{N_a,m }$ where $N_a:=  \fN(a)= d_y(1+ m^2)+1$. This confirms 
\eref{tripa}. \hfill $\Box$

\begin{rem}
\label{rem:track}
To reduce technicalities when tracking the dependence of constants on
problem parameters we assume from now on that\vspace*{-1.7mm}
\be
\label{LA}
1\le L_t, A \le L  ,\vspace*{-1.5mm}
\ee
because a large $L$ will be seen to have the most adverse effect.
Finally, recall  that, by definition  \vspace*{-1.7mm}
\be
\label{La}
A, L\le \|a\|_{L_\infty(0,\Tin;\cA^{\gamma,\ss} )} =: \|a\|.\vspace*{-1mm}
\ee
where this latter notational abbreviation will be used whenever reference to $\gamma,\ss$ is clear from the context.\vspace*{-1.5mm}
\end{rem} 
\subsection{ Characteristics}\label{ssec:char}
\vspace*{-1.5mm}
The field of characterics,   given by the family of ODEs\vspace*{-1.5mm}
\be
\label{3.3}
\dot{z}(t) = a(t, z(t); y) ,\quad z(0) = x,\vspace*{-1.5mm}
\ee
plays a pivotal role in what follows.
Note   that the characteristics have a natural semi-group property, 
namely that they can be obtained by composing individual characteristic segments. More precisely, suppressing the dependence on $y$ for a moment, we consider the solution of the 
more general initial value problem\vspace*{-1.5mm}
\be
\label{zgen}
\dot{z}(t, \tau;\bar z)= a(t, z(t)) ,\quad z(\tau) = \bar z.\vspace*{-1.5mm}
\ee
Later concatenations of characteristic segments necessitates including a specific
initial time $\tau$ in the notation. If $\tau=0$ and there is no risk of confusion
we  
often abbreviate $z(t,0;\bar z) = z(t;\bar z)$.  
Thus, one has for any $\tau$\vspace*{-1.5mm}
\be
\label{writez}
z(t,x) = z(t,\tau; z(\tau;x))=: (z(\cdot;\tau;\cdot)\circ z(\cdot;0;x))(t).
\vspace*{-1.5mm}
\ee
In slight abuse of terminology we refer to this as {\em composing} characteristic segments.

There is a second angle regarding compositional approximations to characteristics,
namely
  that \eref{3.3} is equivalent to the {\em fixed-point} relation\vspace*{-2mm}
\be
\label{3.8}
z(t,\tau;x;y) = x + \int_\tau^t a(s,\tau;z(s;\tau;x;y);y)ds  .\vspace*{-2mm}
\ee
Both, the semi-group property and the fixed-point relation will be combined to
construct compositional approximations to the characteristics.

Under the above assumptions characteristics don't cross, i.e., the value 
of the solution to \eref{parconv} can be determined by tracing back along
characteristics.
In fact, in view of \eqref{3.3}, one has  for the solution $u$ of \eqref{3.1} (suppressing again the dependence on $y$ for a moment)
$\frac{d}{dt} u(t, z(t)) 
 = f(t,z(t))$.
Hence, recalling that $u(0,z(0, x;y);y)= u_0(x;y)$,\vspace*{-1.5mm}
\be
\label{3.6}
u(t, z(t,x;y),y)= u_0(x;y) + \int_0^t f(s, z(s,x;y))ds ,\vspace*{-1.5mm}
\ee
or equivalently, using \eref{writez} and noting that when $x=z(t,0;\bar x)$
one has $z(s,0;\bar x) = z(-(t-s),t;x) = z(s-t,t;x)$, \eref{3.6}
takes the form\vspace*{-2mm}
\be
\label{3.6a}
u(t,x,y)= u_0(z(-t,t;x,y),y) - \int_0^t f(s,z(s-t,t, x;y))ds.\vspace*{-2mm}
\ee
  In summary, if the charateristics have ``good (pointwise) compositional approximability'', for $f=0$, the solution results from one additional compostion. 

The central objective in what follows is 
to construct finitely
parametrized surrogates $\cN(t,z,y)$ for the map
\be
\label{map}
(t,z,y)\in \Hom  \mapsto u(t,z,y),
\ee
that are determined by possibly few degrees of freedom.  
The general flavor of the following results is: membership of the problem data
(convection field, initial conditions, right hand side) to an approximation
class (see \S~\ref{ssec:appr}) imply membership of characteristics and solution
to a certain approximation class.
 \vspace*{-2mm}
\section{Main Results}\label{sec:main}
 In favor of an easier interpretability
we focus, in view of \eref{3.6a}, on the exemplary 
 types of growth functions $\gamma \sim$ (alg) and $\gamma\sim$ (exp), defined 
 in \eref{assa}. 

Since the spatial dimension $m$ is fixed and at most three we do not always mark   the dependence 
of estimates on $m$.


The proofs of  the following results can   all be found in \S~\ref{sec:proofs}.
The point of the first result is that dimension-sparsity of the convection
field is inherited by the characteristic field. 

\begin{theorem}
\label{thm:char1}
Let $\In:=[0,\Tin]$ and assume that the convection field $a$ satisfies
\eref{a-sparse} and in addition \eref{av} or \eref{Lipt}
  for some growth function $\gamma$ of either type in \eref{assa}. 
Abbreviating as before $\|a\|:= \|a\|_{L_\infty(\In; \cA^{\gamma,\ss} )}$,
one has
 \be
 \label{Mz}
  z \in \Lip_1(\Hom) \cap L_\infty(\In;\cA^{\wt\gamma,\ss} ),\quad
 \|z\|_{L_\infty(\In;\cA^{\wt\gamma,\ss} )}\lsim 
e^{\|a\|\Tin}.
 \ee
where
\be
\label{tgamma}
\wt\gamma(r):= \left\{
\begin{array}{ll}
\Big(\frac{r \Ca^{1/\alpha}}{A\Tin}\Big)^{\frac{\alpha}{1+\alpha}}
\Big(\log_2\Big(\frac{r \Ca^{1/\alpha}}{A\Tin}\Big)\Big)^{-\frac{\alpha}{1+\alpha}},&\gamma \sim \mbox{(alg)},\\
\frac{\alpha r}{A\Tin } \Big(\log_2\Big(\frac{\alpha r}{A\Tin }\Big)\Big)^{-2},& \gamma \sim \mbox{(exp)}.
\end{array}\right.
\ee
In particular,
the parameter dependent characteristic field
satisfies
\be
\label{zrate}
\inf_{\substack{\cC\in\Co_{N,\ss }\\
\tripnorm{\cC}_{N,\ss}\le e^{\Tin \|a\|}}}\|z- \cC\|_{L_\infty(\Hom)}\lsim e^{\Tin \|a\|}\wt\gamma
(N)^{-1},\quad N\in \N.
\ee
\end{theorem}

The tamed compositional  approximations in  \eref{zrate}
  are not yet characterized by a finite number of degrees of freedom which is
done in a next step similar to Theorem \ref{thm:implant}.

In what follows we adopt a generous understanding of deep neural networks
regarding the dependence on the time variable $t$. We allow in essence layers
that are piecewise affine in $t$ and hence still enjoy the basic properties of DNNs
regarding evaluation and back-propagation.

\begin{theorem}
\label{thm:dimsparse}
Under the same assumptions on the convection field $a$ and growth functions $\gamma$
according to \eref{assa} there exists for each $\e>0$ a deep neural network (DNN)
$\cN_\e$ such that
\be
\label{DNN1}
\|z- \cN_\e\|_{L_\infty(\Hom;\R^m)}\le \e,
\ee
and 
\be
\label{DNN2}
\#\cN_\e  \lsim A \Tin 2^{\ss}\|a\|^{2\ss}
\left\{
\begin{array}{ll}
 \Ca^{-\frac{1}{\alpha}} \Big(\frac{e^{\|a\|\Tin}}{\e}\Big)^{\frac{(1+\ss)(1+\alpha)}{\alpha}}
\Big |\log_2\Big(\frac{e^{\|a\|\Tin}}{\e}\Big)\Big|^2,  & \gamma \sim \mbox{(alg)},\\
\alpha^{-(1+\ss)} \Big(\frac{e^{\|a\|\Tin}}{ \e}\Big)^{1+\ss} \Big |\log_2\Big(\frac{e^{\|a\|\Tin}}{\e}\Big)\Big|^{3+\ss},& \gamma \sim \mbox{(exp)}.
\end{array}\right.
\ee

\end{theorem}

It is instructive to reformulate these findings in terms of 
{\em convergence rates}. 
\begin{cor}
\label{cor:dimsparse}
For each $N\in\N$,   there exists a DNN $\cN_N$
of complexity $\#\cN_N\le N$,
such that
\be
\label{cNest}
\|z-\cN_N\|_{L_\infty(\Hom;\R^m)}\lsim 
e^{\Tin\|a\|} \wt\gamma(N)^{-1},
\ee
where
 \be
\label{tgamma}
\wt\gamma(r)\eqsim 
\left\{
\begin{array}{ll} 
 B\, r^{\frac{\alpha}{(1+\alpha)(1+\ss)}}
\big|\log_2 r\big|^{- \frac{2\alpha }{(1+\alpha)(1+\ss)}},& \mbox{when }\,\gamma(r)\sim r^\alpha,\\
C  \alpha\, r^{\frac{1}{1+\ss}}\big|\log_2 r\big|^{-\frac{3+\ss}{1+\ss}}, &  \mbox{when }\,\gamma(r)\sim e^{\alpha r},
\end{array}\right.
 \ee
 with
 \be
 \label{BCconv}
 B= \Ca^{\frac{1}{(1+\ss)(1+\alpha)}}(A\Tin 2^\ss \|a\|^{2\ss})^{-\frac{\alpha}{(1+\ss)(1+\alpha)}},\quad
  C=  \big(A\Tin 2^\ss\|a\|^{2\ss}\big)^{-\frac{1}{1+\ss}}. 
 \ee
\end{cor}
\begin{rem}
\label{rem:stab}
It seems that one cannot  expect in general a uniform bound on the composition
norms $\tripnorm{\cN_\e}_{\#\cN_\e,\ss }$, see Theorem \ref{thm:implant} and  the comments preceding Lemma
\ref{lem:Liptxya} in \S~\ref{ssec:implant}, unless the compositional approximations 
of $a(t,\cdot;\cdot)$ have uniformly bounded depth. For the growth-types (alg) and (exp) in \eref{assa},  the following holds\vspace*{-1.3mm}
\be
\label{tripstab}
|\cN_\e|_{\Lip_1([0,\Tin];\R^m\times \cY)}  \lsim\max\{1,\|a\|\},\quad \tripnorm{\cN_\e}_{\#\cN_\e,\ss,\Lip}\lsim e^{L_\e \Tin},\quad \e>0,\vspace*{-1.3mm}
\ee
where
\be
\label{Lebound}
L_\e \lsim  
\left\{
\begin{array}{ll}
(c_3(1+A)\|a\|)^{(2/\Ca)^{\frac{1}{\alpha}}\e^{-1/\alpha} e^{\|a\|\Tin/\alpha}},&
\mbox{in case $\gamma\sim$ (alg)},\\
(c_3(1+A)\|a\|)^{ \frac{1}{\alpha}\big(|\ln \e|+\|a\|\Tin\big)}
,&
\mbox{in case $\gamma\sim$ (exp)}.
\end{array}\right.
\ee
Here $c_3$ is the constant from \eref{Lipbound}.
Thus, Lipschitz continuity with respect to $x,y$ degrades when $\e$ decreases,
  the less though, the stronger the growth order of $\gamma$.
\end{rem}

  In both theorems
the exponential case can be seen as a formal ``limit $\alpha \to \infty$'' of
algebraic rates.
 For $\ss=m+d_y$ the obtained rate would
reflect the full CoD. So, for a non-trivial dimension sparsity $\ss \ll m+1+d_y$ the CoD does not show in the convergence rates.

Nevertheless, some adverse dependence on (the potentially large) 
parametric dimension $d_y$ may still be hidden in $\|a\|=\|a\|_{L_\infty(I;\cA^{\gamma,\ss} )}$.
This will be illustrated by   specializing $a$ to affine 
representations of the form \eref{parconv}.

Recalling from  Remark \ref{rem:a}, that 
$a$ belongs to $C(\In;\cA^{\gamma,m} )$ for any growth function
$\gamma$ a first result follows from Theorem \ref{thm:dimsparse}
and Corollary \ref{cor:dimsparse} by a judicious choice of $\gamma$. 
\begin{cor}
\label{cor:special}
Assume that $a$ is of the form \eref{parconv} and satisfies 
\eref{affinebound} and \eref{aspecial}. Then,
\be
\label{aapprox}
\|a\|= \|a\|_{L_\infty(\In;\cA^{\gamma,m} )}\le 
2A +  \Lambda |\uomega|_1,
\ee
and the characteristic field belongs   to $\Lip_1(\Hom )\cap 
C(\In;\cA^{\wt\gamma,m} )$ where
\be
\label{rough}
\wt\gamma(r)\eqsim    \frac{r}{d_yA\Tin} \Big|\log_2 \frac{r}{d_yA\Tin}\Big|^{-2},
\quad \|z\|_{L_\infty([0,\Tin];\cA^{\wt\gamma,m} )}\le e^{(2A+\Lambda |\uomega|_1)\Tin}. 
\ee
Moreovery, for each $N\in\N$ there exists a   network $\cN_N$ such that
\be
\label{zratesp}
\|z-\cN_N\|_{L_\infty(\Hom\times \cY)}
\lsim d_y F 
 e^{(2A+\Lambda|\uomega|_1)\Tin} 
 N^{-\frac{1}{m+1}}\big|\log_2 N\big|^{\frac{3+m}{1+m}}.
\ee
where $F= (A\Tin 2^m(2A+\Lambda|\uomega|_1)^{2m})^{\frac{1}{m+1}}$.
\end{cor}

While the {\em rates} do not suffer from the CoD, 
to gain traction, $N$ has to exceed $d_y^{m+1}$. Although this delay effect is  
only algebraic in $d_y$, this dependence is not optimal since the choice of any
growth function for $a$ does not fully exploit the special structure \eref{parconv},
see the proof in \S~\ref{ssec:parconv2}. 
A more direct 
reasoning yields the following   better results  with regard to the 
stability of the   networks, the scaling in $\Tin$, and the dependence on $d_y$. 

\begin{theorem}
\label{thm:parconv}
Assume that $a$ is of the form \eref{parconv} and satisfies 
\eref{affinebound} and \eref{aspecial}. Recall from \eref{aL} that
\be
\label{Lspecial}
L:= A+\Lambda|\uomega|_1.
\ee
  Then, for any $\e>0$ there exists a DNN
$\cN_\e$ such that
\be
\label{Nappre}
\|z - \cN_\e\|_{L_\infty(\Hom;\R^m)} \le \e,\quad \#\cN_\e
\lsim  d_y m^2A\Tin  \Big(\frac{e^{  L\Tin}}{\e}\Big)^{m+1}  \Big|\log_2\frac{e^{  L \Tin}}{\e}\Big|^2 .
\ee 
%
Moreover, there exists a DNN $\cN_N$
with complexity $\#\cN_N\le N$ such that
\be
\label{CNNerr}
\|z-\cN_N\|_{L_\infty(\Hom;\R^m)} \lsim e^{  L\Tin} \wt\gamma(N)^{-1},\quad N\in \N,
\quad  \wt\gamma(r) =  
C \Big(\frac{r}{d_y}\Big)^{\frac{1}{m+1}}\Big|\log_2\frac{r}{d_y}
\Big|^{-\frac{2}{m+1}},
\ee
where $C=(A\Tin m^2)^{-\frac{1}{m+1}}$.
 The networks belong to $\Lip_1([0,\Tin];C(\Hom;\R^m))$ and are stable with 
$\tripnorm{\cN_N}_{N,m}\lsim e^{\hat L\Tin}$ where $\hat L \le A+\Tin^{-1}+c_3(1+A^\circ)\Lambda|\uomega|_1$ whenever
$\e\le 1$. 
%
\end{theorem}

If, on the other hand, we consider regime (R2) $|\uomega|_1=d_y$ the Lipschitz constant $\|a\|\ge L$ scales like $d_y$ so that the constant
$e^{\Tin \|a\|}$ depends exponentially on $d_y$ (see \eref{aL}). Hence,   the Curse
of dimensionality still strikes through an exponential delay in gaining accuracy.

 We  discuss next approximability of the parameter dependent
solutions themselves.  \vspace*{-1mm}
\begin{theorem}
\label{thm:still}
Under the same hypotheses on the concevtion field $a$ as in Theorem \ref{thm:parconv} 
  assume  that the data $u_0, f$ satisfy  
 \be
\label{dataappr}
u_0\in \cA^{\gamma,m} ,\,\, f\in L_\infty(\In;\cA^{\gamma,m} )\cap \Lip_1(\In;C(\R^m\times\cY)),\quad \gamma(r)\sim r^{\alpha},\vspace*{-1mm}
\ee
and let
\be
\label{where}
 \beta:= \max\Big\{1,\frac{m+1}{\alpha}\Big\}.
\ee
Then,  for any $\e >0$ there exsits a DNN $\cN_{u,\e}$ such that
for the exact solution $u$ of the transport equation \eref{3.1} \vspace*{-1.9mm}
\be
\label{uappr}
\|u - \cN_{u,\e}\|_{L_\infty(\Hom)}\le  \e,\quad
\#\cN_{u,\e}\lsim B d_y \Big(\frac{e^{ \Tin L}}{   \e}\Big)^{(m+1+\beta)}\Big|\log_2\frac{e^{L\Tin}}{\e}\Big|^2,\vspace*{-0.5mm}
\ee
where $B$ depends on $m,L,\alpha,  \max\{1,\|u_0\|,\|f\|_0\}$ 
with $\|u_0\|:= \|u_0\|_{\cA^{\gamma,m} }$,  $\|f\|:= \|f\|_{L_\infty(\In;\cA^{\gamma,m})}$.

Moreover,
for $N\in\N$ there exists a stable DNN $\cN_N$ with $\#\cN_N\le N$ such that
\vspace*{-1.5mm}
\be
\label{urate}
\|u- \cN_N\|_{L_\infty(\Hom)}\lsim e^{ L\Tin }
\underbrace{(d_yB)^{\frac{1}{m+1+\beta}}
   N^{-\frac{1}{m+1+\beta}}
\Big|\log_2\frac{N}{B}\Big|^{\frac{2}{m+1+\beta}}}_{:= \hat\gamma(N)^{-1}},  
\vspace*{-1.5mm}
\ee
with $
\tripnorm{\cN_{N}}_{N,m }\lsim \max\{1,\|u_0\|,\|f\|\} e^{L\Tin}.$
  Thus,  for $\hat\gamma$ defined in \eref{urate},
 we have $u\in C(\In;\cA^{\hat\gamma,m} )\cap \Lip_1(\In;C(\R^m\times \cY))$.
\end{theorem}

%
\begin{rem}
\label{rem:exp}
Although irrespective of the CoD the rate \eref{urate}
becomes arbitrarily bad when the algebraic order $\alpha$ gets small below 
the space-time dimension $m+1$. The best value $\beta=1$ for $\alpha\ge m+1$,
as opposed to a value tending to zero when $\alpha$ grows as in \eref{tgamma},
is due to the additional time-integration on the source field $f$. If one replaces
the algebraic order $\gamma(r)\sim r^\alpha$ in \eref{dataappr} by an exponetial
growth order $\gamma(r)\sim e^{\alpha r}$  one can show that 
$$
\|u - \cN_{u,\e}\|_{L_\infty(\Hom)}\le  \e,\quad
\#\cN_{u,\e}\lsim B d_y \Big(\frac{e^{ \Tin L}}{   \e}\Big)^{(m+2)}\Big|\log_2\frac{e^{L\Tin}}{\e}\Big|^2,
$$
with the same dependencies of $B$ on problem parameters. Finally,
we could have replaced $f\in \Lip_1(\In;C(\R^m\times\cY))$ by an assumption like \eref{av}.\vspace*{-2.5mm}
\end{rem}
%
\newcommand{\SX}{\mathbb{X}}
%
\section{Comments and Outlook}\label{sec:outlook}
\vspace*{-1.5mm}
The common trait of the above results is a uniform approximation rate 
for the characteristic field close to $\e^{-\frac{1}{m+1}}$ (the closer the stronger the approximability
order of the convection field) which is the rate one can expect for a ball
in $\Lip_1(\In\times D;\R^m)$.   For the solutions themselves it seems that one cannot
quite benefit from increasing algebraic   growth orders  beyond $\alpha=m+1$. Since again the range of possible solutions $u$ is   dense 
in a Lipschitz ball of $\Lip_1(\In\times D)$ the obtained rate in Theorem \ref{thm:still} seems to be close to optimal.
Moreover, whenever the problem data have some
compositional dimension-sparsity, in all scenarios the constructed approximations  
avoid the CoD. In general, emphasis has been on weak dependence on $d_y$ {\em not} 
on high order rates.

We conclude with indicating some ramifications of the preceding findings
whose detailed treatment is postponed to forthcoming work.
Let $\cM(a,\cY)$
denote the set of characteristic fields $z(\cdot,\cdot;y)$ obtained when $y$ traverses $\cY$ for a fixed given convection field $a$, while 
the solution
manifold $\cM(a,u_0,f_0,\cY)$   is comprised of all solutions to \eref{3.1}    for fixed data $a, u_0,f$. To capture stability with respect to those 
data as well, 
let $\mathfrak{A}$ denote the
set of all convection fields with fixed bounds for $L,L_t,A, \|a\|_{L_\infty(\In;\cA^{\gamma,\ss})}$. Likewise  let $\mathfrak{F}$ denote
the set of all $(u_0,f)$ with $\|u_0\|_{\cA^{\gamma,m}}, \|f\|_{L_\infty(\In;\cA^{\gamma,m})}\le M$. Obviously, $\mathfrak{A}$, $\mathfrak{F}$ are compact in $C(\Hom)$. For the Lipschitz-regularizer $\cR$ from \eref{RLip}
the preceding results say that all sets\vspace*{-1.5mm}
\be
\label{Ms}
\begin{array}{c}
\cM(a,\cY),\,\, \cM(a,u_0,f,\cY),\quad \cM(\mathfrak{A},\cY):= \bigcup_{a\in \mathfrak{A}}\cM(a,\cY),\\
  \cM(\mathfrak{A}\times 
\mathfrak{F},\cY):= \bigcup_{(a,u_0,f)\in \mathfrak{A}\times \mathfrak{F}}\cM(a,u_0,f_0,\cY),
\end{array}\vspace*{-1.5mm}
\ee
are contained in bounded balls of  spaces of the type $L_\infty(\In;\cA^{\wt\gamma,\ss})\cap \Lip_1(\Hom)$ for some growth function $\wt\gamma$.

A common way to characterize the complexity of these collections is to determine
their {\em metric entropy} or
suitable versions of {\em widths}, among those so-called (nonlinear) {\em manifold widths},
introduced in \cite{DHM}. More precisely, for a compact set $\cK$ in a Banach space $X$
they are defined as
\be
\label{manwidths}
\delta_n(\cK)_X := \inf_{D,E}\sup_{v\in \cK}\|v- D(E(v))\|_X,
\ee
where the infinimum is taken over all encoder-decoder pairs $E: \cK\to \R^n$, $D:\R^n\to X$, that are at least continuous. Denoting by $\theta_{N.\cdot}\in \R^N$ the collection
of weights defining the respective DNN approximations $\cN_{N,z}, \cN_{N,u}$  in Theorems \ref{thm:dimsparse},
\ref{thm:parconv}, \ref{thm:still}, respectively, the functions \vspace*{-1.5mm}
$$
D_N(t,x;\theta_N(a;y)):= \cN_{N,z}(t,x;y;\theta_{N,z}),\quad D_N(t,x;\theta_{N}(a,y,u_0,f)):= \cN_{N,u}(t,x;y;\theta_{N,u}(u_0;f)),\vspace*{-1.5mm}
$$
are valid candidates for encoder-decoder pairs $D_N\circ E_N$, where 
$E_N(a,y)= \theta_N(a;y)\in \R^N$, $E_N(a,y,u_0,f)= \theta_{N}(a,y,u_0,f))\in \R^N$,  are    the mappings that take $z(\cdot,\cdot;y)$, respectively
$u(\cdot,\cdot;y)$ into $\theta_N(a;y), \theta_{N}(a,y,u_0,f)$. 
Confining the discussion to $a$ according to \eref{parconv},  
for $\cK\in \{\cM(a,\cY),\cM(\mathfrak{A},\cY)\}$ 
the continuity 
of $E_N,D_N$ can be established based on the presented results. In fact, continuity in $y$ follows from the constructive proofs which is all that is needed for fixed $a\in \mathfrak{A}$. As a next step, continuity in $a\in \mathfrak{A}$  follows from the continuity of the construction of the implanted Lipschitz stable
 networks from Proposition \ref{prop:Lip}, as can be seen by inspecting the proof in
 Appendix A. 
 To extend these arguments  to the remaining sets in \eref{Ms}, one yet has to establish the existence of continuous metric (or near metric) selections on the level
 of dimension sparse compositional approximations prior to implanting Lipschitz stable DNNs. In particular, this would yield bounds for the manifold widths of 
compact sets of the type \eref{K}.

Knowing the manifold-widths does not allow one to infer directly on the
entropy numbers of the sets in \eref{Ms} (and hence on the number of bits needed to encode the centers of respective $\e$-covers). For a strenghtened version of
manifold widths, so called {\em stable widths}, introduced in \cite{CDPW},
a version of Carl's inequality is known which asserts that an algebraic order
of stable widths   implies the corresponding algebraic order of the 
entropy numbers. These stable widths require both factors $E_N, D_N$ to
be Lipschitz continuous. For fixed $a,u_0,f$, the above findings 
assert (uniform) Lipschitz continuity of the compositions $D_N\circ E_N$
(for $a$ of the type \eref{parconv}). It is known that DNNs
are Lipschitz continuous with respect to the weights under size constraints on the 
weights, see e.g. \cite{PW}. In general corresponding Lipschitz constants are expected to be very large which impedes an inference from approximation rates to
entropy numbers. This gives rise to the notion of Lipschitz widths studied in \cite{PW}. There, among other things, bounds on entropy numbers are derived from DNN approximation rates which are (necessarily) somewhat weaker than those in Carl's 
inequality, see \cite[\S~6.2]{PW}. Since they are derived under specific architecture constraints (either widths or depths stay bounded) they do not apply directly
to the scenarios discussed here. Specifying (and perhaps refining) such results to the current situation would be interesting as they may shed light on how the entropy
numbers of the
solution manifolds in \eref{Ms} relate to those of the accommodating balls of type
\eref{K}. Finally, it would be interesting to quantify the difference between DNN approximation spaces and compositional approximation spaces.

In a different direction, in principle, the framework allows us to treat  even less regular data leading to solutions that are no longer Lipschitz
continuous. This may require weaker regularizations than \eref{RLip} 
or refined notions of approximation classes that allow  gradually increasing 
Lipschitz constants in compositional approximations, as indicated in Remarks \ref{rem:R0}, \ref{rem:general}.
 Remark \ref{rem:scope} already indicates a wider scope of applications.
For instance, it would be interesting to apply the above concepts
to nonlinear conservation laws by exploiting their equivalent kinetic formulations
as linear parametric transport equations, see \cite{LPT}. An obvious obstacle here is that
the right hand sides are measure-valued. However, solutions do satisfy linear
transport equations with zero right hand side on regions separated by shocks.
Alternatively, one may consider constructing compositional approximations
generated through the fixed-point iterations considered in \cite{V}.
Splitting methods for more involved kinetic models may serve as another starting point for generating compositional approximations.
Finally, the above concepts apply as well to high-dimensional transport equations  and solution manifolds induced by source terms and initial conditions.
Aside from their role in Fokker-Planck equations, the correspondence bewteen 
nonlinear high-dimensional dynamical systems and linear transport PDEs opens another
interesting perspective.

Finally, one may consider the case of smooth data for which one could 
expect better rates. However, in the end one may have to resort to training concepts,
typically based on point samples to determine DNN approximations, perhaps 
in combination with pre-structured architectures suggested by the constructive proofs. It has been shown, however, in \cite{GV} that there is no hope then to
realize higher convergence orders.
\vspace*{-2.9mm} 
\section{Proofs for \S~\ref{sec:main}}\label{sec:proofs} 
\vspace*{-1.5mm}
This section consists on two major parts. In \S~\ref{ssec:prels} we collect
several technical preliminaries that will find repeated use in subsequent
sections \S~\ref{ssec:char1} -- \S~\ref{ssec:still} containing the actual   the proofs
the above theorems. Moreover, some of the elementary technical arguments, needed in these sections, are exported to an Appendix for completeness. The key idea pursued in  \S~\ref{ssec:char1} -- \S~\ref{ssec:still} is
to generate first dimension-sparse compositional approximations to the characteristic fields by combining the semi-group property \eref{writez}   with their fixed point property \eref{3.8}. Note that the lengths 
of the underlying characteristic segments - macro time steps, so to speak -  depends only on problem parameters but not
on the target accuracy which is controlled by the number of fixed point iterations. It also
shows that a straightforward discretization based on the ODEs \eref{zgen}
would provide less favorable complexity estimates.  
\vspace*{-1.9mm}
\subsection{Some Technical Prerequisits}\label{ssec:prels}
\subsubsection{``Inverting'' growth functions}\label{ssec:tech}
\vspace*{-1.5mm}
\begin{rem}
\label{rem:order}
Given $g\in \mathbb{X}$, suppose we have found for each $\e>0$    an approximation $g_\e$, depending on at most $N_\e$ degrees of freedom, that   satisfies $\|g-g_{N_\e}\|_{\mathbb{X}}\le \e$. If $N_\e\eqsim \phi(Q/\e)$
for some strictly increasing function $\phi$ of at most algebraic growth, then one has\vspace*{-1.5mm}
\be
\label{Qbound}
\|g- g_N\|_{\mathbb{X}}\lsim Q\gamma(N)^{-1},\quad N\in \N, 
\ee
where $\gamma(r)$ is any growth function satisfying\vspace*{-1.5mm}
\be
\label{constants}
\gamma(\phi(s))\eqsim s.\vspace*{-1.5mm}
\ee
We often briefly write then $\gamma \eqsim \phi^{-1}$. 
This will be  repeatedly used as follows:
Suppose the $g_N$ in \eref{Qbound} belong to $\Co_{N,\ss }$ and
$\tripnorm{g_N}_{N,\ss }\le Q$ for all $N\in \N$. Then 
$g\in \cA^{\gamma,\ss} $ with $\|g\|_{\cA^{\gamma,\ss} }\lsim Q$.
\end{rem}
To see the last conclusion, just note that
$\,\,
\gamma(N)\Big\{\|g- g_N\|_{\mathbb{X}}+ \gamma(N)^{-1}\tripnorm{g_N}_{N,\ss,\Lip}\Big\}\lsim 2Q  
$.\hfill $\Box$\\

Appropriate ``near-inverses'' $\gamma$  will be needed for growth functions $\phi$ of the following form.
\begin{lemma}
\label{lem:algrowth}
Assume that for positive $b_1, b_2, \zeta$ and real $\beta$
\be
\label{assgrowth}
\phi(s) = b_1 s^\zeta \big|\log_2 b_2 s\big|^\beta,\quad s\ge s_0>0.
\ee
Then
\be
\label{growthinverse}
\phi^{-1}(r) \eqsim b_1^{-1/\zeta}\zeta^{\beta/\zeta}r^{\frac{1}{\zeta}}
\big|\log_2(b_2^\zeta r/b_1)|^{-\frac{\beta}{\zeta}},\quad r\ge r_0>0.
\ee
\end{lemma}
{\bf Proof:} 
Making the ansatz $\phi^{-1}(r) \eqsim F r^{\frac{1}{\zeta}}|\log_2(Qr)|^\theta$, we have
\begin{align}
s &= \phi^{-1}(\gamma(s))\eqsim F \big(b_1 s^\zeta \big|\log_2(b_2 s)\big|^\beta
\big)^{\frac{1}{\zeta}}\Big|\log_2\big(Qb_1 s^\zeta\,\big|\log_2(b_2 s)\big|^\beta\big)\Big|^\theta\nonumber\\
&= F \big(b_1 s^\zeta \big|\log_2(b_2 s)\big|^\beta
\big)^{\frac{1}{\zeta}}\zeta^\theta \big|\log_2\big((Qb_1)^{\frac{1}{\zeta}}s\big)\big|^\theta \left| 1 + \frac{\big|\log_2 \big|\log_2(b_2s)\big|^\beta}{
\big|\log_2(Qb_1s^\zeta)\big|}
\right|^\theta.\vspace*{-1.5mm}
\end{align}
Equating coefficients yields
$
F=b_1^{-\frac{1}{\zeta}}\zeta^{\frac{\beta}{\zeta}},\quad Q= \frac{b_2^\zeta}{b_1}$, $ \theta = -\frac{\beta}{\zeta},
$
which confirms the claim.\hfill $\Box$\\

Note that in the above situation the proportionality constants in \eref{constants}
tend to one when the argument increases. For our purposes uniformly bounded 
proportionality constants suffice so that in later applications we can
drop the constant $\zeta^{\beta/\zeta}$ in \eref{growthinverse}.

\subsubsection{Fixed-point iterations
 }
\label{ssec:FP}
\newcommand{\uT}{{\underline{t}}}
\newcommand{\oT}{{\overline{t}}}
\newcommand{\re}{\vartheta}
In what follows we denote by $I:= [\uT,\oT]$ a fixed
time interval whose length $|I|:= \oT-\uT$ depends on $L$. 
In what follows we fix the ``macro-time-step'' $|I|$ so that  \vspace*{-1.5mm} 
\be
\label{re}
 |I| \|a\|  = \frac 12. \vspace*{-1.5mm}
\ee
For a given time horizon $\Tin$ one then needs $K:=\left\lceil \Tin/|I|\right\rceil$
such steps and we assume for convenience that $K=\Tin/|I|$ is already an integer.
In addition we denote by\vspace*{-1.5mm}
$$
\Hom(I):= I\times \R^m\times \cY\subset \Hom \vspace*{-1.5mm}
$$
the spatio-parametric time-slab determined by $I$.

To find approximate compositions we recall the fixed-point relation \eref{3.8} and   consider the corresponding mapping
 $\Phi_{x,I} : L_\infty(\Hom(I);\R^m)\to L_\infty(\Hom(I);\R^m)$, defined by,
 \vspace*{-1.5mm}
\be
\label{Phinew0}
  \Phi_{x,I}( t,\bar z;y): = x + \int_\uT^t a(s, \bar z (s),y)ds, \quad t\in I= [\uT,\oT].\vspace*{-1.5mm}
\ee

A natural strategy is to approximate the fixed point of \eref{3.8} by iterates of the mapping $\Phi_{x,I}(\cdot;t,y)$. In this case the arguments $x, \bar z$ 
sometimes depend on each other. In fact, a natural initialization would be
the constant-in-time function\vspace*{-1.5mm}
\be
\label{initial}
\bar z_x(s) = x,\quad s\in I,\vspace*{-1.5mm}
\ee
i.e., the initial value $x$ is frozen in time throughout $I$.
Then, we always use the notational convention
\be
\label{iterates}
\Phi^k_{x,I}(t,\bar z ;y):= \Phi_{x,I}\big({ t},\Phi^{k-1}_{x,I}(\cdot,\bar z ;y);y\big),\quad
\bar z\in L_\infty(I;\R^m).
\ee
Condition \eref{re} and $L\le \|a\|$ say that $\Phi_{x,I}$   
is   a contraction in $\bar z$ since\vspace*{-1.5mm}
\begin{align}
\label{contr}
|\Phi_{x,I}(t;\bar z;y)- \Phi_{x,I}(t;\bar z';y)|&\le \int_\uT^t |a(s;\bar z(s);y)-a(s;\bar z'(s);y)|ds\le (t-\uT)L \|\bar z -\bar z'\|_{L_\infty(I;\R^m)}\nonumber\\
&\le \frac 12 \|\bar z -\bar z'\|_{L_\infty(I;\R^m)}.\vspace*{-1.5mm}
\end{align} 
Since by \eref{Lip},
$
|z(t;x;y)-x|= \Big|\int_\uT^t a(s,z(s;x);y)ds\Big|\le (t-\uT) A\le A|I|\le \frac 12,
$
this implies 
\begin{align}
\label{16}
|z(t;x;y)- \Phi^k_{x,I}(t;\bar z_x;y)|&=|\Phi_{x,I}(t;z(\cdot;x);y)- \Phi_{x,I}\big(t;\Phi^{k-1}_{x,I}(\cdot;\bar z_x;y)\big)|\nonumber\\
&\le  2^{-k}\|z(\cdot;x;y)-\bar z_x(s)\|_{L_\infty(I;\R^m)}\le 2^{-k}A|I|\le 2^{-k-1},
\end{align}
where we have used $A\le L$ (see \eref{LA}) and \eref{re}.
Hence, by \eref{re}, it takes roughly $|\log_2 \eta|$ steps to achieve
 accuracy $\eta$\vspace*{-1.5mm}
\be
\label{mu}
|z(t;x;y)- \Phi^\mu_{x,I}(t;\bar z_x;y)|\le \eta  ,\quad (t,x,y)\in \Hom(I),
\quad   \mu=\mu(\eta )=
\left\lceil  \Big|\log_2( 2\eta)^{-1}\Big| \right\rceil.\vspace*{-1.5mm}
\ee
\newcommand{\ue}{{\underline{\eta}}}
In view of \eref{writez}, it is natural to concatenate next iterates
$\Phi_{\cdot,I}^\mu$ in time for successive time intervals $I$.
 To that end, consider
an equally spaced (for simplicity) partition\vspace*{-1.7mm}
$$
[0,\Tin]= \bigcup_{k=1}^K [t_{k-1},t_k],\quad t_k:= \frac{k\Tin}{K},\vspace*{-1.7mm}
$$
where $K$ is chosen in compliance with \eref{re}. Along with the sequence of intervals consider the vector of tolerances
  with corresponding sufficient iteration numbers
\be
\label{muk}
\ue^k= (\eta_1,\ldots,\eta_k)\in \R^k_+,\quad \mu_k:= \mu(\eta_k ),\,\,\, k=1,\ldots,K.
\ee
Then define, for $j<k$, $w\in\R^m$ \vspace*{-1.7mm}
\be
\label{psinew}
\begin{array}{lll}
\Psi_{[k,j]}(t,w;y):= \Phi^{\mu_k}_{w_{k-1,j},I_k}(t;\bar z_{w_{k-1,j}};y)   ,& w_{k-1,j}:= \Psi_{[k-1,j]}(T_{k-1};w;y),& t\in I_k'\\
\Psi_{[j+1,j]}(t;w;y) := \Phi^{\mu_{j+1}}_{w,I_{j+1}}(t,\bar z_w;y),
\end{array}\vspace*{-1.7mm}
\ee
i.e., $\mu_k$ iterates of $\Phi_{\cdot,I}$ are applied to the result of a 
$\mu_{k-1}$-fold application of $\Phi_{\cdot,I}$ evaluated at the last time-junction 
$t_{k-1}$.

Specifically,
\be
\label{specific}
\Psi_{\ue^k}(t,x;y):= \Psi_{[k,0]}(t,x;y)
\ee
is a natural candidate for approximating $z(\cdot,x;y)$ on $I_k$.

To estimate $|z-\Psi_{\ue^k}|$ on $I_k$ we invoke Remark \ref{rem:error}.
Viewing $\Psi_{\ue^k}$ as a perturbation of the characteristic field,
we   need
bounds for the Lipschitz constants of the exact characteristics $z(t;w;y)$.
Recall
that under the above assumptions on the 
convection field $a$, it follows from a classical Gronwall inequality that one has 
\be
\label{zLip}
\|z(\cdot,x;y)-z(\cdot,\bar x;y)\|_{L_\infty(I;\R^m)} \le e^{L|I|} |x-\bar x|
\le e^{1/2} |x-\bar x|,
\ee
so that in terms of Remark \ref{rem:error} we have 
$
L_{[k,j+1]}\le e^{L(I_{j+1}\cup\cdots
\cup I_k)} = e^{L(t_k-t_j)}=e^{(k-j)/2}.
$
Thus, for $k\le K$, $t\in I_k$\vspace*{-2.5mm}
\be
\label{psikerr}
|z(t;x;y)- \Psi_{\ue^k}(t;x;y)|\le \eta_k+\sum_{j=1}^{k-1}\eta_j e^{(k-j)/2}.\vspace*{-1.5mm}
\ee
It remains to choose the intermediate tolerances $\eta_j$.
The simplest option is to take them all equal\vspace*{-1.5mm}
\be
\label{ej}
\eta_j =\eta (\e) :=(e^{1/2}-1) \e e^{-K/2},\quad j=1,\ldots,K,\vspace*{-1.5mm}
\ee
which yields
\be
\label{uniform}
\|z(\cdot;x;y)-\Psi_{\ue^k}(t;x;y)\|_{L_\infty(I_k;\R^m)}\le 
  \e e^{(k-K)/2}\le \e, \quad k=1,\ldots,  K.\qquad
\ee
In summary, we have 
\be
\label{Ze}
\|z- Z_\e\|_{L_\infty(\Hom;\R^m)}\le \e, \quad\mbox{where}\quad Z_\e(t;x;y):= \sum_{k=1}^K\chi_{I_k}(t)\Psi_{\ue^k}(t;x;y).
\ee
  
\subsubsection{Lipschitz bounds}\label{ssec:Lipbounds}
The mappings $Z_\e$ from \eref{Ze} are still global operators. To analyze
their approximation by pointwise compositions via Remark \ref{rem:error}
we need to bound the Lipschitz constants of partial compositions.
\begin{lemma}
\label{lem:Lipprop}
Under the above assumptions one has for $k\in\N$, $\bar z,\bar z'\in L_\infty(I;\R^m), \, y,y'\in \cY$
\be
\label{globlip}
|\Phi^k_{x,I}(t,\bar z;y)- \Phi^k_{x',I}(t,\bar z';y')|
\le e^{1/2} |x-x'| + \frac{2^{-k}}{k!}
 \max\{\|\bar z- \bar z'\|_{L_\infty(I;\R^m)},|y-y'| \}.
\ee
In particular, one has for all $(t,x,y), (t',x',y')\in \Hom(I)$
\be
\label{globlipx}
|\Phi^k_{x,I}(t,\bar z_x;y)- \Phi^k_{x',I}(t',\bar z_{x'};y')|
\le A|t-t'|+ \max\{|x-x'|,|y-y'| \} e^{1/2}, \quad
  k\in \N. 
\ee
\end{lemma}
{\bf Proof:}
By our assumptions \eref{Lip} on the convection field we conclude that for $t\in I$
\be
\label{Phicontr}
|\Phi_{x,I}(t,\bar z;y) -\Phi_{x,I}(t,\bar z';y')| \le   (t-\uT)L \max\{\|\bar z - \bar z'\|_{L_\infty(I;\R^m)},|y-y'| \},
\ee
so that
$
|\Phi_{x,I}(t,\bar z;y)-\Phi_{x',I}(t,\bar z';y')|\le |x-x'|+ (t-\uT)L \max\{\|\bar z - \bar z'\|_{L_\infty(I;\R^m)},|y-y'| \}.
$ 
Hence
\begin{align*}
&\Phi_{x,I}(t;\Phi_{x,I}(\cdot;\bar z;y );y)-\Phi_{x',I}(t ;\Phi_{x',I}
(\cdot ;\bar z';y' );y' )\\
&\, =  x + \int_\uT^t a\Big(s,\Big(x+ \int_\uT^sa(s', \bar z(s');y)ds'\Big)\Big)ds  
  - \Big\{x' + \int_\uT^t a\Big(s,\Big(x'+ \int_\uT^sa(s', \bar z'(s');y')ds'\Big)\Big)ds\Big\} .
  \end{align*}
This yields
\begin{align*}
&|\Phi_{x,I}(t;\Phi_{x,I}(\cdot;\bar z;y );y)-\Phi_{x',I}(t;\Phi_{x'}(\cdot;\bar z';y' );y' )|
\\
&\qquad \le |x-x'| + \int_\uT^t \Big|a\Big(s,\Big(x+ \int_\uT^sa(s', \bar z(s');y)ds'\Big)\Big)-a\Big(s,\Big(x'+ \int_\uT^sa(s', \bar z'(s');y)y')ds'\Big)\Big)\Big|ds
\\
&\qquad \le |x-x'|+ \int_\uT^t L|x-x'| + L\int_\uT^s|a(s',\bar z(s');y)-a(s',\bar z'(s');y')|ds'ds\\
& \qquad \le
(1+(t-\uT)L)|x-x'| + L^2\int_\uT^t\int_\uT^s\max\{|\bar z(s')- \bar z'(s')|, |y-y'|\}ds'ds\\
&\qquad \le (1+(t-\uT)L)|x-x'| + \frac{((t-\uT)L)^2}2 \max\big\{  \|\bar z -\bar z' \|_{L_\infty(t;\R^d)},|y-y'|\big\}. 
\end{align*}

One then easily verifies inductively   that 
\begin{align}
\label{k0}
|\Phi_{x,I}^k(t,\bar z; y)- \Phi_{x',I}^k (t,\bar z';y')|& \le 
\sum_{\nu=0}^{ k-1} \frac{((t-\uT) L)^\nu}{\nu!}|x-x'|\nonumber\\
&\qquad + \frac{(L(t-\uT))^k}{k!}
\max\{ \|\bar z - \bar z'\|_{L_\infty(I;\R^m)},|y-y'| \}, 
\end{align}
which implies \eref{globlip}.
Specifically, when $\bar z(s)=\bar z_x(s)=x$ for $s\in I$, \eref{k0}
gives
\begin{align}
|\Phi^k_{x,I}(t,\bar z_x;y) -\Phi^k_{x',I}(t,\bar z_{x'};y')| &\le  
\sum_{\nu=0}^k\frac{((t-\uT)L)^\nu}{\nu!}\max\{|x-x'|,|y-y'|\}
\nonumber\\
&\le \max\{|x-x'|,|y-y'|\}
e^{L|I| }.
\end{align}
from which \eref{globlipx} follows for $t=t'$ since $L|I|\le 1/2$.

Since (for $t'< t$), keeping \eref{Lip} in mind 
$\,
\big|\Phi_{x,I}(t,\bar z;y)-\Phi_{x,I}(t',\bar z, ;y)\big|\le
\int_{t'}^t|a(s,\bar z(s);y)|ds \le A|t-t'| 
$,
we have\vspace*{-1.9mm}
\begin{align}
\label{totalLip}
\big|\Phi^k_{x,I}(t,\bar z;y)- \Phi^k_{x,I}(t',\bar z;y)|& \le \int_{t'}^t\big|
a(s,\Phi^{k-1}_{x,I}(s,\bar z;y);y)|ds \le  A|t-t'|,\quad t,t'\in I,\vspace*{-1.9mm}
\end{align}
proving \eref{globlipx} and hence the assertion.\hfill $\Box$\\

To approximate the $\Psi_{[k,j]}$ by pointwise compositions we
need the following bounds
%
\begin{cor}
\label{cor;Psi}
For $(t,x,y), (t',x',y')\in \Hom(I_k)$,  one has
\be
\label{bumm}
|\Psi_{[k,j]}(t,w;y)- \Psi_{[k,j]}(t',w';y')|\le |t-t'|+ e^{(k-j)/2}\max\{|w-w'|,|y-y'|\}.
\ee
Moreover, for $Z_\e$ defined by \eref{Ze}, one has
\be
\label{ZLip}
|Z_\e(t,x,y)- Z_\e(t',x',y')|\le A|t-t'|+ \max\{|x-x'|,|y-y'|\} e^{\|a\|\Tin}.
\ee
\end{cor} 
{\bf Proof:} 
Since
$$
\Psi_{[k,j]}(t,w;y)- \Psi_{[k,j]}(t',w';y')
= \Phi^{\mu_{k}}_{w_{k-1,j},I_k}(t,\bar z_{w_{k-1,j}}; y)- \Phi^{\mu_k}_{w'_{k-1,j},I_k}(t' ;\bar z_{w'_{k-1,j}};y'),
$$
where $w_{k-1,j}= \Psi_{[k-1,j]}(t_{k-1};w;y)$,
we infer from  \eref{totalLip} that
\begin{align*}
&|\Psi_{[k,j]}(t,w;y)- \Psi_{[k,j]}(t',w';y')|\le A|t-t'|\\
&\qquad\qquad + \max\big\{\big|
\Psi_{[k-1,j]}(t_{k-1};w;y) -\Psi_{[k-1,j]}(t_{k-1};w';y')|,|y-y'|\big\}
e^{1/2}.
\end{align*}
Again one concludes inductively that
\begin{align*}
&e^{1/2}\max\{|\Psi_{[k-1,j]}(t_{k-1};w;y) -\Psi_{[k-1,j]}(t_{k-1};w';y')|,|y-y'|\big\}\\
&\quad \le e^{1/2}
\max\Big\{e^{1/2}\max\big\{|\Psi_{[k-2,j]}(t_{k-2};w;y) -\Psi_{[k-2,j]}(t_{k-2};w';y')|,|y-y'|\big\},|y-y'|\Big\}\\
&\quad \le e^{(k-j-1)/2}\Big\{|\Psi_{[j+1,j]}(t_{j+1};w;y)-\Psi_{[j+1,j]}(t_{j+1};w';y')|,|y-y'|\Big\},
\end{align*}
and since
\begin{align*}
&|\Psi_{[j+1,j]}(t_{j+1};w;y)-\Psi_{[j+1,j]}(t_{j+1};w';y')| \\
&\qquad  =| \Phi^{\mu_{j+1}}_{w,I_{j+1}}(t_{j+1},t_j;\bar z_w;y)- \Phi^{\mu_{j+1}}_{w',I_{j+1}}(t_{j+1},t_j;\bar z_{w'};y')| 
\le e^{1/2}\max\{|w-w'|,|y-y'|\},
\end{align*}
\eref{bumm} follows.  

Concerning \eref{ZLip}, recall from \eref{re} that $k\le K\le 2\|a\|\Tin$.
Then, \eref{ZLip} follows for any $t,t'\in I_k$, $k\le K$, from \eref{bumm}.
The general case is again obtained by using the triangle inequality, inserting intermediate time-segments.  
This completes the proof.
\hfill $\Box$
\vspace*{-1.9mm}
\subsubsection{Pointwise compositions  }
 We wish to pass from compositions of
global operators (integral operators) to compositions
of pointwise mappings. 
Consider an equidistant partition of $I=[\uT,\oT]$ with breakpoints $\tau_i=  \tau_i(I,q):= \uT+ i |I|/q$, for some $q\in \N$.
Let $\xi_i$ denote the respective midpoints of the intervals $[\tau_{i-1},\tau_i]=: J_i=J_i(I,q)\subset I$, $i=1,\ldots,q$ and define\vspace*{-2mm}
 \be
\label{rhoi}
\rho_{i,I}(t)= \rho_i(t):= \int_\uT^t\chi_{J_i}(s)ds.\vspace*{-2mm}
\ee
The following simple facts will be used frequently.
\begin{lemma}
\label{lem:quadrature}
Adhering to the above notation, the following holds:\\[1.6mm]
(a) For $t\in I_k$
\be
\label{sumrho}
\sum_{i=1}^q \rho_i(t)= 
 \sum_{i=1}^{k-1}|J_i| + t-\tau_{k-1}\le \frac{k|I|}q.
\ee
and  
for $t, t'\in I$
\be
\label{rhodiff}
\sum_{i=1}^q|\rho_i(t)-\rho_i(t')|\le |t-t'|.
\ee
(b)
Assume that $g\in L_\infty(I)$ and let  
$
  g_{J_i}:= |J_i|^{-1}\int_{J_i}g(s)ds$, 
$i=1,\ldots,q$.
Then, 
\be
\label{om0a}
 \Big|\int_\uT^{t} g(s)ds - \sum_{i=1}^q \rho_i(t)g_{J_i}\Big|\le \frac{|I|\|g\|_{L_\infty(I)}}{2q},\quad t\in I.
\ee
 (c) Assume that $g\in \Lip_1(I)$ with Lipschitz constant $L'$. Then
 \be
\label{om0}
 \Big|\int_\uT^{t} g(s)ds - \sum_{i=1}^q \rho_i(t)g(\xi_i)\Big|\le \frac{|I|^2L'}{2q}.
\ee
\end{lemma} 
The proof is elementary and given for completeness in Appendix B.

 We approximate now $\Phi_{x,I}$ in a first step by the piecewise affine-in-time
 function\vspace*{-1.7mm}
\be
\label{Px}
P_{x,I,q}(t;\bar z;y ) := 
 x+  \sum_{i=1}^{q}\rho_i(t)\bar a_i (\bar z(\xi_i);y ),\vspace*{-1.7mm}
\ee
where, depending on our hypothesis on $a$ we set
\be
\label{acases}
\bar a_i(\bar z;y):= \left\{
\begin{array}{lll}
a(\xi_i,\bar z;y),& \mbox{(A1)}&\mbox{in case  \eref{Lipt} holds},\\
a_{J_i}(\bar z;y):=|J_i|^{-1}\int_{J_i}a(s;\bar z;y)ds,& \mbox{(A2)}& 
\mbox{in case  \eref{av} holds}.
\end{array}
\right.
\ee
%
 
We record for later use that, by \eref{sumrho}, the following analog to \eref{contr} holds\vspace*{-1.5mm}
\begin{align}
\label{PLip}
|P_{x,I,q}(t ;w;y)-P_{x,I,q}(t ;\t w;y)|&\le \sum_{i=1}^q\rho_i(t)|\bar a_{i}(w(\xi_i);y)
- \bar a_{ i}( \t w(\xi_i);y)| \nonumber\\
&\le L \|w-\t w\|_{L_\infty(I)}\frac{i(t)|I|}q\le \frac 12 \|w-\t w\|_{L_\infty(I)},
\vspace*{-1.5mm} 
\end{align}
where we have used \eref{re},     $L\le \|a\|$, and the fact that for either version of
$\bar a_i$ Lipschitz constants with respect to $\R^m\times \cY$ are preserved.

Next we estimate the deviation between $\Phi_{x,I}$ and $P_{x,I,q}$.
%
\begin{lemma}
\label{lem:PhiP}
Assume that \eqref{Lip} holds and that $\bar z\in L_\infty(I;\R^m)$ satisfies
\vspace*{-1.7mm}
\be
\label{IAest}
\|\bar z- \bar z(\xi_i)\|_{L_\infty(J_i;\R^m)}\le \frac{A|I|}{2q}.\vspace*{-1.7mm}
\ee
Then one has \vspace*{-1.7mm}
\be
\label{Pxappr}
\Big|\Phi_{x,I}( t;\bar z;y )- P_{x,I,q}(t;\bar z;y)\Big|\le 
\left\{
\begin{array}{ll}
\frac{ (1+A)L|I|^2 }{2q}  ,&  \mbox{when $\bar a_i= a(\xi_i)$},\\
\frac{(L|I|+1)A|I|}{2q} &  \mbox{when $\bar a_i= a_{J_i}$}.
\end{array}\right\}\, \le\frac{A|I|}q\le \frac{1}{2q}.\vspace*{-1.7mm}
\ee 
The second but last inequality is relevant when $L\gg A$ so that $|I|$ is correspondingly small.
\end{lemma}
{\bf Proof:} 
Let $\bar z(\uxi)$ denote the
piecewise constant $\bar z(\uxi)|_{J_i}= \bar z(\xi_i)$ to obtain from \eref{Lip}
and \eref{LA}
\begin{align}
\label{split1}
\Big|\Phi_{x,I}(t,\bar z;y)- P_{x,I,q}(t,\bar z;y)\Big|&\le 
\Big|\Phi_{x,I}(t,\bar z;y)-\Phi_{x,I}(t,\bar z(\uxi);y)\Big|
+\Big|\Phi_{x,I}(t,\bar z(\uxi);y)- P_{x,I,q}(t,\bar z;y)\Big|\nonumber\\
&\le |I|L \|\bar z(\uxi)-\bar z\|_{L_\infty(I;\R^m)} 
+ \Big|\Phi_{x,I}(t,\bar z(\uxi);y)- P_{x,I,q}(t,\bar z;y)\Big|\nonumber\\
&\le \frac{LA|I|^2}{2q}+ \Big|\Phi_{x,I}(t,\bar z(\uxi);y)- P_{x,I,q}(t,\bar z;y)\Big|,
\end{align}
where we have used \eref{IAest}. 
In case (A2), i.e.,   $\bar a_i(\cdot;\cdot)= a_{J_i}( \cdot;\cdot)$, \eref{om0a} yields,
in view of \eref{LA}, 
\be
\label{av3}
\Big|\Phi_{x,I}(t,\bar z(\uxi);y)- P_{x,I,q}(t,\bar z;y)\Big|\le \frac{A|I|}{2q}.
\ee
Thus, in this case
\be
\label{av4}
\Big|\Phi_{x,I}(t,\bar z;y)- P_{x,I,q}(t,\bar z;y)\Big|\le \frac{(L|I|+1)A|I|}{2q}.
\ee

Now suppose (A1), i.e., $\bar a_i(\cdot;\cdot)= a(\xi_i,\cdot;\cdot)$
under assumption \eref{Lipt}. Then, we apply Lemma \ref{lem:quadrature}, (c), to
$g(s)= a(s,\bar z;y)$ and, by \eref{Lipt} and \eref{LA}, ($L'\le L$), obtain
$$
\Big|\Phi_{x,I}(t,\bar z(\uxi);y)- P_{x,I,q}(t,\bar z;y)\Big|\le \frac{L|I|^2}{2q},
$$
which confirms the first inequality. On account of the assumption $1\le A\le L$ (see 
\eref{LA}, \eref{La}),   \eref{re} ensures that the first case is bounded by $A|I|/(2q)$
while the second case is bounded by $A|I|3/(4q)$. Again \eref{re} concludes the proof.
\hfill $\Box$

\begin{rem}
\label{rem:interest}
The hypothesis \eref{IAest} in Lemma \ref{lem:PhiP} is valid in the following cases:
\begin{enumerate}
\item $\bar z = \bar z_w$ for some   $w\in \R^m$ on $I$;
\item
  $\bar z$ results from applying $\Phi_{x,I}$, i.e., $\bar z(t) = \Phi_{x,I}(t;w;y )$
for some $w\in L_\infty(I;\R^m)$, $y\in \cY$.
\item
$\bar z(t)= P_{x,I,q}(t;w;y)$ results from applying $P_{x,I,q}$ to some $w, y$
as above. 
\end{enumerate}
\end{rem}

In fact, in case (i) one has 
  $\bar z_w(s)- \bar z_w(\xi_i)=0$. In case (ii)
 one has
  for $s\in J_i$\vspace*{-1.5mm}
\begin{align*}
|\bar z(s)-\bar z(\xi_i)|&=|\Phi_{x,I}(s;w;y)-\Phi_{x,I}(\xi_i;w;y)| 
\le  \Big|\int_{\xi_i}^s|a(s',w(s');y)|ds'\Big| 
\le \frac{A|I|}{2q}, \vspace*{-1.8mm}
\end{align*}
where we have again used \eref{LA}.
Finally for (iii)
we have for $s\in J_i$, by \eref{rhodiff} and \eref{LA},\vspace*{-1.5mm}
\begin{align*}
|\bar z(s)-\bar z(\xi_i)|&=|P_{x,I,q}(s;w;y)- P_{x,I,q}(\xi_i;w;y)|\le \sum_{k=1}^q|\rho_k(s)-\rho_k(\xi_i)|
|\bar a_k( w(\xi_k);y)|\le \frac{A|I|}{2q}.\vspace*{-1.5mm}
\end{align*}
which confirms the claim.\\


\newcommand{\ram}{\rho}

\noindent
{\em Compositional representation of $P_{x,q}$:}
Note that $P_{\cdot,q}$ 
can be written as a composition\vspace*{-1mm}
\be
\label{Pxcomp}
P_{x,q}(t;\bar z;y) = (g^{2,q}\circ g^{1,q})(t,x,\bar z;y).\vspace*{-1mm}
\ee
In slight abuse of notation we identify a piecewise constant  $\bar z$ with
the vector $
\bar z(\uxi):=(\bar z(\xi_1),\ldots, \bar z(\xi_q))\in \R^{qm}$ when writing
\be
\label{g1}
g^{1,q}: \big( t, x,\bar z;y )\mapsto \big(\rho_1(t),\ldots,\rho_q(t), x ,\bar a_1(\bar z(\xi_1);y),\ldots,
\bar a_q(\bar z (\xi_q);y )\big)\in \R^{(q+1)m +q },
\ee
and the bi-linear map
\be
\label{g2}
g^{2,q}: ( r_1,\ldots,r_q,  x,w^1,\ldots,w^{q})\mapsto x + \sum_{i=1}^{k}r_i w^i\in \R^m.
\ee
\subsubsection{Dimension-sparse approximation }

The approximation $P_{x,I,q}$ to $\Phi_{x,I}$ still involves the functions $\bar a_i(\bar z;y)$ which eventually need to be approximated by finitely parametrized expressions. Here we use the structural assumptions on 
the convection field.
 In case (A1) from \eref{acases}
 $a\in L_\infty([0,\Tin];\cA^{\gamma,\ss} )$ immediately implies that
$\bar a_i(\cdot;\cdot) = a(\xi_i;\cdot;\cdot)$ belong to 
$\cA^{\gamma,\ss}$, uniformly in $i=1,\ldots,q$, $q\in \N$. Hence, for each $i=1,\ldots,q$, $N\in \N$,
there is a composition 
$
\t A_{N,i}\in \Co_{N,\ss}$   such that \vspace*{-1.5mm}
\be
\label{tA}
\max_{i=1,\ldots,q} |\bar a_i(\bar z;y)- \t A_{N,i}(\bar z;y)|
\le
\gamma(N)^{-1}
\|a\|_{L_\infty(I;\cA^{\gamma,\ss})},\quad 
  \tripnorm{\t A_{N,i}}_{N,\ss }\le  \|a\|,
\vspace*{-1.5mm}
\ee
recall $\|a\|:=\|a\|_{L_\infty(I;\cA^{\gamma,\ss} )}$. In case (A2), the same conclusion holds, due to \eref{av}.

\begin{lemma}
\label{lem:compsparse}
We adhere to the definitions
$\bar a_i(\cdot;\cdot)=a(\xi_i,\cdot;\cdot)$  or $a_i= a_{J_i}$ when \eref{Lipt},
respectively \eref{av}, hold and
let \vspace*{-1.9mm}
\be
\label{Fqgen}
A_{x,I,q,N}(t,  \bar z;y ) := 
 x+ \sum_{i=1}^{q} \rho_i( t) \t A_{N,i}(\bar z;y).\vspace*{-1.9mm}
 \ee
Then, for either version of $\bar a_i$ one has \vspace*{-1.5mm}
\be
\label{a2genall}
\Big|\Phi_{x,I}(t;\bar z;y )- A_{x,I,q,N}(t,\bar z;y )\Big|
\le |I|\Big\{\frac{A}{q}+ \frac{\|a\|}{\gamma(N)}\Big\}. \vspace*{-1.5mm}
\ee
In particular, choosing
\be
\label{qetagen}
q= q(\tau):= \left\lceil \frac{2A|I| }{\tau}\right\rceil,\quad N=N(\tau)=
\left\lceil \gamma^{-1}(2|I|\|a\|/\tau)\right\rceil,
\ee
we have \vspace*{-1.5mm}
\be
\label{Phieta}
|\Phi_{x,I}(t;\bar z;y )- A_{x,I,q(\tau),N(\tau)}(t, \bar z;y )|\le \tau,
\ee
In what follows we write briefly  
$A_{x,I,\tau}:= A_{x,I,q(\tau),N(\tau)}$ with dimensionality vector $\Di_{x,I,\tau}$. 
\end{lemma}
{\bf Proof:}
For (A1) it follows  from \eref{Pxappr} and Lemma \ref{lem:quadrature}, (a),
 (see also \eref{Lip})
\begin{align}
\label{a2gen}
\Big|\Phi_{x,I}(t;\bar z;y )- A_{x,I,q,N}(t,\bar z;y )\Big|
&\le   \frac{A|I| }{q} + \sum_{i=1}^{q}\rho_i(t)|\bar a_i(\bar z,y)- \t A_{N,i}(\bar z,y )|
\nonumber\\
&\le  \frac{A|I| }{q}+\frac{|I|\|a\|}{\gamma(N)} = |I|\Big\{\frac{A}{q}+ \frac{\|a\|}{\gamma(N)}\Big\}.
\end{align}
The remainder of the assertion is an obvious consequence.\hfill $\Box$

\begin{rem}
\label{rem:A}
(a)
Suppose that $\Di_{i,q,N}$ is the dimensionality vector of $\t A_{N,i}$. Then,
a corresponding realization of $A_{x,I,q,N}$ results from parallelization 
of the $\t A_{N,i}$. One easily concludes from Remark \ref{rem:comp} that  the resulting
dimensionality vector $\Di_{x,I,q,N}$ of $A_{x,I,q ,N }$ is bounded by
\be
\label{AD}
\fN(\Di_{x,I,q,N})\le q\max_{i=1,\ldots,q}\fN\big(\Di_{i,I,q,N}\big)\le q N.
\ee
Thus, by \eref{qetagen}, one has for an absolute constant (depending only on $m$)
\be
\label{Aeta}
\fN(\Di_{x,I,\tau})\lsim q(\tau) \gamma^{-1}(\|a\|2|I|/\tau)\le \frac{ 2A|I|
\gamma^{-1}(2|I|\|a\|/\tau)}{\tau}\le \frac{ A  
\gamma^{-1}(1/\tau)}{\|a\|\tau},
\ee
where we have used \eref{re}.\\
(b) $A_{\cdot,I,q,N}$ has a compositional representation   analogous 
to \eref{Pxcomp}, obtained by replacing $a(\xi_i;\cdot,\cdot)$ by  
$\t A_{N,i}$. Since by assumption $\ss\ge m$, one  can  see from \eref{Pxcomp},
that $\ss$-dimension sparsity
of the $\t A_{N,i}$ is inherited by the mappings $A_{\cdot,I,q,N}$ and hence by
their compositions.  
\end{rem}

\subsubsection{Lipschitz continuity of pointwise compositions}

As a final prerequisit, to eventually control the stability of compositions of $A_{x,I,q,N}$,
 we need bounds for the Lipschitz constants of such compositions.  To that end, suppose that\vspace*{-1.6mm}
\be
\label{ANi}
\t A_{N,i} = (\t A_{N,i})^{n_i}\circ\cdots \circ (\t A_{N,i})^1,\quad i=1,\ldots,q,
\vspace*{-1.6mm}
\ee
where, by definition of $\Co_{N,\ss}$, each component $(\t A_{N,i})^{j}_\nu$,
$1\le \nu \le d_j$ depends for $j<n_i$ only on at most $\ss$ variables or is at most bilinear.

To proceed, recall also that the Lipschitz constants  of the 
factors in $\t A_{N,i}$ as well as the Lipschitz constants $L_{[n_i,j]}(\t A_{N,i})$ of the partial 
compositions $(\t A_{N,i})^{n_i}\circ\cdots \circ (\t A_{N,i})^j$ are controlled by 
\be
\label{tAicontrol}
\max_{1\le i\le q}\tripnorm{\t A_{N,i}}_{N,{ \Lip}} \le \|a\|:=\|a\|_{L_\infty(I;\cA^{\gamma,\ss}_\Lip)} .
\ee

\begin{lemma}
\label{lem:discrLip}
For any $q\in \N$, $k\in \N$, and $t\in I$, one has
\be
\label{Lipzy}
|A^k_{x,I,q,N}(t,\bar z;y)-A^k_{x,I,q,N}(t,\bar z';y')|\le \frac{(\|a\||I|)^k}{k!}
\max\big\{
|y-y'| ,\|\bar z-\bar z'\|_{L_\infty(I;\R^m)}\big\},
\ee
for $ x\in D, \bar z(\uxi) ,\,\bar z(\uxi)\in \R^{mq}$.
Similarly, when $\bar z = \bar z_x$, $\bar z'=\bar z_{x'}$, one has for
all $x, x'\in D$, 
\be
\label{Lipxy}
|A^k_{x,I,q,N}(t,\bar z_x;y)-A^k_{x',I,q,N}(t',\bar z_{x'};y')|\le \|a\||t-t'|+ \max\big\{|y-y'| ,|x-x'|\big\} e^{\|a\||I|}  .
\ee
Finally, the $A^k_{\cdot,I,q,N}$ belong to $\Co_{CkqN.\ss,\Lip}$, where $C$ is a fixed 
constant.
\end{lemma}
The reasoning is analogous to the proof of Lemma \ref{lem:Lipprop} based on the
smoothing effect of multiple integration, here in terms of multiple summation.
The proof is therefore deferred to Appendinx B.
\vspace*{-2.5mm}
\subsection{Proof of Theorem \ref{thm:char1}
 }\label{ssec:char1} 
\noindent
{\em Step 1 - construction of an $\e$-accurate pointwise composition:}
Given the $\e$-accurate approximation of the characteristic field by compositions
of global operators $Z_\e$ from \eref{Ze}, we construct now a {\em pointwise compositional} counterpart. 
 Specifically, we define approximations $\wt\Psi_{[k,j]}$, $\wt\Psi_{\ue^k}=
 \wt\Psi_{[k,0]}$ to the (global) counterparts  $ \Psi_{[k,j]}$, $\Psi_{\ue^k}$.
 from \eref{psinew}, \eref{specific}. We adhere to the meaning of $\mu_k,\eta_k$ from \eref{muk}, and replace $\Phi_{x,I}$ by $A_{x,I,\tau}=A_{x,I,N(\tau),q(\tau)}$. Precisely, let
 for $j<k$, $w\in\R^m$, and for a new vector of tolerances\vspace*{-1.5mm}
 $$
 \utau=(\tau_1,\ldots,\utau_K),\quad \mbox{with sections}\quad \utau^k=(\tau_1,\ldots,\tau_k),\quad 1\le k\le K,\vspace*{-1.5mm}
 $$
 yet to be chosen. We define for $t\in I_k$\vspace*{-1.5mm}
\be
\label{psinew}
\begin{array}{ll}
\wt\Psi_{[k,j]}(t,w;y):= A^{\mu_k}_{w_{k-1,j},I,\tau_k}(t,t_{k-1};\bar z_{w_{k-1,j}};y)   ,& w_{k-1,j}:= \wt\Psi_{[k-1,j]}(t_{k-1};w;y), \\
\wt\Psi_{[j+1,j]}(t;w;y) := A^{\mu_{j+1}}_{w,I,\tau_{j+1}}(t,t_j;\bar z_w;y).& 
\end{array}\vspace*{-1.5mm}
\ee
We denote as before
\be
\label{specific}
\wt\Psi_{\utau^k}(t,x;y):= \wt\Psi_{[k,0]}(t,x;y).
\ee
We  choose $\tau_k=\tau $ all equal so that for $\eta_k$ given by \eref{ej}
\be
\label{tauk}
|\Phi^{\mu_k}_{w,I_k}(t,t_{k-1};\bar z_w;y)- A^{\mu_k}_{w,I_k,\tau_k}(t,t_{k-1};\bar z_w;y)|
\le \eta_k,\quad k\le K.
\ee
Since, by Lemma \ref{lem:Lipprop}, \eref{globlip}, the Lipschitz constants
$L_{[k,j]}$ of $k-j$ partial compositions of $\Phi_{w,I}$ in Remark \ref{rem:error}
are bounded by $\frac{2^{-(k-j)}}{(k-j)!}$, we conclude\vspace*{-2.4mm}
\be
\label{taucond}
|\Phi^{\mu_k}_{w,I_k}(t,T_{k-1};\bar z_w;y)- A^{\mu_k}_{w,I_k,\tau_k}(t,T_{k-1};\bar z_w;y)|
\le \tau +\tau \sum_{j=1}^{\mu_k-1}\frac{2^{-\mu_k-j}}{(\mu_k-j)!}\le \tau e^{1/2}.
\vspace*{-1.5mm}
\ee
On account of \eref{ej},  choosing $\tau=\tau(\e)$ such that $\tau e^{1/2}
\le \eta(\e)$ from \eref{ej}, i.e.,
\be
\label{tauke}
\tau_k(\e)=\tau(\e)= e^{-1/2}\eta (\e)=  {\e e^{- K/2}}(1-e^{-1/2}),\quad k=1,\ldots,K,
\ee
yields via the same reasoning as in \eref{uniform} \vspace*{-1.5mm}
\be
\label{Aeps}
\|\Psi_{\ue^k(\e)}-\wt\Psi_{\utau^k(\e)}\|_{L_\infty(\Hom(I_k))}\le \e.\vspace*{-1.5mm}
\ee
%

In summary, we obtain as before\vspace*{-1.9mm}
\be
\label{Ze}
\|z- \wt Z_\e\|_{L_\infty(\Hom;\R^m)}\le 2\e, \quad\mbox{where}\quad \wt Z_\e(t;x;y):= \sum_{k=1}^K\chi_{I_k}(t)\wt\Psi_{\utau^k(\e)}(t;x;y).\vspace*{-1.5mm}
\ee

\noindent
{\em Step 2 - Complexity of $\wt Z_\e$:} 
It follows   from Remark \ref{rem:comp} that\vspace*{-1.9mm}
\be
\label{estN1}
\fN(\wt Z_\e)\le \sum_{k=1}^K\mu_k(\e)\fN(A_{\cdot;I_k,\tau_k(\e)}).
\ee
On account of \eref{ej} and \eref{muk} we have (recall $|I_k|=\Tin/K$ and 
$K/2 = \Tin   \|a\| $ by \eref{re})\vspace*{-1.5mm}
\be
\label{mukeps}
\mu_k(\e) =  \Big|\log_2 \Big(\frac{e^{K/2}}{2(e^{1/2}-1)\e}  \Big)\Big| \eqsim  \Big|\log_2 \Big(\frac{ e^{\|a\|\Tin}}{ \e}\Big)\Big|.
\ee
Furthermore,
\eref{tauke} in conjunction with Remark \ref{rem:A} yields (for the range of 
$\gamma$ under consideration, see \eref{ginverse})\vspace*{-1.5mm}
$$
\fN(A_{\cdot;I_k,\tau_k(\e)})\eqsim \frac{Ae^{K/2}}{\|a\|\e} \gamma^{-1}
\Big(\frac{e^{K/2}}{(1-e^{-1/2})\e}\Big)
\eqsim \frac{Ae^{\|a\|\Tin}}{\|a\|\e} 
\gamma^{-1}\Big(\frac{ e^{\|a\|\Tin}}{\e}\Big)\vspace*{-1.5mm}
$$ 
Substituting this into \eref{estN1}, yields\vspace*{-1.5mm}
\begin{align}
\label{Nest2}
\fN(\wt Z_\e)&\eqsim K  \Big|\log_2 \Big(\frac{ e^{\|a\|\Tin}}{\e}\Big)\Big|
\frac{A e^{\|a\|\Tin}}{\|a\|\e} 
\gamma^{-1}\Big(\frac{  e^{\|a\|\Tin}}{\e}\Big)\nonumber\\
& \eqsim A \Tin\Big|\log_2 \Big(\frac{ e^{\|a\|\Tin}}{ \e}\Big)\Big|
\frac{e^{\|a\|\Tin}}{ \e} 
\gamma^{-1}\Big(\frac{  e^{\|a\|\Tin}}{\e}\Big).\vspace*{-1.9mm}
\end{align}
By \eref{ginverse}, we obtain \vspace*{-1.5mm}
\be
\label{recall}
\gamma^{-1}\big(  e^{\|a\|\Tin}\e^{-1}\big) \eqsim
\left\{
\begin{array}{ll}
 \Ca^{-1/\alpha}e^{ \|a\|\Tin/ \alpha}\e^{-\frac{1}{\alpha}},& \mbox{for $\gamma \sim$ (alg)},\\
\frac{1}{\alpha}  \ln \Big(\frac{e^{\|a\|\Tin}}{\Ce\e}\Big) ,& \mbox{for $\gamma \sim$ (exp)}.
\end{array}\right.\vspace*{-1.5mm}
\ee
 We conclude that
 for $\gamma \sim$ (alg), (see \eref{assa}) 
\begin{align}
\label{Nestalg}
 \fN(\wt Z_\e)&\eqsim A \Ca^{-1/\alpha} 
\Tin   \log_2 \Big(\frac{e^{\|a\|\Tin}}{\e }\Big)
  \Big( \frac{e^{ \Tin\|a\| }}{\e}\Big)^{\frac{\alpha +1}{\alpha}} :=\phi_{\rm alg}\big(
  e^{\|a\|\Tin}/\e).
\end{align}
For
  exponential growth $\gamma \sim$ (exp) (since $ \|a\|\ge L\ge 1$, \eref{LA}) 
  we obtain 
\begin{align}
\label{Nestexp}
 \fN(\wt Z_\e)&\eqsim \frac{1}{\alpha}A \Tin\Big(\log_2 \Big(\frac{e^{\|a\|\Tin}}{\e }\Big)\Big)^2
  \Big( \frac{e^{ \Tin\|a\| }}{\e}\Big):=\phi_{\rm exp}\big(
  e^{\|a\|\Tin}/\e)
\end{align}
provided that $\Ce\gsim 1$.
\begin{rem} 
\label{rem:tZe}
On account of Remark \ref{rem:A}, (b) and Lemma 
\ref{lem:discrLip}, we conclude that $\wt Z_\e\in
\Co_{CN_\e,\ss }$ for some uniform constant $C$ and 
$\tripnorm{\wt Z_\e}_{N_\e,\ss } \le e^{\|a\|\Tin}$ for $N_\e$, defined by
the respective right hand sides in 
\eref{Nestalg}, \eref{Nestexp}.
\end{rem}

\noindent
{\em Step 3 - Convergence rates:} To determine the convergence rates, corresponding to
\eref{Ze}, we apply Remark \ref{rem:order} and Lemma \ref{lem:algrowth} to 
$N_\e = \phi_{\rm alg/exp}\Big(e^{\|a\|\Tin}/\e\Big)$ from \eref{Nestalg} and \eref{Nestexp}.  \eref{tgamma} follows then by straightforward calculations.

The following statement follows now from Remark \ref{rem:order} and the 
above observations.

\begin{rem}
\label{rem:summary}
In summary we shown that for each $N\in \N$ there exists  an $\ss$-dimension-sparse  compositional representation
$\wt Z_N$ satisfying
\be
\label{sumrate}
\|z-\wt Z_N\|_{L_\infty(\Hom;\R^m)}\lsim e^{\|a\|\Tin} \wt\gamma(N)^{-1},\quad N\in \N.
\ee
with $\wt\gamma$ is from \eref{tgamma}. This      confirms
\eref{zrate}. 
\end{rem} 

\noindent
{\em Step 4 - Stability of $\wt Z_\e$
:} It follows from Lemma \ref{lem:discrLip},
\eref{Lipxy},
that the Lipschitz constants of partial compositions of the approximations
$\wt\Psi_{\utau^k}$ from \eref{specific} remain uniformly bounded by
$\|a\| + e^{\|a\||I_k|}$. 
Specifically,\vspace*{-1.5mm}
$$
|\wt Z_\e(t,x,y)-\wt Z_\e(t,x'y')|\le \max\{|x-x'|,|y-y'|\}e^{\|a\|\Tin},\quad
 t \in [0,\Tin].\vspace*{-1.5mm}
$$
i.e., \vspace*{-1.5mm}
\be
\label{tripZ}
\tripnorm{\wt Z_\e(t,\cdot,\cdot)}_{N_\e,\ss }\le    e^{\|a\|\Tin},\quad 
t\in [0,\Tin].\vspace*{-1mm}
\ee
Moreover,  the growth functions $\hat\gamma$ satisfy $\hat\gamma(N_\e)\eqsim \e$. 
By Remark \ref{rem:order}, this finishes 
  the proof of Theorem \ref{thm:char1}.\hfill $\Box$.
\vspace*{-2.9mm}
\subsection{Proof of Theorem \ref{thm:dimsparse}}\label{ssec:implant} 
\newcommand{\wh}{\widehat}
\vspace*{-1.5mm}
The main step  is to invoke   Lemma \ref{lem:import} and Proposition \ref{prop:Lip}  to approximate each $\t A_{N,i}(\bar z,y)
\in \Co_{N,\ss}$
in \eref{Fqgen} by a DNN.   In essence we follow the same steps as in the 
preceding section.\\[1.8mm]
\noindent
{\em Approximation of data - accuracy:}
Recall from \eref{ANi} that $\t A_{N,i}$ has a
compositional representation
$\t A_{N,i} = (\t A_{N,i})^{n_i}\circ\cdots \circ (\t A_{N,i})^1$ of some depth $n_i\le \fN(\t A_{N,i})$.

For a given tolerance $\delta>0$, we construct next a 
 network $\cN_{i,\delta,N,q}(\bar z,y)$ approximating $\t A_{N,i}$ with accuracy $\|a\|\delta$. Specifically, we approximate each component 
$(\t A_{N,i})^j_\nu$, $\nu =1,\ldots,d_j$, in the $j$th composition factor 
by a $\delta/n_i$-accurate neural network  $(\cN)^j_\nu$.
 One infers from Lemma \ref{lem:import}, \eref{NapprLip}, \eref{tAicontrol},
that (since $\|a\|\ge 1$)\vspace*{-2.5mm}
\begin{align}
\label{delterr1}
\|\t A_{N,i} - \cN_{i,\delta,N,q}\|_{L_\infty(\R^m\times \cY;\R^m)}
&\le n_i^{-1}\Big\{ \delta  + \sum_{r=1}^{n_i -1}\delta L_{[n_i,r+1]}(\t A_{N,i})\Big\} \nonumber\\
&\le n_i^{-1}\Big\{\delta  +  \|a\|\sum_{r=1}^{n_i -1}\delta\Big\}
\le  \|a\|\delta, \quad i=1,\ldots,q.\vspace*{-1.5mm}
\end{align}
On the other hand, we invoke \eref{Nsharp} in Lemma \ref{lem:import} 
to conclude  \vspace*{-1mm}
\begin{align}
\label{deltasharp2}
\# \cN_{i,\delta,N,q}&\le C_\ss\|a\|^\ss N n_i^\ss\delta^{-\ss}\big||\log_2\delta|+\log_2 n_i\big|,\quad i=1,\ldots,q, \vspace*{-1.5mm}
\end{align}
where the constant  $ C_\ss$  depends only
on  $\ss$.
Now we   define in analogy to \eref{Fqgen}\vspace*{-1.5mm} 
\be
\label{Nqgen}
\cN_{x,I,q,N,\delta}(t,z;y):= x+ \sum_{i=1}^q \rho_i(t) \cN_{i,\delta,N,q}(z,y),\quad 
t\in I ,\vspace*{-1.5mm}
\ee
where we recall that for $I=[\uT,\oT]$,    $J_i=J_i(I,q)
:= \uT + \big[\frac{(i-1)|I|}q,\frac{i|I|}q\big]$ and $\rho_i(t):=\int_{\uT}^t\chi_{J_i}(s)ds$.
Taking\vspace*{-1mm}
\be
\label{Netaq}
 \cN_{x,I,\tau}(t ;z;y):=\cN_{x,I,q(\tau/2),N(\tau/2),\delta(\tau)}(t;z;y)
\ee
with $q(\tau/2),N(\tau/2)$, defined according to \eref{qetagen},
and\vspace*{-2mm}
\be
\label{qetagen2}
\delta = \delta(\tau) := \frac{\tau}{2\|a\|},\vspace*{-1.5mm}
\ee
we conclude that $\|A_{x,I,\tau}- \cN_{x,I,\tau}\|_{L_\infty(\Hom(I);R^m)}\le \tau/2$. Hence, by \eref{Phieta},\vspace*{-1.5mm}
\be
\label{target}
\Big|\Phi_{x,I}(z,t;y) - \cN_{x,I,\tau}(z,t;y)\Big| \le \tau.\vspace*{-1.5mm}
\ee

To estimate the complexity of $\cN_{x,I,\tau}$ we use \eref{deltasharp2} and bound the depths $n_i$ by $N(\tau/2) \sim 
\gamma^{-1}\big(4|I|\|a\|/\tau)= \gamma^{-1}\big(2/\tau)$.
Arguing as in \eref{AD} and recalling \eref{qetagen}, 
 we then have (since $q(\tau)= 2A|I|/\tau =\frac{A}{\|a\|\tau}$)\vspace*{-1.5mm}
\begin{align}
\label{AD2}
\#\cN_{x,I,\tau} &\eqsim q(\tau/2)\max_{i=1,\ldots,q(\tau/2)}\#\cN_{i,\delta(\tau/2),N
(\tau/2)}\nonumber\\
&\eqsim \frac{2A}{\|a\|\tau} \|a\|^\ss \gamma^{-1}(2/\tau)^{1+\ss} \delta(\tau)^{-\ss} \big||\log_2\delta(\tau)| +\log_2 N(\tau/2)\big| \nonumber\\
&\eqsim    A 2^{\ss}\|a\|^{2\ss-1}    \tau^{-(1+\ss)}\gamma^{-1}(2/\tau)^{1+\ss} 
\Big|\log_2 \frac{2\|a\|\gamma^{-1}(2/\tau)}{\tau}\Big|.\vspace*{-1.5mm}
\end{align}
For algebraic growth  
$\gamma(r)=  \Ca r^\alpha$ we obtain 
$
\gamma^{-1}\big(2/\tau\big)
= (2/\Ca)^{1/\alpha}\tau^{-1/\alpha}$, and $\log_2 \gamma^{-1}(2/\tau)\sim \alpha^{-1}\log_2 (2/\Ca \tau) 
$, so that  $\Big|\log_2 \frac{2\|a\|\gamma^{-1}(2/\tau)}{\tau}\Big|=
\frac{\alpha+1}{\alpha}\Big|\log_2\frac{2}{\tau}\Big(\frac{\|a\|}{\Ca}\Big)^{\frac{\alpha}{\alpha+1}} \Big|$. For $\gamma \sim $(exp), we have 
$\Big|\log_2 \frac{2\|a\|\gamma^{-1}(2/\tau)}{\tau}\Big|=\log_2\frac{2\|a\|\ln(2/\Ce\tau)}{\alpha\tau}$.
Thus, \vspace*{-1.5mm}
\be
\label{AD4}
\#\cN_{x,I,\tau}\eqsim A 2^{\ss}\|a\|^{2\ss-1 }
\left\{
\begin{array}{ll}
\Ca^{-1/\alpha}\tau^{-\frac{(1+\alpha)(1+\ss)}{\alpha}}|\log_2\tau|,&
\mbox{when $\gamma \sim$ (alg)},\\
\alpha^{-(\ss+1)}\tau^{-(\ss +1)}
|\log_2 \tau|^{2+\ss},& \mbox{when $\gamma \sim$ (exp)},
\end{array}\right.\vspace*{-1.5mm}
\ee
where we recall that $\Ce\gsim 1$, accepting a logarithmic dependence 
of the proportionality constant on $\|a\|$ (or assume that $\tau\le \tau_0(\|a\|)$).

We can now define $\wh\Psi_{[k,j]}(t,w;y)$ in analogy to \eref{psinew}
with $A_{w,I,\tau_k}(t ;\bar z_w;y)$ replaced by $\cN_{w,I,\tau_k}(t ;\bar z_w;y)$ and likewise $\wh\Psi_{\utau^k}= \wh\Psi_{[k,0]}$ in analogy to
\eref{specific}. With the same tolerances $\tau_k(\e)=\tau(\e)$, given by \eref{tauke}.
The same reasoning as in \S~\ref{ssec:char1} yields (see \eref{Aeps} \vspace*{-1.5mm} 
\be
\label{Neps}
\|\Psi_{\ue^k(\e)}-\wh\Psi_{\utau^k(\e)}\|_{L_\infty(\Hom(I_k))}\le \e, \quad k=1,\ldots,K,\vspace*{-1.5mm}
\ee
and hence\vspace*{-1.5mm}
\be
\label{Neest}
\|z- \cN_\e\|_{L_\infty(\Hom;\R^m)}\le 2\e,\quad \cN_\e(t,x;y):=\sum_{k=1}^K\chi_{I_k}(t)\wh\Psi_{\utau^k(\e)}(t,x;y).
\ee

\noindent
{\em Complexity:} 
 It remains to bound   $\# \cN_\e$. Invoking Remark \ref{rem:comp} as before,   one obtains from \eref{AD4}
 with $\mu_k(\e)$ from \eref{mukeps} (see also \eref{muk}, \eref{ej})
\begin{align}
\label{Nefin}
\# \cN_\e &\eqsim \sum_{k=1}^K\mu_k(\e)\cN_{\cdot,I_k,\tau_k}\eqsim 
K   \log_2\Big(\frac{e^{\|a\|\Tin}}{2\|a\|\e}\Big) \#\cN_{\cdot,I,\tau(\e)}\nonumber\\
& \eqsim A2^{\ss}\|a\|^{2\ss}  \Tin  \log_2\Big(\frac{e^{\|a\|\Tin}}{2\|a\|\e}\Big)   
\left\{
\begin{array}{ll}
\Ca^{-1/\alpha}\tau(\e)^{-\frac{(1+\alpha)(1+\ss)}{\alpha}}|\log_2\tau(\e)|,&
\mbox{$\gamma \sim$ (alg)},\\
\alpha^{-(\ss+1)}\tau(\e)^{-(\ss +1)}|\log_2 \tau(\e)|^{2+\ss},& \mbox{$\gamma \sim$ (exp)}
\end{array}\right.
\end{align}
By  \eref{tauke}, \vspace*{-1.5mm} 
\begin{align*}
\Ca^{-1/\alpha}\tau(\e)^{-\frac{(1+\alpha)(1+\ss)}{\alpha}}|\log_2\tau(\e)|
&\eqsim \Ca^{-1/\alpha}\Big(\frac{e^{\|a\|\Tin}}{ \e}\Big)^{\frac{(1+\alpha)(1+\ss)}{\alpha}} \log_2\Big(\frac{e^{\|a\|\Tin}}{ \e}\Big),\vspace*{-1.5mm}
\end{align*}
while\vspace*{-1.5mm}
$$
\alpha^{-(\ss+1)}\tau(\e)^{-(\ss +1)}|\log_2 \tau(\e)|^{2+\ss}\eqsim
\alpha^{-(\ss+1)}\Big(\frac{e^{\|a\|\Tin}}{\e}\Big)^{-(1+\ss)} \Big|\log_2\Big(\frac{e^{\|a\|\Tin}}{\e}\Big)\Big|^{2+\ss}.\vspace*{-1.5mm}
$$
Inserting these estimates into \eref{Nefin}, confirms \eref{DNN2}.

The convergence rates stated in Corollary
\ref{cor:dimsparse} follow now in the same way as before by applying 
Remark \ref{rem:order} and Lemma \ref{lem:algrowth} to the bounds on
$\fN(\cN_\e)$, given in \eref{DNN2}.


Regarding the stability of the networks $\cN_\e$,
there is a principal obstacle related to the fact that 
the precise dimension-vectors of the compositions 
$
\t A_{N,i} = \t A_{N,i}^{n_i}\circ\cdots \circ \t A_{N,i}^1
$
especially their depths $n_i$ are not known. Although, the Lipschitz constants
of the implanted networks in each factor $\t A_{N,i}^{\nu}$ are controlled
by $\|a\|$ it is not clear whether the Lipschitz constants of their compositions
also remain controlled by $\|a\|$. Such network approximations exist by Proposition \ref{prop:Lip} but need no longer be $\ss$-dimension sparse.  
So, the only guaranteed general 
bound for the Lipschitz constants of the partial compositions
  is 
in view of  Proposition \ref{prop:Lip}  
$L_{[n_i,j+1]}(\cN_{i,\delta,N,q})\le (c_3(1+A)\|a\|)^{n_i-j}$. Hence, 
\be
\label{LN}
 L_{N,\delta} := 
  \max_{j< n_i}L_{[n_i,j+1]}(\cN_{i,\delta,N,q}) \le (c_3(1+A)\|a\|)^{N},
\ee
where we have applied the (perhaps too pessimistic) bound $n_i\le N$.
If on the other hand, the depths $n_i$ remain uniformly bounded by $\bar n$, say,
one obtains
a uniform Lipschitz-bound $ L_{N,\delta}\le (c_3(1+A)\|a\|)^{\bar n}$. 
\begin{lemma}
\label{lem:Liptxya}
The network approximations $\cN_{\cdot,I,q,N,\delta}$ from \eref{Nqgen} 
have the
following Lipschitz continuity properties:
 for $ (t,x,y),\,(t',x',y')\in \Hom(I)$, $\bar z ,\bar z\in L_\infty(I;\R^m)$:
\be
\label{Lipzy0}
|\cN^\ell_{x,I,q,N,\delta}(t,\bar z;y)-\cN^\ell_{x,I,q,N,\delta}(t,\bar z';y')|\le \frac{(L_{N,\delta}|I|)^{\ell }}{\ell !}\big\{
|y-y'| ,\|\bar z-\bar z'\|_{L_\infty(I;\R^m)}\big\}, 
\ee
where $L_{N,\delta}$ is given by \eref{LN}.
Similarly, when $\bar z = \bar z_x$, $\bar z'=\bar z_{x'}$, one has 
\begin{align}
\label{Lipxy0}
&|\cN^\ell_{x,I,q,N,\delta}(t,\bar z_x;y)-\cN^\ell_{x',I,q,N,\delta}(t',\bar z_{x'}, ;y')|\nonumber\\
&\qquad\le (1+\delta)\|a\||t-t'|+\max\big\{|y-y'|, |x-x'|\big\} e^{L_{N,\delta}|I|} .
\end{align}
\end{lemma}
\noindent
{\bf Proof:}
Recall from    \eref{Nqgen} that $\cN_{x,I,q,N,\delta}(t,\bar z;y):= x+  
\sum_{i=1}^{q}\rho_i(t) \cN_{i,\delta,N,q}( \bar z,y)$.
Then, we have for $i=1,\ldots,q$,
\be
\label{boundLam}
 |\cN_{i,\delta,N,q}(\bar z,y)- \cN_{i,\delta,N,q}(\bar z',y')|\le 
 |\cN_{i,\delta,N,q}|_{\Lip_1(\R^m\times \cY)} 
 \max\{\|\bar z-\bar z'\|_{L_\infty(I;\R^m)},|y-y'|\}.
 \ee
By the comments preceding the lemma, we obtain
\begin{align}
\label{Lip2}
&\Big|\cN_{x,I,q,N,\delta}(t,\bar z;y)- \cN_{x,I,q,N,\delta}(t,\bar z';y')\Big|
\le 
  \sum_{i=1}^{q}\rho_{ i}(t)L_{N,\delta}  \max\{\|\bar z-\bar z'\|_{L_\infty(I;\R^m)},|y-y'|\}.
\end{align}
Hence, we are in the same situation as
in \eref{thesame}. 
Finally, 
\vspace*{-1.5mm}
\begin{align*}
&\Big|\cN_{x,I,q,N,\delta}(t,\bar z;y)- \cN_{x,I,q,N,\delta}(t',\bar z;y)\Big|
\\
&\qquad\le \|a\||t-t'|
 +\sum_{i=1}^q|\rho_i(t)-\rho_i(t')||\bar a_i(\xi_i;y)-\cN_{i,\delta,N,q}(\bar z(\xi_i);y)| 
\le (1+\delta)\|a\||t-t'|,\vspace*{-1.5mm}
\end{align*}
where we have used \eref{delterr1}.
Therefore, the claim follows by the same arguments as used
in the proof of Lemma \ref{lem:discrLip}.\hfill $\Box$\\

Regarding Remark \ref{rem:stab},  recall that $N(\e)\sim \gamma^{-1}(2/\tau(\e))$, where $\tau(\e)\eqsim \e e^{-\|a\|\Tin}$, $\delta(\e) = \tau(\e)/2\|a\|$ so that
by \eref{LN},\vspace*{-1.8mm}
\be
\label{Leps}
L_{N(\e),\delta(e)} \eqsim \left\{
\begin{array}{ll}
(c_3(1+A)\|a\|)^{(2/\Ca)^{\frac{1}{\alpha}}\e^{-1/\alpha} e^{\|a\|\Tin/\alpha}},&
\mbox{in case $\gamma\sim$ (alg)},\\
(c_3(1+A)\|a\|)^{ \frac{1}{\alpha}\big(|\ln \e|+\|a\|\Tin\big)}
,&
\mbox{in case $\gamma\sim$ (exp)}.
\end{array}\right.\vspace*{-1.5mm}
\ee
This confirms Remark \ref{rem:stab}.
\hfill $\Box$\vspace*{-1.5mm}
\subsection{Proof of Corollary \ref{cor:special} and Theorem \ref{thm:parconv}}\label{ssec:parconv2}
 
We first prove Corollary \ref{cor:special}. To apply Theorem \ref{thm:dimsparse}
and Corollary \ref{cor:dimsparse}
note first that condition \eref{av} is applicable, i.e., version (A2) can be used.
In fact, $a_j(t;\cdot )=\omega_j a^\circ_j(t,\cdot)\in \Lip_1(\R^m )$, uniformly in $t\in[0,\Tin]$, immediately implies that $\bar a_{j,i} = a_{j,J_i}= \omega_j
a^\circ_{j,J_i}$ belongs
to $\Lip_1(\R^m )$, for $i=1,\ldots,q$, with the same Lipschitz constants 
$\Lambda\omega_j$ from \eref{La}.

Next we recall from Remark \ref{rem:a} that $a\in \Co_{N_a,m }$ with
$N_a=1+ d_y(1+m^2)$.
 Hence, the simplest compositional approximations $A_N(t,\cdot)$ to $a$ is
 \vspace*{-1.5mm}
$$
A_N(t,x,y)= \left\{
\begin{array}{ll}
0, & N< N_a,\\
a(t;x;y), &N\ge N_a.
\end{array}
\right.\vspace*{-1.5mm}
$$
In view of Remark \ref{rem:a}, \eref{tripa}, one obtains  for $\gamma(r)= \Ce e^{\alpha r}$ and all $t\in [0,\Tin]$, 
\begin{align*}
\|a\|_{L_\infty(\In;\cA^{\gamma,m}_\Lip)}\le & \max_{N\in \N}\gamma (N)\Big\{\|a(t)- A_N(t)\|_{L_\infty(\R^m\times \cY;\R^m)} +\gamma(N)^{-1}
\tripnorm{A_N(t)}_{N,m }\Big\} \\
&\le  
\left\{
\begin{array}{ll}
A \Ce e^{\alpha N}+ (A+\Lambda |\uomega|_1),& N<N_a,\\
A+\Lambda |\uomega|_1 ,& N\ge N_a.
\end{array}\right.
\end{align*}
Taking   
$
\,\alpha := N_a^{-1}, \quad \Ce:= 1,\,
$
yields $\|a\|\le \wh L:= 2A+  \Lambda  |\uomega|_1$.
 By \eref{Coa}, one has $\alpha \eqsim d_y^{-1}$ (with $m$-dependent proportionality).
Theorem \ref{thm:char1} yields then $z\in L_\infty(\In;\cA^{\wt\gamma,m}_\Lip)$ with 
$$
\wt\gamma(r)\eqsim \frac{ r}{d_yA\Tin}\Big|\log_2 \frac{r}{d_yA\Tin}\Big|^{-2},
\quad \|z\|_{L_\infty([0,\Tin];\cA^{\wt\gamma,m}_\Lip)}\le e^{(2A+\Lambda |\uomega|_1)\Tin},
$$
and hence \eref{aapprox}. Now \eref{zratesp} and the expression for $F$ follow from Corollary \ref{cor:dimsparse}. \hfill $\Box$\\

We now turn to the proof of Theorem \ref{thm:parconv} approximating
  $\Phi_{x,I}$ 
  in a first step by \vspace*{-1.5mm} 
\be
\label{Px}
P_{x,I,q}(t,\bar z;y) := x+ \sum_{j=1}^{d_y} y_j\sum_{i=1}^{q}\rho_i(t)a_{j,J_i}(\bar z(\xi_i)),\quad a_{j,J_i }(\cdot):= \frac{\omega_j}{|J_i|}\int_{J_i}a^\circ_j(s;\cdot)ds,\vspace*{-1.5mm}
\ee
recalling that
$J_i=J_i(I,q)= \uT+ \big[\frac{(i-1)|I|}{q}, \frac{i|I|}{q}\Big]$ and
$\rho_i(t)= \int_{J_i}\chi_{s\le t}(s) ds$.
For $L$ from \eref{aL} and $\bar z$ as in Lemma \ref{lem:PhiP},
we infer from \eref{Pxappr} that $|\Phi_{x,I}(t,\bar z;y)- P_{x,I,q}(t,\bar z;y)|
\le \frac{A|I| }{q}$. 
Invoking Proposition 
\ref{prop:Lip}, we approximate the low-dimensional functions $a^\circ_{j,J_i}(\cdot)$
by finitely parametrized functions such as neural networks.   Specifically, there exist networks (suppressing the reference to $I$)
$\cN_{j,i,\delta}$ of depth $\lsim \log_2\delta^{-1}$ such that
for $j=1,\ldots,d_y$, $i=1,\ldots,q$,\vspace*{-1.5mm}
\be
\label{locappr}
\|a^\circ_{j,J_i} - \cN_{j,i,\delta}\|_{L_\infty(\R^m;\R^m)}\le \delta,\quad
 \#\cN_{j,i,\delta}\lsim \Lambda^m\delta^{-m}|\log_2\delta| .\vspace*{-1.5mm}
\ee 
Then\vspace*{-1.5mm}
\be
\label{Nx}
\cN_{x,I,q,\delta}(t,\bar z;y):= x+ \sum_{j=1}^{d_y} y_j\sum_{i=1}^{q}\rho_i(t)\omega_j\cN_{j,i,\delta}(\bar z(\xi_i)),\vspace*{-1.5mm}
\ee
is indeed an $m$-dimension-sparse neural network.
To that end, we keep viewing $t$ as a parameter and the input-variables $x,y$ are   passed across layers,
formally in a ``skip-connection'' format. Thus, formally we have\vspace*{-1.5mm}
\be
\label{Ncompos}
\cN_{x,I,q,\delta}(t,\bar z,y) = (G_{2}\circ G_{1})(t,x,\bar z,y).\vspace*{-1.5mm}
\ee
For better readability the following representation groups variables in a formally incorrect way and should be viewed as a $t$-dependent mapping into $\R^{1+m+d_y
+ m q d_y}$\vspace*{-1.5mm}
\be
\label{G1}
G_{1}: (t,x,\bar z,y) \mapsto 
\left(
\begin{array}{c}
x,y,\rho_1(t),\ldots,\rho_q(t)\\
 \omega_1\cN_{1,1,\delta }(\bar z),\ldots  \omega_{d_y}\cN_{d_y,1,\delta }(\bar z),\\
 \vdots \\
\omega_1 \cN_{1,q  ,\delta }(\bar z),\ldots, \omega_{d_y} \cN_{d_y,q ,\delta }(\bar z)
\end{array}\right) 
 \in \R^{q+m+d_y
+ m q d_y}, \vspace*{-1.5mm} 
\ee
which is obviously $m$-dimension sparse.
Hence $G_1$ itself is a neural network whose depth is bounded by $\log_2\delta^{-1}$. 
 The tri-linear factor $G_2$ reads then\vspace*{-1.5mm}
\be
\label{G2}
G_2: (x,y,r_1,\ldots, r_q,\zeta^{1,1},\ldots, 
\zeta^{d_y,q }) \mapsto \Big( x+ \sum_{i=1}^{q}r_i
\sum_{j=1}^{d_y}y_j \zeta^{j,i}\Big)\in \R^{m}.\vspace*{-1.5mm}
\ee
Assessing the accuracy of $\cN_{x,I,q,\delta}$ follows in essence the same lines
as before.
In view of \eref{Pxappr} and \eref{locappr},
\begin{align}
\label{Nerror}
&|\Phi_{x,I}(t,\bar z;y)- \cN_{x,I,q,\delta}(t,\bar z;y)|\le |\Phi_{x,I}(t,\bar z;y)-
P_{x,I,q}(t,\bar z;y)|\nonumber\\
&\qquad + \sum_{j=1}^{d_y} |y_j|\sum_{i=1}^{q}\rho_i(t)\omega_j|\cN_{j,i,\delta}(\bar z(\xi_i))-a^\circ_{j,J_i}(\bar z(\xi_i))| 
\le \frac{A|I| }{q} + |\uomega|_1|I|\delta.
\end{align}
Thus, given any target tolerance $\tau>0$, choosing
\be
\label{qeta}
q(\tau)= \left\lceil \frac{2A|I| }{\tau}\right\rceil,\quad \delta(\tau)= \frac{\tau}{2|I||\uomega|_1},
\ee
and abbreviating $\cN_{x,I,\tau}:= \cN_{x,I,q(\tau),\delta(\tau)}$, we   obtain
\be
\label{etaac}
|\Phi_{x,I}(t,\bar z,y)- \cN_{x,I,\tau}(t,\bar z;y)|\le \tau,\quad x\in \R^m,\,\bar z\in L_\infty(I;\R^m), \, t\in I .
\ee

Regarding the complexity of $\cN_{x,I,\tau}$ we see from \eref{G1} 
that $\fN(G_1)= m+d_y+q +qd_y m^2$, $\fN(G_2)=1$ because of bilinearity.   Thus,
$\fN(G_2\circ G_1)\eqsim m^2qd_y$ so that by \eref{locappr}  and \eref{qeta}  
\begin{align}
\label{Netacompl}
\#\cN_{x,I,\tau}(t,\cdot,\cdot) &\lsim m^2 d_y q(\tau)\Lambda^m \delta(\tau)^{-m}|\log_2\delta(\tau)| \nonumber\\
&\lsim m^2d_y  \frac{2A|I|}{\tau} \frac{(2|I|\Lambda|\uomega|_1)^m}{\tau^m}\Big|\log_2\frac{2|I||\uomega|_1}{\tau}\Big|, 
\quad t\in I.
\end{align}
%
In view of \eref{tripa},  the earlier role of $\|a\|$   will now be played  
by (see \eref{aL})
$
 \max_{t\in\In}\tripnorm{a}_{N_a,m }= A+ \Lambda |\uomega|_1 =:L,
$
and since by \eref{re} $
|I|  L = \frac 12$, we obtain from \eref{Netacompl}\vspace*{-1.5mm}
\be
\label{Netacompl2}
\#\cN_{x,I,\tau}(t,\cdot,\cdot)\lsim m^2d_y A|I|
   \tau^{-(m+1)}\Big|\log_2\frac{2|I||\uomega|_1}{\tau}\Big|.\vspace*{-1.5mm}
\ee
 Thus, with
the same number   $\mu=\mu(\eta)\ge |\log_2(2\eta)^{-1}|$ from \eref{muk} we get 
  $|z(t;x;y)-\Phi^\mu_{x,I}(t;\bar z_x;y)|\le \eta$. 
Hence, the same $Z_\e$, defined by \eref{Ze},
based on tolerances $\eta(\e)$ from \eref{ej}, provide $\e$-accuracy of time-catenated iterates of $\Phi_{x,I_k}$ where the number $K$ of macro-time-steps
still equals $2\Tin   L$. 
We can therefore choose the same vectors of tolerances $\utau^k(\e)$ from \eref{tauke} (i.e., $\tau(\e)\eqsim \eta(\e)\eqsim \e e^{-K/2}$) as well as
tolerances $\delta(\eta(\e))\eqsim \e e^{-K/2}/(2|I||\uomega|_1)$. 
We then define $\cN_\e$ in complete analogy to $\wt Z_\e$ from \eref{Ze}, with 
$A_{w,I,\tau(\e)}$ replaced by $\cN_{w,I,\tau(\e)}$ to obtain $\|z-\cN_\e\|_{L_\infty(\Hom;\R^m)}\le 2\e$.

Hence, on account of \eref{Netacompl2}, since $|I|K=\Tin$, $K/2=   L\Tin$\vspace*{-1.5mm}
\begin{align*}
\#\cN_\e & \lsim m^2d_y K \mu(\eta(\e))A|I|
   \tau(\e)^{-(m+1)}\Big|\log_2\frac{2|I||\uomega|_1}{\tau(\e)}\Big| 
   \le A \Tin m^2d_y
\tau^{-(m+1)}\Big|\log_2\frac{2|I||\uomega|_1}{\tau(\e)}\Big|   \\
&\lsim   A \Tin m^2d_y\Big(\frac{e^{L\Tin}}{\e}\Big)^{m+1}
\Big|\log_2\frac{2|I|e^{L\Tin}|\uomega|_1}{ \e}\Big| 
\lsim  A \Tin m^2d_y\Big(\frac{e^{L\Tin}}{\e}\Big)^{m+1}
\Big|\log_2\frac{ e^{L\Tin} }{ \e}\Big|\\
&=: \phi\Big(\frac{ e^{  L\Tin}}{\e}\Big),\vspace*{-1.5mm}
\end{align*}
where we have used that, by \eref{re}, $1\lsim |I|\Lambda|\uomega|_1\le 1/2$.
This confirms \eref{Nappre}.

We apply Remark \ref{rem:order} and Lemma \ref{lem:algrowth} to conclude that
(surpressing a logarithmic dependence on $A\Tin m^2$)\vspace*{-1.5mm}
\be
\label{pargrowth}
 \wt \gamma(r):= (A\Tin m^2)^{-\frac{1}{m+1}}\Big(\frac{r}{d_y}\Big)^{\frac{1}{m+1}}\Big|\log_2\frac{r}{d_y}
\Big|^{-\frac{2}{m+1}},  \vspace*{-1.5mm}
\ee
satisfies $\wt\gamma  (\phi(s))\eqsim s$. Hence,
there exixts a network $\cN_N$ with $\#\cN_N\le N$ such that
$$
\|z-\cN_N\|_{L_\infty(\Hom;\R^m)} \lsim e^{  L\Tin} \wt\gamma(N)^{-1},\quad N\in \N,
$$
where the constant depends only on $m$. This proves \eref{CNNerr}.

It remains to estimate  
  $\tripnorm{\cN_\e}_{N_\e,m}$ where $N_\e := \#\cN_\e$.
\begin{lemma}
\label{lem:Liptxy}
For $\cN_{\cdot,I,q,\delta}(t, \cdot;\cdot)$ from \eref{Nx} 
the
following statements hold: Let\vspace*{-1.5mm}
\be
\label{Ldelta}
L_\delta :=  \big(A+(\bar\Lambda +\delta)|\uomega|_1\big),\quad \bar \Lambda:= c_3(1+A^\circ)\Lambda,\vspace*{-1.5mm}
\ee  
with $c_3$ from \eref{Lipbound} (see also \eref{aspecial}).
Then,
 for $ (t,x,y),\,(t',x',y')\in \Hom(I)$,\vspace*{-1.5mm}
\be
\label{Lipzy0}
|\cN^\ell_{x,I,q,\delta}(t,\bar z;y)-\cN^\ell_{x,I,q,\delta}(t,\bar z';y')|\le \frac{(  L_\delta|I|)^\ell}{\ell !}\big\{
|y-y'| ,\|\bar z-\bar z'\|_{L_\infty(I;\R^m)}\big\}.\vspace*{-1.5mm}
\ee
Similarly, when $\bar z = \bar z_x $, $\bar z'=\bar z_{x'}$, one has for
all $x, x'\in \R^m$, $t,t'\in I$,  \vspace*{-1.5mm}
\begin{align}
\label{Lipxy0}
|\cN^\ell_{x,I,q,\delta}(t,\bar z_x ;y)-\cN^\ell_{x',I,q,\delta}(t',\bar z_{x'};y')|&\le (A+|\uomega|_1\delta)|t-t'|\nonumber\\
&\qquad + \max\big\{|y-y'|     , |x-x'| \big\} e^{L_\delta|I|} .\vspace*{-1.5mm}
\end{align}
\end{lemma}
\noindent
{\bf Proof:}
Recall from    \eref{Nx} 
that $\cN_{x,I,q,\delta}(t,\bar z;y):= x+  
\sum_{i=1}^{q}\rho_i(t)\sum_{j=1}^{d_y}y_j\omega_j\cN_{j,i,\delta}(\bar z(\xi_i))$,
where $\cN_{j,i,\delta}$ are Lipschitz stable DNNs  approximating
 $a^\circ_{j,J_i}( \cdot)$ (Proposition \ref{prop:Lip}). Thus, by \eref{affinebound}
 \vspace*{-1.5mm} 
\begin{align}
\label{barLambda}
\sum_{j=1}^{d_y}\omega_j|\cN_{j,i,\delta}(\bar z(\xi_i))|  &\le \sum_{j=1}^{d_y}
 |a_{j,J_i}( \bar z(\xi_i))|  +\omega_j|a^\circ_{j,J_i}( \bar z(\xi_i))-\cN_{j,i,\delta}(\bar z(\xi_i))|  
 \le A +  |\uomega|\delta,\vspace*{-1.5mm}
\end{align}  
while, by Proposition \ref{prop:Lip}, \eref{Lipbound},\vspace*{-1.8mm}
\be
\label{boundLam}
  |\cN_{j,i,\delta}(\bar z)- \cN_{j,i,\delta}(\bar z')| \le  \bar \Lambda \|\bar z-\bar z'\|_{L_\infty(I;\R^m)},\quad \bar z,\bar z'\in L_\infty(I;\R^m),\vspace*{-1.5mm}
\ee
 uniformly in $i=1,\ldots,d_y,\,\, j=1,\ldots,d_y$, with $\bar\Lambda$ from \eref{Ldelta}where.
 Then\vspace*{-1.8mm}
\begin{align}
\label{Lip2}
&\Big|\cN_{x,I,q,\delta}(t,\bar z;y)- \cN_{x,I,q,\delta}(t,\bar z';y')\Big| 
\le 
  \sum_{i=1}^{q}\rho_{i}(t)
\sum_{j=1}^{d_y}\omega_j\Big|y_j\cN_{j,i,\delta}(\bar z(\xi_i))- y_j'
\cN_{j,i,\delta}(\bar z'(\xi_i))\Big| \nonumber\\
&  \le  \sum_{i=1}^{q}\rho_{i}(t)\Big\{  \sum_{j=1}^{d_y}\omega_j|y_j-y_j'|  
|\cN_{j,i,\delta}(\bar z(\xi_i))|  +  
\sum_{j=1}^{d_y}\omega_j|y_j'|
|\cN_{j,i,\delta}(\bar z(\xi_i))- \cN_{j,i,\delta}(\bar z'(\xi_i))| \Big\}\nonumber\\
& \le \sum_{i=1}^{q}\rho_{i}(t)\Big\{(A +|\uomega|_1\delta) |y-y'|   + \bar \Lambda
|\uomega|_1 \|\bar z-\bar z'\|_{L_\infty(I;\R^m)}\Big\}\nonumber\\
&\le \sum_{i=1}^{q}\rho_{i}(t)
L_\delta \max\{
|y-y'|  ,\|\bar z-\bar z'\|_{L_\infty(I;\R^m)}\} . 
\end{align}

Similarly, for $t,t'\in I$\vspace*{-1.7mm}
\begin{align*}
|\cN_{x,I,q,\delta}(t,\bar z;y)- \cN_{x,I,q,\delta}(t',\bar z;y)| 
&\le \sum_{i=1}^q |\rho_i(t)-\rho_i(t')| \sum_{j=1}^{d_y}|y_j|\omega_j|\cN_{j,i,\delta}(\bar z(\xi_i))| \\
&\le (A+|\uomega|_1\delta)|t-t'|,\vspace*{-1.5mm}
\end{align*}
where we have used \eref{rhodiff} and \eref{barLambda}. Applying this to $\bar z = \cN^{\ell-1}_{\cdot,I,q,\delta}$ extends this to iterates of $\cN_{\cdot,I,q,\delta}$.

Hence, we are in the same situation as
in \eref{thesame}. Therefore, \eref{Lipzy0} and \eref{Lipxy0} follow by the same arguments as used
in the proof of Lemma \ref{lem:discrLip}.\hfill $\Box$\\

Now recall that $\delta(\e)\eqsim \frac{\e e^{-K/2 }}{2|I||\uomega|_1}$.   Hence,
$
 \delta(\e)|\uomega|_1 \lsim   
\frac{\e     e^{-K/2}}{2|I|}\lsim  \Tin^{-1}\e Ke^{-K/2}\le \e/\Tin.
$ 
Then  \vspace*{-1.5mm}
$$
L_{\delta(\e)}\le A  +\bar\Lambda|\uomega|_1 + \Tin^{-1}{\e}\le \wh L
:= \wh A + \bar\Lambda|\uomega|_1  , \quad \wh A:= \max_{\e\le 1}
A+ \Tin^{-1}\e   ,\vspace*{-1.5mm}
$$
and  the same arguments as used earlier provide \vspace*{-1.5mm}
\be 
\label{PsiLip}
|\wh\Psi_{\utau^k(\e)}(t,x,y)-\wh\Psi_{\utau^k(\e)}(t',x',y')|_\infty
\le \wh A |t-t'|+ \wh L \max\{|x-x'| ,|y-y'| \}.\vspace*{-1.5mm}
\ee
 From these observations it follows that   $\tripnorm{\cN_\e}_{N_\e,m }\le e^{\wh L\Tin}$ and for $  \gamma(r):= \Big(\frac{r}{d_y}\Big)^{\frac{1}{m+1}}\Big|\log_2\frac{r}{d_y}
\Big|^{-\frac{2}{m+1}}$ \vspace*{-1.5mm}
$$
\gamma(N_\e)K_{m }(z,N,\gamma(N)^{-1})\le \gamma(N_\e)\|z-\cN_\e\|_{L_\infty(\Hom)}
+ \tripnorm{\cN_\e}_{N_\e,m } \lsim e^{LT}+ e^{\wh L\Tin},\vspace*{-0.5mm}
$$
which confirms the remainder of the assertion.\hfill $\Box$
\vspace*{-1.5mm} 
\subsection{Proof of Theorem \ref{thm:still}}\label{ssec:still}
\newcommand{\te}{{\t\e}}
\vspace*{-1.5mm}
By assumption, given $\te>0$, there exists an $f_\te\in L_\infty(\In;\Co_{N_\te(f),m})
\cap \Lip_1(\In;C(\R^m\times \cY))$, piecewise affine in time, and a composition $u_{0,\te}\in \Co_{N_\te(u_0),m} $, so that   (identifying for
notational convenience in what follows mappings and representations)  
\be
\label{suchthat0}
\|u_0- u_{0,\te}\|_{L_\infty(\R^m\times \cY)}\le \te,\quad \fN(u_{0,\te})\le \gamma^{-1}\big(\|u_0\|_{\cA^{\gamma,m} }/\te\big)\eqsim \|u_0\|_{\cA^{\gamma,m} }^{1/\alpha}\te^{-1/\alpha},  
\ee
and likewise, since $\fN(f_\te(t,\cdot))\le \gamma^{-1}\big(\|f(t,\cdot)\|_{\cA^{\gamma,m} }/\te\big)$ for $t\in \In$
\be
\label{suchthatf}
\|f- f_{\te}\|_{L_\infty(\In;C(\R^m\times \cY)))}\le \te,\quad \fN(f_\te(t,\cdot))
\lsim 
\|f(t,\cdot)\|_{\cA^{\gamma,m} }^{1/\alpha}\te^{-\frac{1}{\alpha}}. 
\ee
Next, we use that the compositional factors in $u_{0,\te}, f_\te$ are Lipschitz
continuous with constants controlled by 
$
\|u_0\|:=\|u_0\|_{\cA^{\gamma,m} }\,, 
\|f\|:=\|f\|_{L_\infty(\In;\cA^{\gamma,m} )},
$
respectively. We employ Lemma \ref{lem:import} to implant $\eta$-accurate Lipschitz controlled DNNs into $u_{0,\te}, f_{\te}$, respectively. Invoking Remark
\ref{rem:error} we obtain\vspace*{-1.5mm} 
\be
\label{uf0}
\|u_{0,\te}- \cN_{u_{0,\te},\eta}\|_{L_\infty(D)}\le \eta \Big\{1+ \sum_{j=1}^{n(\Di(u_{0,\te}))-1}  \|u_0\|\Big\}\le \eta \fN(u_{0,\te})\|u_0\|,\vspace*{-1.5mm}
\ee 
as well as\vspace*{-1.5mm}
\be
\label{uf1}
\|f_\te(t,\cdot) - \cN_{f_\te,\eta}(t,\cdot)\|_{L_\infty(\R^m\times \cY)} 
\le \eta \Big\{1+ \sum_{j=1}^{n(\Di(f_{\te}))-1}  \|f_0\|\Big\}\le \eta \fN(f_\te(t))\|f_0\|.
\quad t\in I,\vspace*{-1.5mm}
\ee 
By Proposition \ref{prop:Lip}, \eref{suchthat0},  and \eref{suchthatf}, it follows that\vspace*{-1.5mm}
\be
\label{ufcompl0}
\# \cN_{u_{0,\te},\eta}\le \tripnorm{u_{0,\te}}^m_{\fN (u_{0,\te}),m }\fN (u_{0,\te})\eta^{-m}|
\log_2\eta| \le \|u_0\|^{m+\frac{1}{\alpha }}  \te^{-\frac{1}{\alpha }}\eta^{-m}|\log_2\eta| ,\vspace*{-1.5mm}
\ee
and, uniformly in $t\in \In$,
\be
\label{ufcompl1}
\# \cN_{f_\te, \eta}(t,\cdot)\le \tripnorm{f_{\te}(t,\cdot)}^m_{\fN(f_\te(t)),m,\Lip}\fN (f_\te(t))\eta^{-m}|
\log_2\eta| \le \|f\|^{m+\frac{1}{\alpha }}  \te^{-\frac{1}{\alpha }}\eta^{-m}|
\log_2\eta| .
\ee
Employing again a time-discretization of size $q$, \eref{3.6a} suggests
the following DNN approximation to $u$ which yields, on account of Lemma
\ref{lem:quadrature},\vspace*{-1.9mm}
\be
\label{fint}
\Big|\int_0^t f_\te(s,w 
,y))ds - \sum_{i=1}^q \rho_i( t) f_\te(\xi_i,w,y)
\Big| \le \frac{\|f\|\Tin^2}{2q}.\vspace*{-1.5mm}
\ee 
Finally,
we know from Theorem \ref{thm:parconv} that there exists a DNN $\cN_{z,\te_z}$
that approximates the characteristic field $z$ within accuracy $\te$, i.e.,
 in view of \eref{Nappre},\vspace*{-1.5mm}
\be
\label{weknow}
\|z- \cN_{z,\te}\|_{L_\infty(\Hom)}\le \te ,\quad
\#\cN_{z,\te }\lsim  d_y\Tin  \Big(\frac{e^{L\Tin }}{\te }\Big)^{m+1}\Big|\log_2 \Big(\frac{e^{L\Tin }}{\te }\Big)\Big|^2 ,\vspace*{-1.5mm}
\ee
where we suppress in what follows the dependence on $A,m$ and where $L$ is given by \eref{Lspecial}. 
 
In summary, the network $\cN_{u,\te,\eta,q}$ formed by composing
 the DNNs 
$\cN_{u_0,\te,\eta}, \cN_{f_\te ,\eta}(\xi_i,\cdot)$, $i=1,\ldots,q$, with
the approximate characteristics $\cN_{z,\te_z}$
satisfies, on account of \eref{Gsum},  \eref{Nplus}, \eref{ufcompl0},  and
\eref{ufcompl1},
\begin{align}
\label{totalcompl}
\# \cN_{u,\te,\eta,q}&\le  \#\big(\cN_{u_0,\te,\eta}\circ \cN_{z,\te}) + 
q \max_{i=1,\ldots,q} \#\big(\cN_{f_\te,\eta}(\xi_i,\cdot)\circ \cN_{z,\te })
\nonumber\\
&= \# \cN_{u_0,\te,\eta} + \# \cN_{z,\te} + q\big(\max_{i=1,\ldots,q} \# \cN_{
f_\te\,\eta}( \xi_i,\cdot) + \# \cN_{z,\te}\big)\nonumber\\
&\lsim \Big(\|u_0\|^{m+\frac{1}{\alpha}}\te^{-\frac{1}{\alpha}} + q\|f\|^{m
+\frac{1}{\alpha}} 
\te^{-\frac{1}{\alpha}}\Big)\eta^{-m}|
\log_2\eta| \nonumber\\
&\qquad\quad + (1+q)d_y\Tin  \Big(\frac{e^{L\Tin }}{\te}\Big)^{m+1}\Big|\log_2 \Big(\frac{e^{L\Tin }}{\te}\Big)\Big|^2\nonumber\\
&\lsim q \Big\{ M^{m+\frac{1}{\alpha}}\te^{-\frac{1}{\alpha}} \eta^{-m}|\log_2 \eta|
+ d_y\Tin \Big(\frac{e^{L\Tin }}{\te}\Big)^{m+1}\Big|\log_2 \Big(\frac{e^{L\Tin }}{\te}\Big)\Big|^2\Big\},
\end{align}
where we we set
$
M:=\max\{1,\|u_0\|,\|f\|\}.
$
To determine $\eta$, we have by \eref{3.6a},
\begin{align*}
&\|u(t, \cdot) - \cN_{u,\te,\eta,q}(t,\cdot)\|_{L_\infty(\R^m\times \cY)}
\le \|u_0(z(-t,t,\cdot;\cdot))- \cN_{u_0,\te,\eta}\circ \cN_{z,\te}(-t,t,\cdot;\cdot)\|_{L_\infty(\R^m\times \cY)} \\
&\qquad\qquad + \sup_{x,y}\left|\int_0^t f(s,z(t-s,t,x,y) ds -\sum_{i=1}^{q}
\rho_i(t) \cN_{\xi_i,f_\te,\eta}\circ \cN_{z,\te}(t -\xi_i,x,y)\right|\\
&\qquad\qquad =: Q_1 + \sup_{x,y} Q_2(x,y).
\end{align*}
Regarding $Q_1$, let $L_0:= |u_0|_{\Lip_1(\R^m\times \cY)}\le \|u_0\|$.
Because of \eref{uf0} and \eref{suchthat0},
\begin{align}
\label{Q1}
Q_1 &\le \|u_0(z(-t,t,\cdot;\cdot))- u_0 (\cN_{z,\te}(-t,t,\cdot;\cdot))\|_{L_\infty(\R^m\times \cY)}\nonumber\\
&\qquad  + \| u_0 (\cN_{z,\te}(-t,t,\cdot;\cdot))- u_{\te,0}(\cN_{z,\te}(-t,t,\cdot;\cdot))\|_{L_\infty(\R^m\times \cY)} \nonumber\\
&\qquad + \| u_{\te,0}(\cN_{z,\te}(-t,t,\cdot;\cdot))
- \cN_{u_0,\te,\eta}(\cN_{z,\te}(-t,t,\cdot;\cdot))\|_{L_\infty(\R^m\times \cY)}\nonumber\\
&\qquad
\le (1+L_0)\te + \|u_0\|\eta \fN(u_{0,\te})
{\lsim}(1+L_0)\te + 
\|u_0\|^{1+\frac{1}{\alpha}}\eta \te^{-\frac{1}{\alpha}}\nonumber\\
&\le \|u_0\|\Big\{2\te + \eta\te^{-\frac{1}{\alpha}}\|u_0\|^{\frac{1}{\alpha}}\Big\},
 \end{align}
where we have used $\|u_0\|\ge 1$ and \eref{dataappr} in the last step. Similarly,
by \eref{suchthatf}
\eref{uf1}, and \eref{fint}, abbreviating $L_f:= |f|_{\Lip_1}(\In;C(\R^m\times \cY))$,\vspace*{-2mm}
\begin{align}
\label{Q20}
Q_2(x,y) &\le \Tin L_f\te  +  
   \int_0^t |f(s,\cN_{z,\te}(t-s,t,x,y)-
f_\te(s,\cN_{z,\te}(t-s,t,x,y))| ds\nonumber\\
&  + \int_0^t \Big|f_\te(s,\cN_{z,\te}(t-s,t,x,y)) - 
\sum_{i=1}^{q}
\rho_i(t) f_\te(\xi_i,\cN_{z,\te}(t -\xi_i,x,y))\Big|\nonumber\\
&+  \sum_{i=1}^{q}\rho_i(t)\big|f_\te(\xi_i,\cN_{z,\te}(t -\xi_i,x,y))-
 \cN_{f_\te,\eta}(\xi_i,\cN_{z,\te}(t -\xi_i,x,y))\big|. 
\end{align}
On account of \eref{om0} and the assumption $L_f\le \|f\|$,    this gives
\begin{align}
\label{Q2}
Q_2(x,y)& \le (1+\|f\| )\Tin \te +\frac{\|f\| \Tin^2}{2q}+ \Tin \eta \fN(f_\te)\|f\| 
\le  2\|f\| \Tin \te +\frac{\|f\|\Tin^2}{2q}+ \Tin \frac{\eta}{\te^{1/\alpha}}\|f\|^{1+\frac{1}{\alpha }}
\nonumber\\
&
\le \Tin \|f\|\Big\{2\te +\frac{\Tin}{2q}+\|f\|^{\frac{1}{\alpha}}\eta \te^{-\frac{1}{\alpha}}\Big\}.
\end{align}
Now recall  that $M= \max\{1, \|f\|,\|u_0\|\} $ and let
\be
\label{parchoices}
 q(\te)=\frac{\Tin  }{2\te  },\quad \eta(\te)=M^{-\frac{1}{\alpha}} \te^{1+
 \frac{1}{\alpha  }},
\ee
to conclude that $\max_{x,y}Q_2(x,y)\le 4\Tin\|f\|\te$.
Hence, we derive  from \eref{Q1} and \eref{Q2} that the network
$\cN_{u,\te}:= \cN_{u,\te,\eta(\te),q(\te)}$ satisfies (recall that by assumptions   $\|u_0\|,\|f\|\ge 1$)
\begin{align}
\label{finestu}
\|u- \cN_{u,\te}\|_{\infty}&\le \{4\Tin\|f\|+3\|u_0\|\}\te \le 7\Tin M\te.
\end{align}
This confirms the first part of \eref{uappr} with $\e:= 7\Tin M\te$.

Now
 we infer from \eref{totalcompl}, \eref{suchthat0}, \eref{suchthatf} that
\begin{align*}
\#\cN_{u,\te}&\lsim \Tin \te^{-1}\Big\{ M^{\frac{\alpha m+ 1}{\alpha}} 
\te^{-\frac{1}{\alpha}} M^{\frac{m}{\alpha}}\te^{-\frac{m(1+\alpha)}{\alpha}}\Big|\log_2 \frac{M}{\te^{ \alpha +1} }\Big|
+ d_y\Tin \Big(\frac{e^{L\Tin }}{\te}\Big)^{m+1}\Big|\log_2 \Big(
\frac{e^{L\Tin }}{\te}\Big)\Big|^2\Big\}\\
&\lsim \Tin \te^{-1}\Big\{ M^{\frac{\alpha m+ 1}{\alpha}} 
\te^{-\frac{1}{\alpha}} M^{\frac{m}{\alpha}}\te^{-\frac{m(1+\alpha)}{\alpha}}
+ d_y\Tin e^{\Tin L(m+1)} \te^{-(m+1)}\Big\}\Big|\log_2 \Big(
\frac{e^{L\Tin }}{\te}\Big)\Big|^2\\
&= \Tin  \Big\{M^{\frac{(\alpha+1) m+ 1}{\alpha}}\te^{-\frac{(1+\alpha)(m+1)}{\alpha}}
+ d_y\Tin e^{\Tin L(m+1)}\te^{-(m+2)}\Big\}\Big|\log_2 \Big(
\frac{e^{L\Tin }}{\te}\Big)\Big|^2. 
\end{align*}
Introducing
 $
\beta:= \max\{1, (m+1)/\alpha\},
$ (see \eref{where})
and substituting $\te = \e/(7\Tin M)$, yields upon elementary calculations
\begin{align*}
\#\cN_{u,\te} &\lsim   \Big\{M^{\frac{(\alpha+1)m+1}{\alpha}}\Tin^{m+2+\beta}
e^{-L\Tin(m+1+\beta)} + d_y\Tin^{m+4} e^{-L\Tin\beta}\Big\} \Big(\frac{Me^{L\Tin}}{\te}\Big)^{m+1+\beta}\Big|\log_2 \Big(
\frac{e^{L\Tin }}{\te}\Big)\Big|^2.
 \end{align*}
The terms $\Tin^{m+2+\beta}
e^{-L\Tin(m+1+\beta)}$, $\Tin^{m+4} e^{-L\Tin\beta}$ remain uniformly bounded
for all $\Tin >0$ with a constant that actually decreases when $L$ gets large.
Thus, fixing $M$, a large parametric dimension in the second summand dominates, giving
\be
\label{phi} 
\#\cN_{u,\te}\lsim  \max\{M^{\frac{(\alpha+1)m+1}{\alpha}},d_y\}
\Big(\frac{Me^{L\Tin}}{\te}\Big)^{m+1+\beta}\Big|\log_2 \Big(
\frac{e^{L\Tin }}{\te}\Big)\Big|^2=: \phi(M e^{L\Tin }/\e),
\ee
which proves \eref{uappr}. 
 
Regarding the remainder of the claim, recall from Theorem \ref{thm:parconv} that the approximations $\cN_{z,\e}$ have uniformly bounded composition norms $\tripnorm{\cN_{z,\e}}_{\#\cN_{z,\e},m }\lsim e^{L\Tin}$, see also Lemma \ref{lem:Liptxy}. To bound the composition norms of $\cN_{u,\te}$, we recall
that the composition norms of the network approximations to $u_0$ and $f$ 
are bounded by $M= \max\{1,\|u_0\|,\|f\|\}$. We then infer from Remark \ref{rem:closedness} (see also 
\eref{G1G2Lip} and \eref{maxplus}), applied to the first line of \eref{totalcompl},
that 
$$
\tripnorm{\cN_{u,\te}}_{\#\cN_{u,\te},m }\lsim M e^{L\Tin},
$$
which is the asserted stability estimate.
The coonvergence rate \eref{urate} follows from   Remark \ref{rem:order} and Lemma \ref{lem:algrowth} applied to to the growth function $\phi$ in \eref{phi}. \hfill $\Box$\\

Regarding Remark \ref{rem:exp},
The same reasoning applies (with slightly simpler
technicalities), replacing  \eref{parchoices} by $\eta(\te):= \te/(\ln((\|u_0\|+\|f\|)/\te))$ while keeping $q(\te)$ the same,  \hfill $\Box$


\vspace*{-2mm}
\section*{Appendix A}\label{sssec:propLip}
\paragraph{Proof of Proposition \ref{prop:Lip}}\label{sssec:propLip}
 %
 In this section we build mainly on findings from \cite{GKP,Yar}.
Consider    
  the ``hat-function'' $\phi(x):= (1-|x|)_+=\max\{0,1-|x|\}$, $x\in \R$,
  as well as the scaled and shifted versions $\phi_{i,h}(x):=  \phi(h^{-1}x-i)$, $i\in \Z$,
 with support $S_i= [(i-1)h,(i+1)h]$. 
 We let $h=1/q$ for some integer $q\in\N$, so that the restrictions of the $\phi_{i,h}$ to $[0,1]$ form a stable basis for all piecewise linears on $(0,1)$
 subordinate to the partition induced by the nodes $\{ih=i/q:i=0,\ldots, q\}$.
 Since each $\phi_{i,h}$ is a second order divided difference of the ReLU rectifier
 $\sigma(x):= x_+$ with respect to the nodes $ih$ it has an exact representation as a univariate neural network of fixed finite depth and a fixed finite number of weights.
 The Lipshitz constant of   $\phi_{i,h}$ and hence of this network is 
 clearly $h^{-1}$. Abbreviating $\bi:=(i_1,\ldots,i_\ss)\in \{0,\ldots,q\}^\ss$
 we consider next for $x=(x_1,\ldots,x_\ss)\in \R^\ss$ the tensor products 
 $\,
 \phi_{\bi,h}:= \phi_{i_1,h}(x_1)\cdots \phi_{i_\ss,h}(x_\ss),\,
 $
 which obviously satisfy $\phi_{\bi,h}(\bi')=\delta_{\bi,\bi'}$ for any
 $\bi,\bi'\in \cI_h:= \{0,h,\ldots,hq\}^\ss$ while we still have 
 $\|\partial_j\phi_{\bi,h}\|_\infty\le h^{-1}$. 
 The next step consists in approximating each $\phi_{\bi,h}$, viz. a product
 of univariate ReLU networks of fixed depth and number of weights by
 a ReLU network of input dimension $\ss$. This is where one uses that the function
 $M:\nu= (\nu_1,\ldots,\nu_\ss)\mapsto \prod_{j=1}^\ss \nu_s$ can be approximated 
 by a ReLU network $\cN_{M,\delta}$ according to \vspace*{-1.5mm}
 \be
 \label{prod}
 \|M -\cN_{M,\delta}\|_{W^k(L_\infty((0,1)^\ss)} \le \delta, \quad k\in \{0,1\},
 \vspace*{-1.5mm}
 \ee
 where  the depth of $\cN_{M,\delta}$ as well as   $\#\cN_{M,\delta}$
 is bounded by a constant multiple of $\log_2\delta^{-1}$, with constants
 depending only on $\ss$. Moreover, $\cN_{M,\delta}(0)=0$.
 The case $k=0$ in \eref{prod} appears already in \cite{Yar}. 
 A key observation in \cite[\S~C]{GKP} is that $k=1$ still holds under the same
 complexity bounds. This is then used to show that for each $\bi\in \cI_h$
 there exists a ReLU network $\cN_{\bi,\delta}$ such that\vspace*{-2.8mm}
 \be
 \label{Nbi}
 \begin{array}{c}
 \|\phi_{\bi,h}- \cN_{\bi,\delta}\|_{W^k(L_\infty((0,1)^\ss))}\infty \le c^k\delta h^{-1},\quad k\in\{0,1\},\\
   \|\cN_{\bi,\delta}\|_{\Lip_1}\le  c\delta^{-1},\quad {\rm supp}\,\cN_{\bi,\delta}\subseteq{\rm supp}\,\phi_{\bi,h},\vspace*{-1.5mm}
   \end{array}
 \ee
 (with a constant $c$, depending on $\ss$) and
 \be
 \label{complss}
 \#\cN_{\bi,\delta}, \,\mbox{depth of } \cN_{\bi,\delta}\lsim \log_2\frac{1}{\delta},
 \ee
 with constants depending only on $\ss$.
  Now, given $g$, consider the interpolant\vspace*{-1.5mm}
\be
\label{gh}
g_h := \sum_{{\bf i}\in \cI_h} g({\bf i}h)\phi_{{\bf i},h}.\vspace*{-1.5mm}
\ee
Obviously $\|g_h\|_\infty\le \|g\|_\infty$. We claim that $g_h$ is also Lipschitz
continuous. To see this, let $\cI_h(x):= \{\bi\in \cI_h: \phi_{\bi,h}(x)\neq 0\}$
denote the collection of those nodes whose basis functions contain $x\in (0,1)^\ss$
in the interior of their support. Then for $x\in (0,1)^\ss$ let $\bi(x)\in \cI_h(x)$ denote the node closest to the Chebyshev center of the convex hull $\big[\cI_h(x)\big]$ of $\cI_h(x)$. Then, since $g(\bi(x)h)\sum_{\bi\in\cI_c(x)}\phi_{\bi,h}(x)$  is constant in a neighborhood of $x$, one has\vspace*{-1.5mm}
$$
|\partial_j g_h(x)|= \Big|\sum_{\bi\in \cI_h(x)}(g(\bi h) -g(\bi(x)h))\partial_j
\phi_{\bi,h}(x)\Big|\le \max_{\bi\in\cI_h(x)}|g(\bi h) -g(\bi(x)h)|\frac 1h
\#(\cI_h(x)).\vspace*{-1.5mm}
$$
Since ${\rm diam}\,\cI_h(x)\le c h$ for a constant depending on $\ss$, this yields\vspace*{-1.5mm}
\be
\label{pgh}
|\partial_j g_h(x)|\le c^{-1}\|g\|_{\Lip_1},\quad x\in (0,1)^\ss,\vspace*{-1.5mm}
\ee
from which it follows that (weakly)\vspace*{-1.9mm}
\be
\label{ghLip}
|g_h(x)-g_h(x')|= \Big|\int_0^1 \nabla g_h(x+t(x'-x))\cdot(x'-x)dt\Big|
\le |x-x'|\sqrt{\ss}c^{-1}\|g\|_{\Lip_1}.\vspace*{-1.5mm}
\ee
Moreover, 
\be
\label{gdiff}
|g(x)-g_h(x)|=\Big|\sum_{i\in \cI_h(x)}(g(x)-g(\bi h))\phi_{\bi,h}\Big|
\le C h\|g\|_{\Lip_1},
\ee
since $\max\{|x-\bi h|: \bi\in\cI_h(x)\}\le Ch$ whith $C$ depending only on $\ss$.
Given $\delta>0$, the choice $h=h(\delta)\le \frac{\delta}{2C\|g\|_{\Lip_1}}$ ensures
$|g(x)-g_h(x)|\le \frac{\delta}2$. Now approximate each $\phi_{\bi,h(\delta)}$
by a ReLU network $\cN_{\bi,\delta}$ with accuracy $\|\phi_{\bi,h(\delta)}-
\cN_{\bi,\delta}\|_\infty \le c^*\delta$, with $c^*$ to be determined in a moment. We obtain, by \eref{Nbi}\vspace*{-1.5mm}
\begin{align*}
\Big\|g- \sum_{\bi\in \cI_{h(\delta)}} g(\bi h)\cN_{\bi,\delta}\|_\infty
&\le \|g- g_h\|_\infty + \sup_{x\in (0,1)^\ss}\sum_{\bi\in \cI_{h(\delta)(x)}}
|g(\bi h)|\big|\phi_{\bi,h(\delta)}(x)- \cN_{\bi,\delta}(x)\big|\\
&\le \frac{\delta}2 + \sup_{x\in (0,1)^\ss}\#\cI_{h(\delta)}(x) c^*\delta.\vspace*{-1.5mm}
\end{align*}
Thus, choosing $c^*= (2\sup_{x\in (0,1)^\ss}\#\cI_{h(\delta)}(x))^{-1}$,
we have confirmed \vspace*{-1.5mm}
\be
\label{confirmed}
\Big\|g- \sum_{\bi\in \cI_{h(\delta)}} g(\bi h)\cN_{\bi,\delta}\|_\infty\le \delta.
\vspace*{-1.5mm}
\ee
Moreover, defining
$
\cN_{\delta}:= \sum_{\bi\in \cI_h}g(\bi h)\cN_{\bi,\delta},
$
we obtain, again by \eref{Nbi} and \eref{pgh},
\begin{align*}
\big|\partial_j \cN_\delta (x)\big| &\le \big|\partial_j (\cN_\delta (x)-g_h(x))\big| + c^{-1}\|g\|_{\Lip_1}=
\Big|\sum_{\bi\in \cI_{h(\delta)}}  g(\bi h) \big(\partial_j\cN_{\bi,\delta}(x)-
\partial_j \phi_{\bi,h}\big)\Big| + c^{-1}\|g\|_{\Lip_1}\\
&\le c'\|g\|_\infty \#(\cI_{h(\delta)})h^{-1}\delta + c^{-1}\|g\|_{\Lip_1} 
\le c''(\|g\|_\infty +1)  \|g\|_{\Lip_1},
\end{align*}
where $c''$ depends only on $\ss$ and where we have used that $\frac{\delta}h
\le 4C\|g\|_{\Lip_1}$. This confirms \eref{Lipbound}.
 Regarding the complexity \eref{Ndcompl} of $\cN_\delta$, we have 
 $
 \#(\cI_h(\delta))=h(\delta)^{-\ss}\le c \delta^{-\ss}\|g\|_{\Lip_1}^\ss,
 $
 which completes the proof because $\#\cN_\delta \lsim \#(\cI_h(\delta))
 \log_2\frac{1}{\delta}$.\hfill $\Box$
 \vspace*{-3mm}
\paragraph{Proof of Proposition \ref{prop:comp}}\label{sssec:propcomp}
\renewcommand{\XX}{\mathbb{X}}
For each finite $N$ there is only a finite number of feasible
dimensionality vectors $\Di$ with $\fN(\Di)\le N$ which representations of $G\in \Co_N$  may have.
For each such $\Di$
 consider first
the following auxiliary classes.
Let $\F_\ell\subset C(\R^{d_{\ell-1}};\R^{d_\ell})$ be  compact  
and let 
$
\Co(\Di,\F) := \{G\in \mathbb{X}_0, G=G_\bg, \, \Di(\bg)=\Di:   \mbox{  with }\,g^j\in \F_j,\,j=1,\ldots, n(\Di)\}.
$ 
\begin{lemma}
\label{lem:F}
The collection 
\be
\label{Co}
\Co_N(\F) := \bigcup_{\fN(\Di)\le N} \Co(\Di,\F).
\ee
is compact in $C(D_0;\R^{d_{n(\Di)}})$ (and so are the subsets $\Co(\Di,\F)$).
\end{lemma}
To establish first this Lemma, it is enough to confirm compactness 
of $\Co(\Di,\F)$ for each of the eligible $\Di$. To this end, let $\F,\F'$, 
be dimensionally compatible compact subclasses so that for 
$g\in \F, h\in \F'$ compositions $h\circ g$ are defined. Then $\{h\circ g: h\in \F',\,g\in \F\}$ is compact in $\XX$.
In fact, let $(g_n)_{n\in\N}$, $(h_n)_{n\in\N}$ be uniformly bounded sequences in $\F$, $\F'$, respectively.
By Arzela-Ascoli's Theorem they are equicontinuous and have a convergent subsequence with continuous limits $g, h$, say.
Denote these subsequences again by  $(g_n)_{n\in\N}$, $(h_n)_{n\in\N}$). Then,  
\begin{align*}
\|h_n\circ g_n - h\circ g\|_{\XX(D_g)}&\le \|h_n\circ g_n - h_n\circ g\|_{\XX(D_g)}+ \| h_n\circ g - h\circ g\|_{\XX(D_g)}. 
\end{align*}
By uniform convergence of the $g_n$ and equicontinuity of the $h_n$ the first summand becomes arbitrarily small for $n$ large
enough. By uniform convergence of the $h_n$ the second summand gets small as well. Iterating this argument, shows that
$\Co(\Di,\F)$ is compact. Since $\Co_N(\F)$ is a finite union of such sets the assertion of Lemma \ref{lem:F} follows.

The proof of Proposition \ref{prop:comp} follows now from noticing
that membership to $\Co_{N}( B )$ requires all composition factors to have
uniformly bounded Lipschitz norm and hence belong to compact classes.
\hfill $\Box$
\vspace*{-3mm}
\paragraph{Proof of Remark \ref{rem:Acomp}:} 
Fix $B< \infty$ and let  
$(f_j)_{j\in \N}$ be a sequence in $\cK_{\gamma,\ss,\cR^\circ }(B)$.  Take a sequence $(\e_k)_{k\in\N}$ of numbers decreasing monotonically
to zero. For each $f_j$ let $\bg_{j,k}$ denote a compositional representation 
of some function in   $\Co_{N_{\e_k},\ss}$ such that (see \eref{estA})
$\tripnorm{G_{\bg_{j,k}}}_{N_{\e_k},\ss,\cR^\circ}\le B$ and
$\|f_j - G_{\bg_{j,k}}\|_\XX \le \e_k$. The complexity function $\fN(\bg_{j,k})$ is controlled uniformly in $j$ by $N_{\e_k}$,
defined by \eqref{Ne}. For fixed $k$ the class $\Co_{N_{\e_k},\ss}$ is compact
(Proposition \ref{prop:comp}). Therefore, for fixed $k$,   $(G_{\bg_{j,k}})_{j\in\N}$ contains a subsequence (again denoted 
by $(G_{\bg_{j,k}})_{j\in \N}$), converging uniformly to some $G_k\in \Co_{N_{\e_k},\ss}$.  Now one can take a diagonalization argument, 
letting $k$ tend to infinity, extracting a convergent subsequence from $(f_j)_{j\in\N}$.\hfill $\Box$
\vspace*{-2mm}
\section*{Appendix B}
\vspace*{-2mm}
\paragraph{Proof of Lemma \ref{lem:quadrature}} 
As for  (a), let $i(t):= \argmin_{i=1,\ldots,q}|t-\xi_i|$ and since
\be
\label{rhovalues}
\rho_i(t)= 
\left\{\begin{array}{ll}
0, &t< \tau_{i-1},\\
t-\tau_{i-1},& t\in J_i,\\
|I|/q,& t> \tau_i.
\end{array}\right.
\ee
\eref{sumrho} follows.  Regarding \eref{rhodiff}, without loss of generality 
 assume that $t\le t'$ so that\\[1.5mm]
{\scriptsize 
 \begin{tabular}{c|c|c|c|c|c}
 &  $t, t'\le \xi_{i-1}$ or $ t, t'\ge \xi_{i}$ & $t<\xi_{i-1}, t'\in J_i$ 
 & $t\le \xi_{i-1}, t'> \xi_i$ & $t, t'\in J_i$ & $t\in J_i, t'>\xi_i$\\
 \hline\\[-2.5mm]
$\rho_i(t)-\rho_i(t')$ & $0$ & $\xi_{i-1}-t'$ & $|J_i|= \frac{|I|}q$ & $t-t'$ &
$t-\xi_{i-1} - |I|/q$
\end{tabular}
}

\vspace*{2.5mm}
%
\noindent
Hence $\rho_i(t)\ge \rho_i(t')$,
$i=1,\ldots,q$. Specifically, assume that $t'\in J_\nu$, $t\in J_\ell$.
Then
\begin{align*}
\sum_{i=1}^q|\rho_i(t)-\rho_i(t')| & = \sum_{i=\ell}^\nu \rho_i(t)-\rho_i(t')
 = |J_\ell|-(t-\tau_{\ell-1}) + |J_{\ell+1}| +\cdots + |J_{\nu-1}|+ t'-\tau_{\nu-1}
\\
&=t'-t +\tau_\nu-\tau_\ell -(\tau_{\nu-1}-\tau_{\ell-1})= t'-t,
\end{align*}
confirming  claim (a).

As for (b), we obtain
\begin{align*}
  \int_\uT^{t} g(s)ds - \sum_{i=1}^q \rho_i(t)g_{J_i} &=\sum_{i=1}^q
 \int_{J_i}\Big(\chi_{s\le t}(s)- \frac{\rho_i(t)}{|J_i|}\Big)g(s)ds.
 \end{align*}
Now suppose  that $t\in J_k$. By \eref{rhovalues},  we have $\Big(\chi_{s\le t}(s)- \frac{\rho_i(t)}{|J_i|}\Big)|_{J_i}=0$ for $i\le k-1$ while elementary calculations 
yield\vspace*{-1.7mm}
\begin{align*}
\Big|\int_{\tau_{k-1}}^t\Big(g(s)- \frac{\rho_k(t)}{|J_k|}\Big)ds\Big| &
= \Big|\frac{\tau_k-t}{\tau_k-\tau_{k-1}}\int_{\tau_{k-1}}^t g(s)ds -
\frac{t-\tau_{k-1}}{\tau_k-\tau_{k-1}} \int_t^{\tau_k}g(s)ds\Big|\\
&\le \frac{\tau_k-t}{\tau_k-\tau_{k-1}}\Big\{(t-\tau_{k-1})\|g\|_{L_\infty(J_k)}
+ (t-\tau_{k-1})\|g\|_{L_\infty(J_k)}\Big\}\\
&\le 2\|g\|_{L_\infty(J_k)}\frac{(t-\tau_{k-1})(\tau_k-t))}{\tau_k-\tau_{k-1}}
\le \frac{|J_k| \|g\|_{L_\infty(J_k)}}{2},\vspace*{-1.7mm}
\end{align*}
which is \eref{om0a}.

Concerning (c),
By assumption $\int_{J_i}|g(s)-g(\xi_i)|ds\le |J_i|^2 L' /2$ 
 so that\vspace*{-1.7mm}
\begin{align*}
 \Big|\int_0^{t} g(s)ds - \sum_{i=1}^q \rho_i(t)g(\xi_i)\Big|&\le
   \sum_{i=1}^{i(t)-1} \int_{J_i } |g(s)- g(\xi_i)|ds 
   + \int_{\tau_{i(t)-1}}^t|g(s)- g(\xi_{i(t)})|ds\\
   & \le i(t)\frac{|J_i|^2L'}2 =\frac{i(t)}{q} \frac{|I|^2L'}{2q}.\qquad\Box
 \end{align*}
\paragraph{Proof of Lemma \ref{lem:discrLip}}
We consider first the case of fixed $x$. 
Using \eref{Fqgen} and \eref{tAicontrol},
one finds for $t\in J_\ell$\vspace*{-1.5mm}
\begin{align}
\label{thesame}
|A_{x,I,q,N}(t,\bar z;y)-A_{x,I,q,N}(t,\bar z';y')|
&\le  \sum_{i=1}^\ell \rho_i(t) \|a\|\max\{ |y-y'|  , |\bar z(\xi_i)-\bar z'(\xi_i)|\},\vspace*{-1.5mm}
\end{align}
where $\|a\|$ plays the role of $L$. To see the pattern,\vspace*{-1.5mm}
\begin{align*}
|A^2_{x,I,q,N}(t,\bar z;y)-A^2_{x,I,q,N}(t,\bar z';y')|
&\le  \sum_{i_1=1}^\ell \rho_i(t)\|a\| |A_{x,I,q,N}(\bar z(\xi_{i_1}),\xi_{i_1};y)-
A_{x,I,q,N}(\bar z(\xi_{i_1}),\xi_{i_1};y')|\\
&\le  \sum_{i_1=1}^\ell \rho_{i_1}(t)\|a\|\Big\{ \sum_{i_2=1}^{i_1}\rho_{i_2}(t)\|a\|
\max\{|y-y'|, |\bar z(\xi_{i_2})-\bar z'(\xi_{i_2})|\}\Big\}.\vspace*{-1.5mm}
\end{align*}
Inductively it follows that for $t\in J_\ell$ and 
$k\in \N$\vspace*{-1mm}
\begin{align*}
|A^k_{x,I,q,N}(t,\bar z;y)-A^k_{x,I,q,N}(t,\bar z';y')|&\le  
\|a\|^k \sum_{\ell\ge i_1\ge i_2\ge\cdots\ge i_k\ge 1} 
\rho_{i_1}(t)\cdots\rho_{i_k}(t)\max\{|y-y'|,|\bar z(\xi_{i_k})-\bar z'(\xi_{i_k})\}|.\vspace*{-1.5mm}
 \end{align*}
Invoking \eref{om0}, yields\vspace*{-1.5mm}
$$
\sum_{\ell\ge i_1\ge i_2\ge\cdots\ge i_k\ge 1} 
\rho_{i_1}(t)\cdots\rho_{i_k}(t)\le |I|^k \sum_{\ell\ge i_1\ge i_2\ge\cdots\ge i_k\ge 1} 
 \frac{i_1\cdots i_k}{q^k}\le \frac{|I|^k}{k!},\vspace*{-1.5mm}
$$
providing
\begin{align}
\label{smooth2}
|A^k_{x,I,q,N}(t,\bar z;y)-A^k_{x,I,q,N}(t,\bar z';y')|&\le  
 \frac{(\|a\||I|)^k}{k!}\max\{|y-y'| ,\|\bar z-\bar z'\|_{L_\infty(I;\R^m)}\}.%
\end{align}
Similarly, for $\bar z= \bar z_x, \bar z'= \bar z_{x'}$ (see \eref{initial}), we obtain
\begin{align*}
|A_{x,I,q,N}(t,\bar z;y)-A_{x',I,q,N}(t,\bar z';y')|
&\le |x-x'| + \sum_{i=1}^\ell \rho_i(t)\|a\|\max\{ |y-y'|  , \|\bar z-\bar z'\|_{L_\infty(I;\R^m)}\},
\end{align*}
and hence inductively
\be
\label{Lipall}
|A^k_{x,I,q,N}(t,\bar z;y)-A^k_{x',I,q,N}(t,\bar z';y')|
\le   \max\big\{|y-y'|,  |x-x'|\big\} \sum_{\nu=0}^k\frac{(\|a\||I|)^\nu}{\nu!}.
\ee
Finally, by \eref{Fqgen} and \eref{rhodiff}, one has for any $t,t'\in I $
\begin{align*}
|A_{x,I,q,N}(t,\bar z;y)-A_{x,I,q,N}(t',\bar z;y)|
&\le  \sum_{i=1}^q |\rho_i(t)-\rho_i(t')| 
 |\t A_{N,i}(\bar z;y)| 
\le \|a\| |t-t'|, 
\end{align*}
where we have used \eref{rhodiff}  and  the definition of
$\|a\|$.
This confirms \eref{Lipxy}. The remaining claim follows from 
Remark \ref{rem:A}, \eref{AD}.  \hfill $\Box$


\begin{thebibliography}{9}

\bibitem{Barron}
Andrew R. Barron. Universal approximation bounds for superpositions of a sigmoidal function. IEEE Transactions on Information Theory, 39(3) (1993),930--945. 

\bibitem{BD}
M. Bachmayr, W. Dahmen, Adaptive Low-Rank Approximations for Operator Equations: Accuracy Control and  Computational Complexity, 
Contemporary Mathematics, volume 754,
American Mathematical Society,  https://doi.org/10.1090/conm/754/15151, 
 http://arxiv.org/abs/1910.07052.


\bibitem{BCDD}
A.Barron, A. Cohen, W. Dahmen, R. DeVore,
Approximation and learning
  by greedy algorithms,  
  Annals of Statistics,  3(No 1)(2008), 64--94.

\bibitem{BBGJJ}
C. Beck, S. Becker, P. Grohs, N. Jaafari, A. Jentzen,
Solving the Kolmogorov PDE by means of deep learning,
Journal of Scientific Computing 88 (3), 1--28.

\bibitem{CDDN}
A. Cohen, W. Dahmen, R. DeVore, J. Nichols, Reduced Basis Greedy Selection Using Random Training Sets, 
ESAIM: M2AN, 54 (no 5) (2020), 1509--1524, 
	https://doi.org/10.1051/m2an/2020004, 


\bibitem{CDPW}
A. Cohen, R. DeVore, G. Petrova, P. Wojtaszczyk,  Optimal stable nonlinear approximation,  (2020),  arXiv:2009.09907.

\bibitem{DDGS}
W. Dahmen, R. DeVore, L. Grasedyck, E. S\"uli, Tensor Sparsity  of Solutions to High-Dimensional Elliptic Partial Differential Equations,
Foundation of Computational Mathematics, 16 (No 4) (2016), 813--874.
 (DOI) 10.1007/s10208-015-9265-9.

 

\bibitem{DGM}
W. Dahmen, F. Gruber, O. Mula, An Adaptive Nested Source Term Iteration for Radiative Transfer Equations, 
Mathematics of Computation, 89 (No 324) (2020), 1605--1646,\\
https://doi.org/10.1090/mcom/3505.

\bibitem{DDFHP}
I. Daubechies, R. DeVore, S. Foucart, B. Hanin, G. Petrova, Nonlinear approximation and
(deep) ReLU networks, Constructive Approximation, 55 (2022), 127--172.

\bibitem{CDacta}
A. Cohen and R. DeVore, Approximation of high-dimensional parametric PDEs,
Acta Numer. 24  (2015), 1--159.

\bibitem{DHP}
R. DeVore, B. Hanin, G. Petrova, Neural Network Approximation,
Acta Numerica,  Volume 30 , May 2021,pp. 327 -- 444.
DOI: https://doi.org/10.1017/S0962492921000052

\bibitem{DHM}
R. DeVore, R. Howard, C. Micchelli,  Optimal non-linear approximation,
Manuscripta Math. 4 (1989), 469--478.


\bibitem{E}
W. E, C. Ma, L.  Wu,   The Barron Space and the Flow-Induced Function Spaces for Neural Network Models, 
Constructive Approximation, 
 55 (2022), 369--406.

 \bibitem{GP}
 F. Girosi,  T. Poggio, Representation properties of networks: Kolmogorov's Theorem is irrelevant, iKun neural computation,  (1), (no. 4) (1989),  465--469, 
 doi: 10.1162/neco.1989.1.4.465.


\bibitem{GKNV}
R. Gribonval, G. Kutyniok, M. Nielsen, F. Voigtlaender,  Approximation spaces of deep neural networks, Constructive Approximation,
55(2022),
259--367.  

\bibitem{GH}
P. Grohs, L. Herrmann, Deep neural network approximation for high-dimensional
parabolic Hamilton-Jacobi-Bellman equations,  arXiv:2103.05744v1 [math.NA], March 2021.

\bibitem{GHJW}
P. Grohs, F. Hornung, A. Jentzen, P. Von Wurstemberger,
A proof that artificial neural networks overcome the curse of dimensionality in the numerical approximation of Black-Scholes partial differential equations,
Memoirs of the American Mathematical Society,
https://doi.org/10.48550/arXiv.1809.02362

\bibitem{GPEB}
P. Grohs,  D. Perekrestenko,D. Elbr\"achter, H. B\"olcskei, Deep neural
network approximation theory, IEEE Transactions on Information Theory,
(2019),   arXiv:1901.02220.

 
 \bibitem{GV}
 P. Grohs,  F. Voigtlaender,
 Proof of the Theory-to-Practice Gap in Deep Learning via Sampling Complexity bounds
for Neural Network Approximation Spaces, arXiv:2104.02746v1 [cs.LG], 04/06/2021.


\bibitem{GKP}
I. G\"uhring, G. Kutyniok, P. Petersen,
Error bounds for approximations with deep ReLU neural networks in $W^{s,p}$-norms,
Analysis and Applications, 18 (5)(2020), 803--859,
https://doi.org/10.1142/S0219530519410021,

 \bibitem{HS}
 M. Hansen, C. Schwab, Sparse Adaptive Approximation of High Dimensional
Parametric Initial Value Problems,
 Vietnam J. Math., 41 (2013),181--215,
DOI 10.1007/s10013-013-0011-9


 
\bibitem{Kur}
V. Kurkov\'{a},  Kolmogorov's Theorem is relevant, https://doi.org/10.1162/neco.1991.3.4.617

\bibitem{KPRS}
G. Kutyniok, P. Petersen, M. Raslan, R. Schneider, A Theoretical Analysis of Deep Neural Networks and Parametric PDEs, 
Constructive Approximation, 55 (2022), 73--125.
https://doi.org/10.48550/arXiv.1904.00377


\bibitem{LPT}
J.P. Lions, B. Perthame, E. Tadmor, A kinetic formulation of multidimensional scalar conservation laws
and related equations, Journal of the American Mathematical Society, 7 (no. 1)(1994), 169--191.

\bibitem{Lorentz}
 G.G. Lorentz,    Metric entropy, widths, and superpositions of functions. American Mathematical Monthly. 69 (6)(1962), 469--485. doi:10.1080/00029890.1962.11989915.



\bibitem{Shen}
J. Lu, Z. Shen, H. Yang, S. Zhang,  Deep network approximation for
smooth functions, SIAM J. Math. Anal., 53 (5) (2020),   arXiv:2001.03040.

\bibitem{MasPogg}
H. N. Mhaskar,  T.  Poggio,  Function approximation by deep networks,
Communications on Pure and Applied Analysis,  19(8)(2020), 4085--4095. doi: 10.3934/cpaa.2020181
(2019).  arXiv:1905.12882. 

\bibitem{NW}
E.\ Novak and H.\ Wo\'zniakowski, {\it Approximation of infinitely differentiable
multivariate functions is intractable}, J.\ Complexity {\bf 25} (2009), 398--404.


\bibitem{OpSchwZe}
J.A.A. Opschoor, C. Schwab,  J. Zech,  Exponential ReLU DNN expression of
holomorphic maps in high dimension, Constructive Approximation, 55 (2022),
537--582.


 \bibitem{PP}
 F. Laakmann, P. Petersen, Efficient Approximation of Solutions of Parametric Linear Transport Equations
by ReLU DNNs, Advances in Computational Mathematics, 47 (11) (2021), 
https://doi.org/10.1007/s10444-020-09834-7.

\bibitem{PW}
G. Petrova, P. Wojtaszczyk, Limitations on approximation by deep and shallow neural networks, Dec 2022, arXiv:2212.02223v1 [stat.ML].
 

\bibitem{PV}
P. Petersen,  F. Voigtlaender,  Optimal approximation of piecewise smooth
functions using deep relu neural networks, Neural Networks 108, (2018), 296--330.

\bibitem{Schmid-Hieber}
J. Schmidt-Hieber,  Nonparametric regression using deep neural networks
with ReLu activation, Ann. Statist. 48(4)(2020), 1875--1897.
 DOI: 10.1214/19-AOS1875(2020), 
 

\bibitem{SX}
J.W. Siegel, J. Xu, Approximation rates for neural networks with general activation functions, Neural Netw., 128
(2020), 313--321.  doi: 10.1016/j.neunet.2020.05.019.  




\bibitem{V}
A. Vasseur, Kinetic semi-discretization of scalar conservation laws and convergence by using averaging lemmas,
SIAM J. Numer.Anal., 36 (No.2)(1999), 465--474.

\bibitem{Yar}
D. Yarotsky,   Error bounds for approximations with deep ReLu networks, Neural Networks 94, (2017), 103--114.


\end{thebibliography}
\end{document}